\documentclass{article}
\usepackage{caption}
\usepackage{booktabs}
\usepackage{verbatim}
\usepackage{mathtools}
\usepackage{hyperref}
\usepackage{float}
\usepackage{amssymb}
\usepackage{graphicx}
\usepackage{amsmath}
\usepackage{amsthm}
\usepackage{amsfonts}
\usepackage{tikz}
\usetikzlibrary{positioning}
\usetikzlibrary{shapes.geometric}
\usepackage{scalefnt}
\usepackage[pagewise]{lineno}
\usepackage{appendix}
\def \Rone{\uppercase\expandafter{\romannumeral1}}
\def \Rtwo{\uppercase\expandafter{\romannumeral2}}
\def \Rthree{\uppercase\expandafter{\romannumeral3}}
\def \Rfour{\uppercase\expandafter{\romannumeral4}}
\def \Rfive{\uppercase\expandafter{\romannumeral5}}
\def \Rsix{\uppercase\expandafter{\romannumeral6}}
\numberwithin{figure}{section}
\numberwithin{equation}{section}
\newtheorem{theorem}{Theorem}
\newtheorem{lemma}{Lemma}
\newtheorem{conjecture}{Question}
\newtheorem*{corollary}{Corollary}
\newtheorem{proposition}[lemma]{Proposition}
\newtheorem{main}{Proposition}

\theoremstyle{remark}
\newtheorem*{remark}{Remark}

\theoremstyle{definition}
\newtheorem{definition}{Definition}

\usepackage{fancyhdr}

\fancypagestyle{firstpage}{
  \fancyhf{} 
  \fancyfoot[L]{\textbf{Keywords:} outer billiard, Birkhoff normal form, KAM theory, elliptic island, escaping orbit }
}

\begin{document}

\title{Elliptic islands and zero measure escaping orbit in a class of outer billiards}
\author{Zaicun Li\\ \small Department of Mathematics, University of Maryland, College Park, USA\\ \small zli12324@umd.edu}
\date{}
\maketitle
\thispagestyle{firstpage}

\noindent\textbf{Disclaimer:} This is the version of the article before peer review or editing, as submitted by an author to \textit{Discrete and Continuous Dynamical Systems} (\url{https://www.aimsciences.org/DCDS}). AIMS is not responsible for any errors or omissions in this version of the manuscript or any version derived from it.
\begin{abstract}
We study outer billiard systems around a class of circular sectors.
 For semi-discs, we prove the existence of elliptic islands occupying a positive proportion of the plane. Combined with known results, this shows the coexistence of stability and diffusion for this system.

On the other hand, we show that there exists a countable family of circular sectors
for which the outer billiard system has zero measure of escaping orbits.

\end{abstract}

\section{Introduction and
background}\label{introduction-and-background}
An outer billiard map $F$ is defined on the plane outside a closed convex curve $\Gamma$: it maps each point $z$ to $F(z)$, so that the line passing through these two points is the supporting line from $z$ to $\Gamma$, the contact point divides the line segment $z-F(z)$ in half, and $\Gamma$ lies on the right of the line (see Figure \ref{semi-disc sketch}).

Moser introduced outer billiard systems as simple models to illustrate KAM theory (after Kolmogorov, Arnold and Moser). 
In \cite{MR0442980}, he observed that when $\Gamma$ is smooth and strictly convex, then the map $F$ is approached by a twist map near infinity in an adequate coordinate system. 
He then suggested that KAM theory should be applicable and that there exist smooth quasi-periodic closed invariant curves
for $F$. As a consequence, all orbits would be bounded. Later, Douady gave a complete proof of this phenomenon for strictly convex $C^6$ curves \cite{Douady}. In fact, he showed that there are invariant closed curves that are arbitrarily close to the convex shape and to the infinity.

However,  the smoothness condition is required by KAM theory. Moser asked what the behavior of the outer billiard would be for a non-smooth or a non-strictly convex shape $\Gamma$ such as polygons and the semi-disc. He raised the problem of existence of escaping orbits for outer billiards around such shapes. Phillip Boyland proved in \cite{boyland1994dualbilliardstwistmaps} that there are examples of outer billiard systems of some strictly convex closed curve $G$ such that there are orbits converging to a point in $G$ when the radius of curvature is not $C^1$.

Moser's question was addressed by Dmitry Dolgopyat and Bassam Fayad, in their work \cite{dolgopyat_unbounded_2009} in 2009, where they showed that there is an open ball escaping to infinity for the semi-disc outer billiard.

In this paper, we continue the study of stability phenomena in outer billiards around non-smooth and non-strictly convex shapes. 
We consider a class of shapes that we call circular sectors that generalize the semi-disc by allowing to cut the disc with a line that is not necessarily the diameter (see Figure \ref{circular outer billiard sketch}). We study the dynamics of these outer billiards systems using the first return map to a fundamental domain and using the fact that near infinity the map can be regarded as a perturbation of piecewise linear map in a cylinder. Such piecewise linear map is also known as a "sawtooth map" and has been extensively studied; see \cite{de2012dynamics,CHEN1990217,BIRD1988164,bullett_invariant_1986,QChen_1989,PhysRevA.24.2664,PERCIVAL1987373,Wojtkowski_1982}.

In a first part of the paper, we address the question of coexistence of stability and diffusion for the outer billiard around the semi-disc. 
We show that this is the case, by proving that there are infinitely many elliptic fixed points with stable
neighborhoods (elliptic islands) around them (Theorem \ref{Elliptic islands near fixed points}) (as mentioned above, Dolgopyat and Fayad proved diffusion exists). In fact, we show that the union of the elliptic islands occupies a positive proportion of the plane, which is also the case for the escaping balls that were exhibited in \cite{dolgopyat_unbounded_2009}. The proof relies on KAM theory that can be applied around elliptic periodic points that are found away from the singularity lines of the outer billiard map. 

In a second part of the paper, we go in the opposite direction of \cite{dolgopyat_unbounded_2009}, and show that there exists a countable family of circular sectors for which the outer billiard has zero measure of escaping orbits (Theorem \ref{main no escape}). We also showed that the method can be applied to small perturbations of integrable systems such as the Fermi-Ulam model \cite{de2012dynamics}.

To finish this introduction, we mention that in \cite{gutkin_dual_1992,vivaldi_global_1987}, it was proved that for the so-called quasi-rational polygons all trajectories are bounded. However, Schwartz proved the existence of an unbounded orbit for kite shapes with vertices $\left(-1,0\right)$,  $\left(0,1\right)$,  $\left(0,-1\right)$ and  $\left(A,0\right)$ for irrational $A$ in \cite{schwartz2007unboundedorbitsouterbilliards,schwartz2008outerbilliardskites}.

\bigskip 

\noindent {\bf Plan of the paper.} The two main results, Theorem \ref{Elliptic islands near fixed points} and \ref{main no escape} are stated in Section \ref{main-results}. 

In Section \ref{Proof of elliptic island}, we give the complete proof of existence of elliptic islands (Theorem \ref{Elliptic islands near fixed points}). We first prove Theorem \ref{Elliptic islands near fixed points} using KAM theory and the Birkhoff normal form near a well-chosen sequence of period orbits of $F$, which is given by Proposition \ref{normal form near fixed point}. The rest of section \ref{Proof of elliptic island} is devoted to the proof of Proposition \ref{normal form near fixed point}.

In Section \ref{Proof no escape}, we prove that there is zero measure of escaping orbits in a family of circular sector outer billiard systems (Theorem \ref{main no escape}). We first show that the linear part of the first return map to a fundamental domain has no escaping orbits. The first return map is given by Lemma \ref{keylemma}. The process of finding this first return map is technical and is put in the appendix (see Appendix \ref{proof of key lemma}).

Finally, in section \ref{Open questions}, some open questions are discussed.

\section{Notations}
In this paper, for the norms of the $C^5$ mappings, we use the following notation.
\begin{definition}[$C^5$ norm for mappings]
   A $C^5$ mapping $f: \mathbb{R}^2 \supset U \rightarrow \mathbb{C}^2$ has the following norm:
   
   $$\left\|f\right\|=\sup_{(x,y)\in U,m+n\leq 5}\left \vert\frac{\partial f}{\partial^mx\partial^ny}\right\vert$$
\end{definition}
For two mappings (real or complex) defined on the same domain, we write 
$$f = \mathcal{O}(g),$$
if there is a positive constant $M$ such that 
$$\| f \| \leq {M}{\|g \|}.$$

If for a family of maps $f_n: \mathbb{R}^2 \supset U \rightarrow \mathbb{R}^2$ if there is a positive constant $M$ such that
$$\| f_n \| \leq \frac{M}{n^k},$$ then we write 
$$f_n = \mathcal{O}(n^{-k}) \ $$
If moreover, $f_n$ satisfy the condition that for any $r+s \leq l$,
{\color{red}
 $$\frac{\partial f_n}{\partial^rx\partial^sy}(0,0)=0$$
 then we use the notation
 $$f_n=\mathcal{O}_l(n^{-k}).$$
}
\section{Main results}\label{main-results}
\subsection{Outer billiard for the semi-disk}
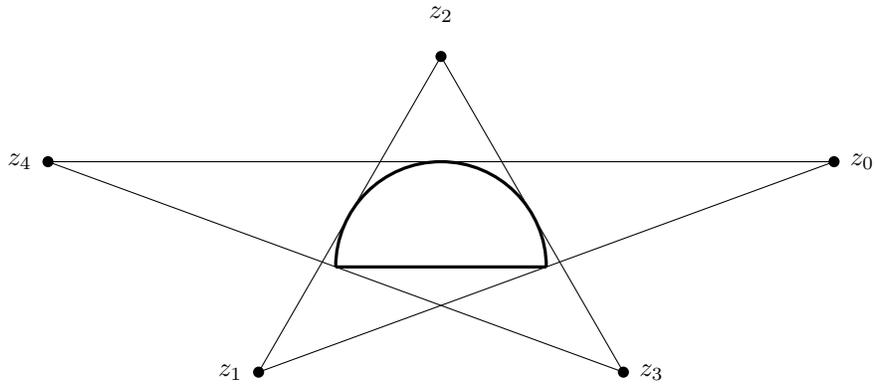
\begin{figure}[htbp]
	\centering
	\begin{tikzpicture}[scale=1.4]
		\draw[very thick] (1,0) arc (0:180:1);
		\draw [very thick] (1,0)--(-1,0);
		\draw  ({2+sqrt(3)},1)--({-sqrt(3)},-1);
		\draw  (-{sqrt(3)},-1)--(0,2);
		\draw  (0,2)--({sqrt(3)},-1);
		\draw  ({-2-sqrt(3)},1)--({sqrt(3)},-1);
		\draw ({-2-sqrt(3)},1)--({2+sqrt(3)},1);
		\fill ({2+sqrt(3)},1) circle(1.5pt);
		\fill ({-2-sqrt(3)},1) circle(1.5pt);
		\fill (0,2) circle(1.5pt);
		\fill ({sqrt(3)},-1) circle(1.5pt);
		\fill ({-sqrt(3)},-1) circle(1.5pt);	
		\node  at (4,1) { $z_0$};
		\node  at (-2,-1) { $z_1$};
		\node  at (0,2.4) { $z_2$};
		\node  at (2,-1) { $z_3$};
		\node  at (-4,1) { $z_4$};
	\end{tikzpicture}
	\caption{A $5$-periodic point for semidisc outer billiard}
	\label{semi-disc sketch}
\end{figure}

Let us first consider the outer billiard map $F$ for a semi-disc. Figure \ref{semi-disc sketch} gives an example of a periodic point of period $5$ for the outer billiard map. 

\begin{figure}[H] 
	\centering 
	\includegraphics[width=0.7\textwidth]{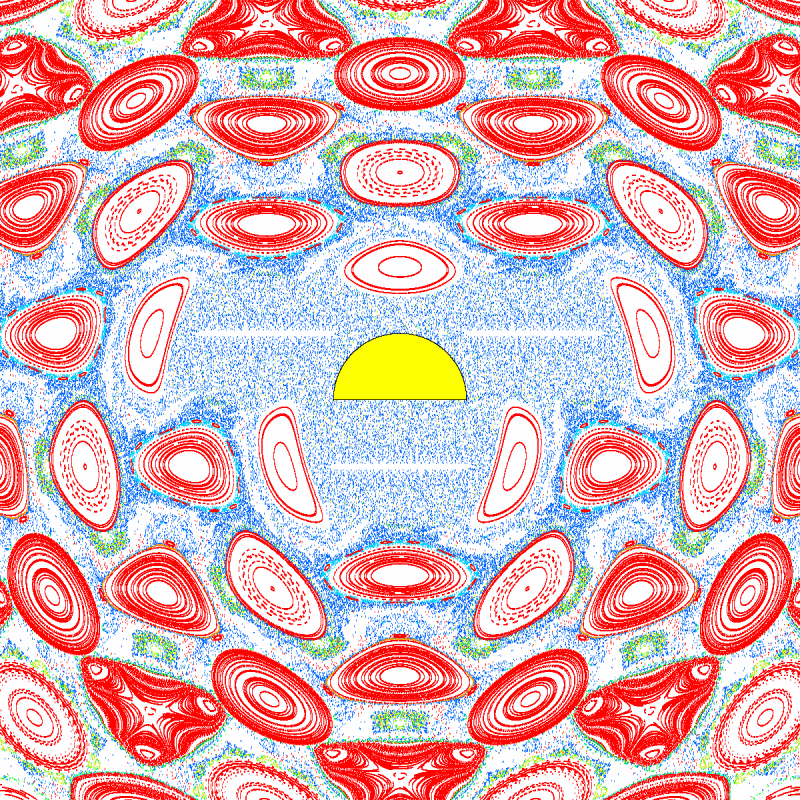} 
	\caption{Elliptic islands for the semi-disc outer billiard map from \cite{mathematikExteriorBilliards}} 
	\label{elliptic islands} 
\end{figure}

As shown in Figure \ref{elliptic islands}, there are periodic orbits with stability regions around them. Our goal is to show that there are infinitely many such stable periodic orbits. Let us start with some precise definitions.

We will prove the existence of elliptic islands around elliptic periodic points of the outer billiard map $F$, which we will first define before stating the exact result. We consider a map $T$ from an annulus {\color{red}$ \mathbb{R}^2  \supset \textbf{A} $} to $\mathbb{R}^2$.

\begin{definition}[Self-intersection property]
	The map $T$ is said to have a self-intersection property if for any closed curve $\gamma \in \textbf{A}$ , $\gamma \cap T(\gamma) \neq \emptyset$. 
\end{definition}

\begin{definition}[KAM tori]
KAM tori refer to invariant simple closed curves of the map $T$.
\end{definition}
In classical KAM theory, if $T$ is a perturbation of a non-degenerate integrable map of $\textbf{A}$ that has the self-intersection property and is $\mathcal{O}(\epsilon)$ close to the original map in $C^5$ norm, then KAM tori will fill up $1-\mathcal{O}(\sqrt{\epsilon})$ proportion of Lebesgue measure of the annulus. In fact, an invariant circle on which the original map has the Diophantine rotation number will be preserved after perturbation, and each KAM torus will be $\mathcal{O}(\epsilon)$ close to the original circle with the same Diophantine rotation number in the $C^1$ norm.

To be precise, consider a map of $\mathcal{F}:\textbf{A} \rightarrow   \mathbb{R}^2$ in polar coordinates:
	\begin{align*}
	 \mathcal{F}&: \left(r,\theta\right)\mapsto \left(r',\theta'\right)\\
	\theta'&=\theta+\alpha+\alpha_2r+f,\\
	r' &= r +g,
    \end{align*}	
where $\alpha_2$ is nonzero. From the classical KAM theory, see for example the proof on page 50 in \cite{Moser+2001}, if there exists an $\epsilon$ positive and small enough, such that 
{\color{red}$$\|f\|,\|g\| \leq |\alpha_2|\epsilon$$} 
and the mapping $\mathcal{F}$ on $\textbf{A}$ has the self-intersection property, then $\mathcal{F}$
possesses an invariant curve of the form 
$$r=c+u(t),\ \ \theta=t+v(t)$$
in $\textbf{A}$ where $u,v$ are continuously differentiable, of period $2\pi$ and satisfy
$$\left\|u\right\|_1+\left\|v\right\|_1<\mathcal{O}\left(\epsilon\right).$$
Moreover, the induced mapping on this curve is given by:
$$t \rightarrow t+\omega,$$
where $\omega$ is a Diophantine number. Such invariant curves will cover over $1-\mathcal{O}\left(\sqrt{\epsilon}\right)$ proportion of the annulus $\textbf{A}$.
\begin{definition}[$(\epsilon,\delta)$-Elliptic Island]
For a map $\mathcal{F}$ defined in some open set in $\mathbb{R}^2$, if $z$ is an elliptic $q$-periodic point of $\mathcal{F}$, then an elliptic island $I$ is a neighborhood of $z$ fixed by $\mathcal{F}^q$ with a set of KAM tori inside. In this paper, such elliptic islands in the plane are called  $(\epsilon,\delta)$-elliptic Islands if they also contain the disc of radius $\delta$ centered at the fixed point and KAM tori fill up to $(1-\epsilon)$ proportion of the Lebesgue measure of the elliptic island.
\end{definition}

Recall that $F$ denotes the outer billiard map around the semi-disc, then we have the following theorem.
\begin{theorem} [Elliptic islands near periodic points]\label{Elliptic islands near fixed points}
There is a sequence of elliptic periodic points $\left(x_n,y_n\right)$ for $F^2$, where
\[x_{n} = 3n - \frac{3}{4}+\mathcal{O}\left( \frac{1}{n} \right),\]
\[y_{n} = 1 + \mathcal{O}\left( \frac{1}{n} \right).\]
There is a fixed $\delta>0$, such that for any $\epsilon>0$, there exists an integer $N(\epsilon)$ such that for any integer $n \geq N(\epsilon)$, there exists an $(\epsilon,\delta)$-elliptic island for $F$ near $\left( x_{n},y_{n} \right)$.

\end{theorem}

\begin{remark}
	The diameter of the elliptic islands near the periodic points $\left( x_{n},y_{n} \right)$ are not decreasing when $n$ is growing.  Moreover, when $n$ is larger, $\epsilon$ converges to $0$ and the elliptical islands are filled more with KAM tori. 
	
	Let us consider the right $x$-axis and its image under map $F^2$, which is a strip near $x$ -axis with width less than $2$. This is a {\color{red}fundamental domain} for $F^2$. Observe that $\left( x_{n},y_{n} \right)$ lies in this region. Since $\delta$ is uniform for all elliptic islands introduced in Theorem \ref{Elliptic islands near fixed points}, 
	they will occupy a positive proportion of this Poincar\'{e} section. Under the map $F^2$ they will rotate around the origin and take up a positive proportion of the whole plane. This is depicted in Figure \ref{elliptic islands}.
\end{remark}
\subsection{Outer billiard for circular sector}
\begin{figure}[htbp]
	\centering
	\begin{tikzpicture}[scale=1.2]
		\draw[very thick] (cos{-67.5},sin{-67.5})--(cos{247.5},sin{247.5});
		\draw[very thick] (cos{-67.5},sin{-67.5}) arc (-67.5:247.5:1);
		\draw [dashed] (0,0)--(cos{-67.5},sin{-67.5});
		\draw [dashed] (0,0)--(cos{247.5},sin{247.5});
        \draw (0,{sqrt(5)})--({-0.8*sqrt(5)},{-0.6*sqrt(5)});
        \draw ({-0.8*sqrt(5)},{-0.6*sqrt(5)})--({2*cos(-67.5)+0.8*sqrt(5)},{2*sin(-67.5)+0.6*sqrt(5)});
        \node at ({2*cos(-67.5)+0.8*sqrt(5)+0.3},{2*sin(-67.5)+0.6*sqrt(5)}) {$z_0$};
        \node at ({-0.8*sqrt(5)-0.2},{-0.6*sqrt(5)-0.2})  {$z_1$};
        \node at (0,2.5) {$z_2$};
		\node  at (0,-0.7) { $2\beta$};
      	\fill (0,{sqrt(5)}) circle(1.5pt);
      	\fill ({-0.8*sqrt(5)},{-0.6*sqrt(5)}) circle(1.5pt);
      	\fill ({2*cos(-67.5)+0.8*sqrt(5)},{2*sin(-67.5)+0.6*sqrt(5)}) circle(1.5pt);
	\end{tikzpicture}
	\caption{ Outer billiard map for a circular sector with angle $2\pi-2\beta$}
	\label{circular outer billiard sketch}
\end{figure}
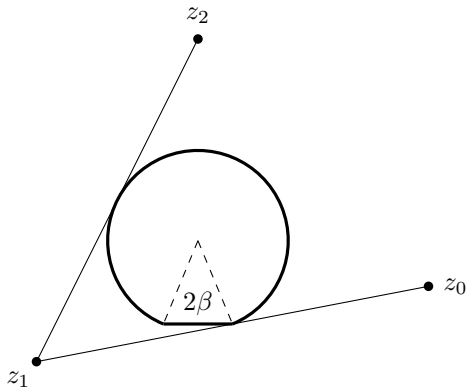

In the paper by Dolgopyat and Bassam \cite{dolgopyat_unbounded_2009}, they proved that for the semi-disc outer billiard map, there exists an escaping open ball. This result can be generalized to the outer billiard maps of circular sectors that are close to a semi-disc. In this section, we consider the outer billiard map $F_\beta$ for the circular sector with angle $2\pi-2\beta$, as shown in figure \ref{circular outer billiard sketch}. Our result goes in the opposite direction and finds a sequence of $\beta$ with zero-measure escaping orbits.

We define the escaping orbit set  \(E_\beta\) = \{\(x\): $\lim_{n \to \infty}{F}_\beta^{n}(x)=\infty$\} and let $\mu$ be the Lebesgue measure on the plane. In this paper, we show that if the value of $\beta$ is chosen adequately, the set of escaping orbits is of zero Lebesgue measure on the plane. 
\begin{theorem}\label{main no escape}
	There is a decreasing sequence of $\beta_m \in (0,\pi/2)$, with
    $\lim_{m \to \infty}\beta_m=0$, such that the Lebesgue measure $\mu(E_{\beta_m})=0$,.
\end{theorem}
\begin{remark}
	The values of such $\beta$ are chosen in order that $\pi$ is an integral multiple of the rotation number of the asymptotic behavior of the mapping $\mathcal{F}_{\beta}^2$ when the orbit goes to infinity. Thus, it is a rare event to have zero Lebesgue measure of escaping orbits.
\end{remark}

\section{Proof of theorem \ref{Elliptic islands near fixed points}}\label{Proof of elliptic island}
\begin{figure}[H]
	\centering
	\includegraphics[width=4in,height=3.2in]{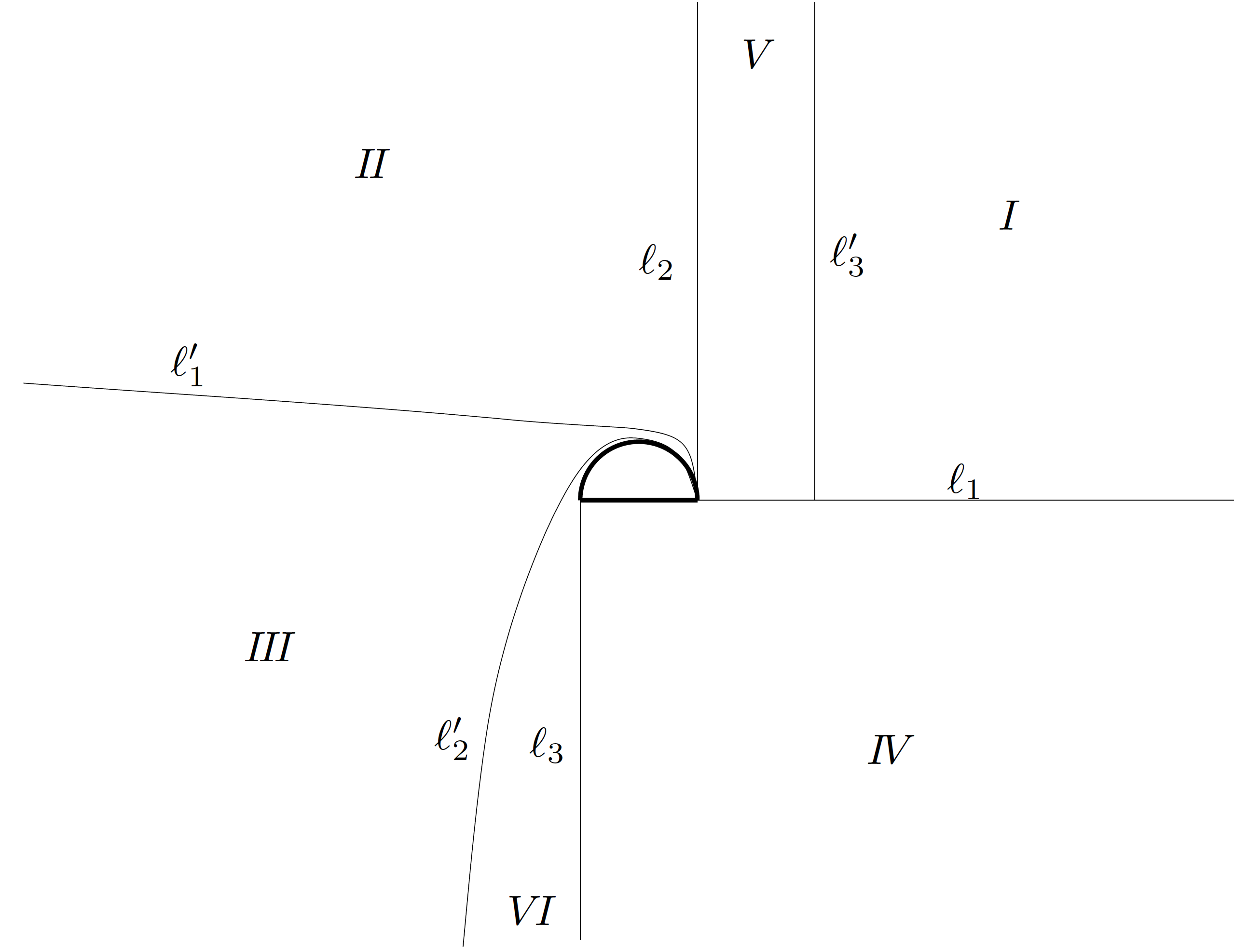}
	\caption{Continuous region  for $F^2$ of a semi-disc outer billiard from \cite{dolgopyat_unbounded_2009}}
	\label{semi-disc}
\end{figure}

Following \cite{dolgopyat_unbounded_2009}, we divide the plane into six
regions such that in each region the two iterations of the outer billiard map $F$ are smooth, as shown in Figure \ref{semi-disc}. Recall that \(F\)
denotes the outer billiard map. Leaving \(x_{0}\) as a large
constants, we consider the infinite region
\(\mathcal{D}\) bounded by \(\ell_{1}\), \(F^{2}\ell_{1}\),
and \(\{ x = x_{0}\}\).  Let \(\mathcal{F}\) denote the first return map to the region
\(\mathcal{D}\).
We need to verify that the first terms of the Birkhoff normal forms of the sequence of elliptic points are nonvanishing and use the KAM twist theorem to prove theorem \ref{Elliptic islands near fixed points}. The following theorem gives the conjugacy we need.

In this section, we use the notation $P(x,y,z,w)$ to denote the parallelogram in the real plane or the complex plane with four vertices $a,b,c,d \in \mathbb{R}^2$ or $a,b,c,d \in \mathbb{C}$. 
We define 
$${\mathcal{E}}=P\left(\left(\frac{1}{256},0\right), \left(-\frac{1}{256},\frac{1}{64}\right), \left(-\frac{1}{256},0\right), \left(\frac{1}{256},-\frac{1}{64}\right)\right)$$
to be a parallelogram in the real plane, and
\begin{align}
	\mathbf{R}&=P\left(\left(0,\frac{3}{512}\right),  \left(\frac{3}{256},-\frac{3}{512}\right), \left(0,-\frac{3}{512}\right),  \left(-\frac{3}{256},\frac{3}{512}\right)\right),\\
	\mathbf{R}'&=P\left(\left(\frac{3}{512},\frac{3}{1024}\right),  \left(\frac{3}{512},-\frac{3}{1024}\right), \left(-\frac{3}{512},-\frac{3}{1024}\right),  \left(-\frac{3}{512},\frac{3}{1024}\right)\right)
\end{align}

to be two rectangles centered at the origin.

\begin{main} [Normal forms]\label{normal form near fixed point}
	There exists a sequence of elliptic periodic points of $F$: 
	$$\left\{\left(x_n,y_n\right)=\left(3n-\frac{3}{4}+\mathcal{O}\left(\frac{1}{n}\right),1+\mathcal{O}\left(\frac{1}{n}\right)\right)\right\}_{n\geq N}$$
for some large $N$. 
	Let $\mathbf{R}_{n}$ be a neighborhood of the periodic point $(x_n,y_n)$, where 
	\[ \mathbf{R}_{n} = \left( 3n - \frac{3}{4} ,1  \right)+\mathbf{R}\]	is obtained by translating $R$ to be centered at $\left(3n-3/4,1\right)$.
	We have a conjugacy $C$, such that
	\begin{align}\label{elllptic formula}
		C&:\ \mathbf{R}_n \rightarrow \mathbb{R}^2,\ \  {\mathcal{E}} \subset C\left(\mathbf{R}_n\right)\notag \\
		&\ \  \left(x,y\right) \mapsto \left(\xi,\eta\right)\notag\\
		\begin{split}
		\xi&=\frac{y}{4}-\frac{1}{4}+\mathcal{O}(n^{-1}),\\
		\eta&=\frac{\sqrt{2}x}{4}-\frac{3\sqrt{2}n}{4}+\frac{3\sqrt{2}}{16}+\mathcal{O}(n^{-1}).
		\end{split}
	\end{align}
	Moreover, \(C\circ \mathcal{F} \circ C^{-1}: {\mathcal{E}} \rightarrow \mathbb{R}^2\) takes the following form
	\begin{align}\label{rotation formula}
	 C\circ \mathcal{F} \circ C^{-1}&:\left(\xi,\eta\right) \mapsto \left(\xi',\eta'\right)\notag\\
	 \begin{split}
	 		\xi' &= \xi\cos w - \eta\sin w + \mathcal{O}_{4}(n^{-3}),\\
	 	\eta' &= \xi\sin w + \eta\cos w + \mathcal{O}_{4}(n^{-3}),\\
	 	w &= \alpha(n) + \alpha_{2}(n)\left( \xi^{2} + \eta^{2} \right),
	 \end{split}
	\end{align}
	where \((\xi',\eta')\) is the image of
	\((\xi,\eta)\), and where
	\[\alpha(n) = \arccos\left( - \frac{7}{9} \right)+\mathcal{O}\left( \frac{1}{n} \right),\]
	\[|\alpha_{2}(n)| \geq \frac{c}{n^{2}},\]
	for some constant \(c > 0\) independent of \(n\).
\end{main}

\begin{remark}
	The proof of {\color{red}Proposition \ref{normal form near fixed point}} will be given in the next section. {\color{red}We only need the perturbation term in \ref{rotation formula} to be $\mathcal{O}\left(n^{-3}\right)$ instead of $\mathcal{O}_{4}(n^{-3})$ for the proof of Theorem \ref{Elliptic islands near fixed points} (the proof does not require the first three derivatives at the origin to vanish; the theorem merely states this fact).}
\end{remark}


\begin{proof}[Proof of theorem \ref{Elliptic islands near fixed points}]
	    
	     Let $\textbf{A}=\left(\epsilon^4,\frac{1}{2^{18}}\right) \times \mathbb{T}$ be an annulus in the plane centered at the origin and $\textbf{B}$ be a ball with radius $\frac{1}{512}$ in the plane centered at the origin.
        We conjugate the normal form in Proposition \ref{normal form near fixed point} using the following polar coordinates:
     	\begin{align*}
     			U: \textbf{A} & \rightarrow \mathcal{E}\\
     			(r,\theta)& \mapsto (\xi,\eta)=\left( \sqrt{r}\cos\theta,\sqrt{r}\sin\theta \right)
     	\end{align*}
     Then we have 	$\mathcal{F}' =U^{-1}\circ C\circ  \mathcal{F} \circ C^{-1} \circ U:\textbf{A} \rightarrow \mathbb{R}^2$ such that
     	\begin{align*}
     		\mathcal{F}' &: (\theta,r) \rightarrow (\theta',r')\\
     		\theta'&=\theta+\alpha(n)+\alpha_2(n)r+f,\\
     		r' &= r +g,
     	\end{align*}	
     where 
     $$f,g=\mathcal{O}(n^{-3}).$$
     	Notice that  $U(\textbf{A})\subset \textbf{B}\subset\mathcal{E}$.
     	
        Since $\left|\alpha_{2}(n)\right| \geq \frac{c}{n^{2}}$,we can choose $n(\epsilon)$ large enough, such that we have $\|f\|,\|g\| \leq |\alpha_2(n)|\epsilon^4$. 
        Moreover, the outer billiard map for a half disc is area preserving on the real plane, hence it has the self-intersection property. This property is preserved under any conjugacy, thus the return map $\mathcal{F}$ and the map $\mathcal{F}'$ also have the self-intersection property.

        On the other hand, from classical KAM theory, see, for example, page 50 in \cite{Moser+2001}. There exists a $\epsilon'$, such that for the mapping $\mathcal{F}'$ in KAM tori will take up the proportion $1-\mathcal{O}\left(\sqrt{\epsilon^4}\right)=1-\mathcal{O}\left(\epsilon^2\right)$ of the annulus $\textbf{A}$, and thus will take up the proportion $1-\mathcal{O}\left(\epsilon^2\right)$ of $\mathbf{B}$.
       
       Considering the fact that the conjugacy $C$ is constant modulo $\mathcal{O}(n^{-1})$,
       we can choose a positive $\delta$, such that $C^{-1}\left(\textbf{B}\right) \subset \mathbf{R}_n$ contains the ball centered at $\left(x_n,y_n\right)$ with radius $\delta$ for any $n$.
       
        The only problem now is that $C^{-1}\left(\textbf{B}\right)$ is not invariant under the map $\mathcal{F}$. To fix this problem, we can look at the area inside a KAM circle in $\textbf{B}$. To be more specific, we choose a KAM torus $\gamma$ that is $\mathcal{O}(\epsilon^2)$ to the outer boundary of $\textbf{B}$. Denote the part of $\textbf{B}$ inside $U(\gamma)$ by $\textbf{B}'$. Then KAM tori will take up over $1-\mathcal{O}\left(\epsilon^2\right)$ proportion of the annulus $\mathbf{B}'$. The number $\epsilon$ is chosen small enough so that $C^{-1}\left(\textbf{B}'\right) \subset \mathbf{R}_n$ still contains the ball centered at $\left(x_n,y_n\right)$ with radius $\delta$ for any $n$.
     	
     	Then $C^{-1}(\textbf{B}')$ is the elliptic island that is filled by KAM tori up to $1-\mathcal{O}\left(\epsilon^2\right)$ and contains a $\delta$ ball inside. If we choose $\epsilon$ small, it would be an $\left(\epsilon,\delta\right)$-elliptic island near $\left(x_n,y_n\right)$, and we obtain theorem \ref{Elliptic islands near fixed points}.
\end{proof}

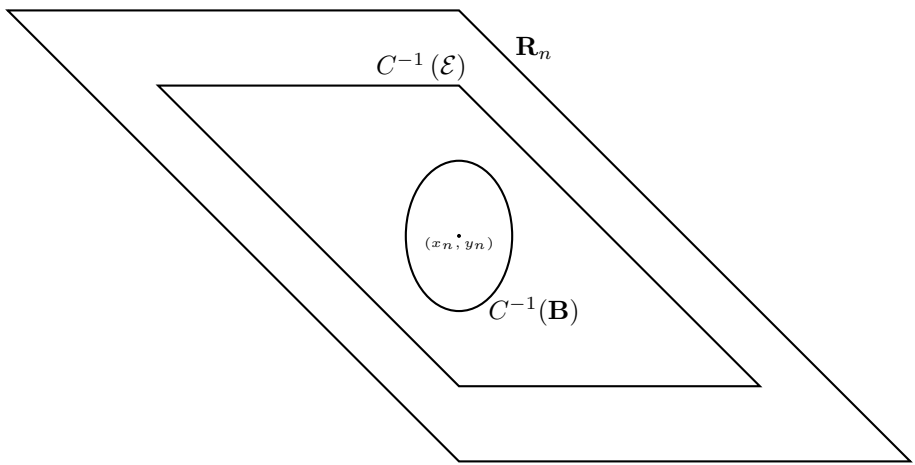
\begin{figure}[H]
	\centering
	\begin{tikzpicture}[scale=0.5]
     \coordinate (A) at (0,6);
  \coordinate (B) at (12,-6);
  \coordinate (C) at (0,-6);
  \coordinate (D) at (-12,6);

  \draw[thick] (A) -- (B) -- (C) -- (D) -- cycle;
  \coordinate (A1) at (0,4);
  \coordinate (B1) at (8,-4);
  \coordinate (C1) at (0,-4);
  \coordinate (D1) at (-8,4);

  \draw[thick] (A1) -- (B1) -- (C1) -- (D1) -- cycle;
          \draw [thick] (0,0) ellipse ({sqrt(2)} and 2 );
          \fill (0,0) circle(1.5pt);
          \node at (-0,-0.2) {\tiny  $(x_n,y_n)$};
          \node at (-1,4.5) { $C^{-1}\left(\mathcal{E}\right)$};
          \node at (2,5) {$\mathbf{R}_n$};
          \node at (2,-2) {$C^{-1}(\textbf{B}) $};
	\end{tikzpicture}
	\caption{Sketch of $\mathbf{R}_n, \mathcal{E}$ and elliptic islands}
	\label{domians figure}
\end{figure}

The rest of this section is devoted to the proof of {\color{red}Proposition \ref{normal form near fixed point}}. The proof will follow from the composition of two conjugacies. 
In subsection \ref{find H1} the first conjugacy $H_1$ is defined near the sequence of elliptic periodic points which we found to simplify the dynamics and translate each fixed point to the origin in $\mathbb{R}^2$. In subsection \ref{find H2} the second conjugacy $H_2$ is defined to give the first terms of the Birkhoff normal form so that we can check that the twist term is not zero.

\subsection{Find the conjugacy \texorpdfstring{$H_1$}{TEXT} near the fixed points of  \texorpdfstring{$\mathcal{F}$}{TEXT}}\label{find H1}
Recall that we define $R$ and $R'$ are parallelogram centered at the origin.
\begin{main}\label{twisted form near fixed points}
	There is a sequence of fixed points for the map $\mathcal{F}$ in the region $\mathcal{D}$, $$\left\{\left(x_n,y_n\right)=\left(3n-\frac{3}{4}+\mathcal{O}\left(\frac{1}{n}\right),1+\mathcal{O}\left(\frac{1}{n}\right)\right)\right\}_{n\geq N}$$
for some large $N$. Around each fixed point there is a neighborhood 
	$$\mathbf{R}_{n} = \left( 3n - \frac{3}{4} ,1  \right)+\mathbf{R}$$
	and there is a neighborhood of the origin, 
	\[\mathbf{R}_{n}'=\mathbf{R}'+\mathcal{O}\left(\frac{1}{n}\right)\]
	There is a conjugacy
	\begin{align*}
		H_1: \ \mathbf{R}_n &\rightarrow \mathbf{R}_n'\\
		(x,y) & \mapsto (x_1,y_1)
	\end{align*}
	such that	
    \begin{equation}\label{xytoxy}
    \begin{split}
        x_1&=x+y-3n-\frac{1}{4}+\mathcal{O}\left(\frac{1}{n}\right),\\
        y_1&=\frac{1}{2}y-\frac{1}{2}+\mathcal{O}\left(\frac{1}{n}\right),
    \end{split}
\end{equation}
and
	\begin{equation}\label{pre_Birk}
		\begin{split}
			H_1\circ \mathcal{F}\circ H_1^{-1}:&\left( x_{1},y_{1} \right)\mapsto ({\widehat{x}}_{1},{\widehat{y}}_{1})\\
			{\widehat{x}}_{1} = A_{1,1}x_{1} +& A_{1,2}y_{1} + \sum_{j = 2}^{3}{F_{1,j}\left( x_{1},y_{1} \right)}+\mathcal{O}_4\left( n^{- 3} \right),\\
			{\widehat{y}}_{1} = A_{2,1}x_{1} + &A_{2,2}y_{1} + \sum_{j = 2}^{3}{F_{2,j}\left( x_{1},y_{1} \right)}+  \mathcal{O}_4\left( n^{- 3} \right),\end{split}
	\end{equation}
where ${\rm det}(A_{ij})=1$, and
\[A_{1,1} = \frac{1}{9}+\mathcal{ O}\left( n^{- 1} \right),\ \ A_{1,2} = - \frac{8}{3}+\mathcal{O}\left( n^{- 1} \right),\]
\[A_{1,1} = \frac{4}{9}+\mathcal{ O}\left( n^{- 1} \right),\ \ A_{2,2} = - \frac{5}{3}+\mathcal{O}\left( n^{- 1} \right),\]
and $F_{i,j}$ is a homogeneous polynomial of degree $j$ for any $i,j$ and  $F_{i,j}=\mathcal{O}_{j}(n^{-j+1})$.
More specifically, we have:
\begin{equation}\label{Fij}
\begin{split}
      F_{i,2}&=0+\mathcal{O}\left(n^{-2}\right),\ \ i=1,2,\\
      F_{1,3}&=\frac{1}{n^2}\left(-\frac{1520xy^2 }{729}+\frac{2080{x}^2 {y}}{2187}-\frac{3392{x}^3 }{19683}+\frac{1360{y}^3 }{729}\right)+\mathcal{O}\left(n^{-3}\right),\\
      F_{2,3}&=\frac{1}{n^2}\left(-\frac{32{x}{y}^2 }{81}+\frac{64{x}^2 {y}}{243}-\frac{128{x}^3 }{2187}+\frac{16{y}^3 }{81}\right)+\mathcal{O}\left(n^{-3}\right).
      \end{split}
\end{equation}

\end{main}
\begin{remark}
	The explicit form of the polynomials $F_{i,j},j=2,3$ we list here will be needed in the computation of the twist term in the Birkhoff normal form, and we computed them in Appendix \ref{calculation}.
\end{remark}

We prove this proposition in two steps. Firstly, we study the dynamics in the regions $\Rone\sim\Rfour$ to find coordinates that simplify $F^2$. Then we use these coordinates in the regions $\Rone\sim\Rfour$ and change back to the polar coordinates to transition from one region to the next to calculate the first return map. The detailed proof is technical and is put in Appendix \ref{proof of proposition B}.

\subsection{Find the conjugacy \texorpdfstring{$H_2$}{TEXT}, the beginning of the Birkhoff normal form and the twist property}\label{find H2}
\begin{main}\label{FindBirk}
There is a conjugacy \(H_{2}\) giving the quadratic part
of the Birkhoff normal
form of \emph{(}\ref{pre_Birk}\emph{)},

$$H_2:\ \mathbf{R}_n' \rightarrow \mathbb{R}^2$$
where  $H_2\left(R_n'\right)\supset {\mathcal{E}}$ 
and on ${\mathcal{E}}$, the map \(H_2\circ H_1\) conjugates \(\mathcal{F}\)
to
\begin{equation*}
	\begin{split}
	    H_2\circ H_1\circ &\mathcal{F}\circ H_1^{-1}\circ H_2^{-1}:\left( \xi,\eta\right) \mapsto \left( \xi',\eta' \right)\\
		\xi' &= \xi\cos w - \eta\sin w + \mathcal{O}_{4}\left({n^{-2}}\right),\\
		\eta' &= \xi\sin w + \eta\cos w +\mathcal{O}_{4}\left({n^{-2}}\right),\\
		w &= \alpha(n) + \alpha_{2}(n)\left( \xi^{2} + \eta^{2} \right),
	\end{split}
\end{equation*}
where
\[\alpha(n) = \arccos\left( - \frac{7}{9} \right)\mathcal{+ O}\left( \frac{1}{n} \right),\]
\[|\alpha_{2}(n)| \geq \frac{c}{n^{2}},\]
for some constant \(c > 0\) independent of \(n\).
\end{main}

To prove proposition \ref{FindBirk}, we need three conjugacies, which are given by the following lemmas. We follow Moser's proof introduced in his book \cite{moser_invariant_1962} from page 155 to page 174.

The first step is to diagonalize the linear part of the map $H_1 \circ \mathcal{F} \circ H_1^{-1}$. To achieve this goal, we introduce a complex linear conjugacy.
\begin{lemma}
	There is a linear conjugacy $H_{2,1}:\ (x_1,y_1) \rightarrow (a,b)$, diagonalizing the linear part of the map $H_1\circ \mathcal{F}\circ H_1^{-1}$. Let 
	\begin{equation}\label{domain of a}
		U=0.99P\left(	\frac{3}{2048},  -\frac{3}{2048}+\frac{3\sqrt{2}i}{1024}, -\frac{3}{1024},  \frac{3}{2048}-\frac{3\sqrt{2}i}{1042}\right)
	\end{equation}
    be a complex parallel, and let $\widetilde{U}=\left\{(a,b)\in \mathbb{C}^2, a=\overline{b}, a\in U\right\}$.
    Then for $n$ large enough,
	 $$H_{2,1}\left(R_n'\right)\supset \widetilde{U}.$$ 
	On $\widetilde{U}$ we have ,
	\begin{align}\label{formab}
		 T_1=H_{2,1}&\circ H_1\circ \mathcal{F} \circ H_1^{-1}\circ H_{2,1}^{-1}: (a,b)\mapsto (a',b'),\notag\\
		a'&=\lambda a + \sum_{j = 2}^{3}{G_{1,j}(a,b)}+\mathcal{O}_4\left( n^{- 3} \right),\notag\\
		b'&=\lambda^{- 1}b + \sum_{j = 2}^{3}{G_{2,j}(a,b)}+\mathcal{ O}_4\left( n^{- 3} \right),
	\end{align}
	where \(\lambda = exp(i(\alpha + \mathcal{O}\left( n^{- 1} \right)))\),
	\(\cos(\alpha) = -7/9\). Moreover,
	$$G_{1,j}(a,b) = \sum_{k = 0}^{j }{G_{1,j}^{k}a^{k}b^{j - k}},$$
	$$G_{2,j}(a,b) = \sum_{k = 0}^{j }{G_{2,j}^{k}b^{k}a^{j  - k}}$$
	then we get \(G_{1,j}^{k} = \overline{G_{2,j}^{k}}\). \(G_{1,j}\) and
	\(G_{2,j}\) are of order \(\mathcal{O}\left( n^{- j+1} \right)\).
	More specifically, we have: 
	\begin{equation}\label{Gij}
    \begin{split}
		G_{1,2}^{0}&,\ G_{1,2}^{0},\ G_{1,2}^{0}=0+\mathcal{O}\left(n^{-2}\right),\\
		G_{1,3}^2&=\frac{1}{n^{-2}}\left(-\frac{32}{81}+\frac{28\,\sqrt{2}\,\mathrm{i}}{81}\right)+\mathcal{O}\left(n^{-3}\right).
        \end{split}
	\end{equation}
	
\end{lemma}

\begin{proof}
To diagonalize the linear part, we conjugate this
system by a linear transform $H_{2,1}$:

\begin{equation}\label{Elp}
\begin{bmatrix}
x_{1} \\
y_{1}
\end{bmatrix} = \left( \begin{bmatrix}
2-\sqrt{2}i& 2 + \sqrt{2}i \\
1 & 1
\end{bmatrix}\mathcal{+ O}\left( n^{- 1} \right) \right)\begin{bmatrix}
a \\
b
\end{bmatrix},
\end{equation}
where $a$ and $b$ lie in the complex plane and $a=\bar{b}$. Straightforward calculation from $R_n'$ will show that when $(x_1,y_1)$ is real the domain of $a$ up to $\mathcal{O}(n^{-1})$ accuracy, is 
$$P\left(	\frac{3}{2048},  -\frac{3}{2048}+\frac{3\sqrt{2}i}{1024}, -\frac{3}{1024},  \frac{3}{2048}-\frac{3\sqrt{2}i}{1042}\right).$$
Then for $n$ large enough,  $$H_{2,1}\left(\mathbf{R}_n'\right)\supset \widetilde{U}.$$ 
Denote the map under this conjugacy by $T_1=H_{2,1}\circ H_1\circ \mathcal{F}\circ H_1^{-1}\circ H_{2,1}^{-1}$, then we obtain (\ref{formab}). Furthermore,
\(G_{1,j}(a,b) = \overline{G_{2,j}(b,a)}\) are homogeneous polynomials of degree $j$.
Let
$$G_{1,j}(a,b) = \sum_{k = 0}^{j }{G_{1,j}^{k}a^{k}b^{j - k}},$$
$$G_{2,j}(a,b) = \sum_{k = 0}^{j }{G_{2,j}^{k}b^{k}a^{j  - k}}$$
then we get \(G_{1,j}^{k} = \overline{G_{2,j}^{k}}\) and they are all of order \(\mathcal{O}\left( n^{- j+1} \right)\), since in (\ref{pre_Birk}) $F_{i,j}$ are homogeneous polynomials of degree $j+1$ for any $i,j$ and their norms are of order $\mathcal{O}(n^{-j})$. Direct calculation using $(\ref{Fij})$ shows $(\ref{Gij})$. 
\end{proof}

The second step is to find the first term of the Birkhoff normal form in complex coordinates.
\begin{lemma}
	There is a conjugacy $H_{2,2}: (a,b)\rightarrow (\xi, \eta)$ given by the following implicit form,
	\begin{equation}\label{conj_conjugacy}
		a = \xi + \sum_{m = 2}^{3}{\phi_{m}(\xi_1,\eta_1)},\ b = \eta + \sum_{m = 2}^{3}{\psi_{m}(\xi_1,\eta_1)},
	\end{equation}
	where \(\phi_{m}\) and \(\psi_{m}\) are homogeneous polynomials of
	degree \(m\), and under this conjugacy, $T_2=H_{2,2}\circ T_1 \circ H_{2,2}^{-1} $ takes the form
	\begin{align}\label{formBirk}
		\begin{split}
			\xi_1'= u\xi_1 + \mathcal{ O}_4\left( n^{- 3} \right),\\
			\eta_1'= v\eta_1 +\mathcal{ O}_4\left( n^{- 3} \right),\\ 
			u = \overline{v} = e^{i(\alpha + \alpha_{2}\xi_1\eta_1)},
		\end{split}
	\end{align}
	where \(e^{i\alpha} = \lambda\).
	
	Moreover, if we let
	\begin{equation}\label{domain of xi}
		V=0.98P\left(	\frac{3}{2048},  -\frac{3}{2048}+\frac{3\sqrt{2}i}{1024}, -\frac{3}{1024},  \frac{3}{2048}-\frac{3\sqrt{2}i}{1042}\right)
	\end{equation}
	be a complex parallelogram, and let $\widetilde{V}=\left\{(\xi,\eta)\in \mathbb{C}^2, \xi=\overline{\eta},\xi\in V\right\}$.
	Then for any $n$ large enough,
	$$H_{2,2}\left(\widetilde{U}\right)\supset \widetilde{V},$$ and we have,
	\begin{equation}\label{complex_conj}
		{\overline{\phi}}_{m}(\eta_1,\xi_1) = \psi_{m}(\xi_1,\eta_1).
	\end{equation}
	When $a$ and $b$ are complex conjugate numbers, $\xi_1$ and $\eta_1$ will also be complex conjugate numbers.
	
	Furthermore, the polynomials $\phi_m$ and $\psi_m$ are of order $\mathcal{O}(n^{-m+1})$. So, this conjugacy is almost identity when $n$ goes to infinity.
\end{lemma}
\begin{proof}
We assume the conjugacy is in the form of (\ref{conj_conjugacy}) and let
\(\phi_{m}(\xi_1,\eta_1) = \sum_{l = 0}^{m}a_{m}^{l}\xi_1^{l}\eta_1^{m - l}\),
\(\psi_{m}(\xi_1,\eta_1) = \sum_{l = 0}^{m}b_{m}^{l}\xi_1^{m - l}\eta_1^{l}\). Further,  we assume that 
\begin{equation}\label{conj_relation}
\overline{a_{m}^{l}} = b_{m}^{l},
\end{equation}
then condition (\ref{complex_conj}) will be satisfied.
Assume that this conjugacy will turn (\ref{formab}) into the form (\ref{formBirk}).

 We then solve the following equation. 
\begin{equation}
	H_{2,2}^{-1}\circ T_2=T_1 \circ H_{2,2}^{-1}.
\end{equation}
To solve this equation, we expand both sides into polynomials of
\((\xi_1,\eta_1)\) to the fourth degree, then the equation involving \(\xi_1\) turns out to be
\begin{equation}\label{Birk_eq}
	\begin{split}
&\lambda\left( \xi_1 + \sum_{m = 2}^{3}{\phi_{m}(\xi_1,\eta_1)} \right) + \sum_{l = 2}^{3}G_{1,l} \left( \xi_1 + \sum_{m = 2}^{3}{\phi_{m}(\xi_1,\eta_1)},\eta_1 + \sum_{m = 2}^{3}{\psi_{m}(\xi_1,\eta_1)} \right)\\
 & = \lambda\left( 1 + \sum_{m=1}^{\infty}{\left(i\alpha_{2}\xi_1\eta_1\right)}^m/m!  \right)\xi_1  \\
 &+ \sum_{m = 2}^{3}{\phi_{m}\left( \lambda\left( 1 + i\alpha_{2}\xi_1\eta_1 \right)\xi_1 +\mathcal{O}_4(n^{-3}), \lambda^{- 1}\left( 1 - i\overline{\alpha_{2}}\xi_1\eta_1 \right)\eta_1 +\mathcal{O}_4(n^{-3}) \right)}\\
 &+\mathcal{O}_4(n^{-3}),
 \end{split}
\end{equation}

It should be pointed out that the equation for \(\eta_1\) will give 
exact the complex conjugate version of (\ref*{Birk_eq}), so there is no need to list
it here. Our goal is to find $a_m^l$ and $b_m^l$ so that all the coefficients of homogeneous polynomials of $\xi_1$ and $\eta_1$ of degree less than four cancel out in (\ref{Birk_eq}). This, together with (\ref{conj_relation}), would lead us to the following equations:
\begin{align}
    {\rm coefficient\ of\ } \xi_1: &\ \lambda=\lambda,\notag \\
	{\rm coefficient\ of\ } \eta_1:&\ 0=0,\notag \\
	\label{2,0} {\rm coefficient\ of\ }\xi_1^2:& \left(\lambda^2-\lambda\right)a_2^2=G_{1,2}^2,\\
	\label{1,1} {\rm coefficient\ of\ }\xi_1\eta_1:& \left(1-\lambda\right)a_2^1=G_{1,2}^1,\\
	\label{0,2} {\rm coefficient\ of\ }\eta_1^2:& \left(\lambda^{-2}-\lambda\right)a_2^0=G_{1,2}^0,\\
	\label{3,0} {\rm coefficient\ of\ }\xi_1^3:& \left(\lambda^3-\lambda\right)a_3^3=G_{1,2}^1b_2^0+G_{1,3}^3+2a_2^2G_{1,2}^2,\\
	\label{2,1} {\rm coefficient\ of\ } \xi_1^2\eta_1:&\\
	  i\lambda\alpha_{2} + \left(\lambda-\lambda \right)a_{3}^{2}&= \left( a_{2}^{2} + b_{2}^{1} \right)G_{1,2}^{1} +2a_2^1G_{1,2}^2+2b_2^0G_{1,2}^0+ G_{1,3}^{2},\\
    \label{1,2} {\rm coefficient\ of\ }\xi_1\eta_1^2:&\\
     \left(\lambda^{-1}-\lambda\right)a_3^1&= \left( a_{2}^{1} + b_{2}^{2} \right)G_{1,2}^{1}+G_{1,3}^1+2a_2^0G_{1,2}^2+2b_2^1G_{1,2}^0,\\
	\label{0,3} {\rm coefficient\ of\ } \eta_1^3:& \left(\lambda^{-3}-\lambda\right)a_3^0=G_{1,2}^1a_2^0+G_{1,3}^0+2b_2^2G_{1,2}^0.
\end{align}
From this series of equations, we can see that from (\ref{2,0}) to (\ref{0,2}), we can directly solve the value of $a_2^0,a_2^1$ and $a_2^2$. Then plug in these values in (\ref{3,0}), (\ref{1,2}) and (\ref{0,3}) (and notice the conjugate relation (\ref{conj_relation})) , we can solve out $a_3^0,a_3^1$ and $a_3^3$ (notice that \(\lambda^{k} \neq 1\) for \(k = 1,2,3,4\) if
\(n\) is large enough). However, this process fails for $a_3^2$ since it has coefficient zero in (\ref{2,1}), and that is why $\alpha_2$ has to be involved in the normal form. We can then just let $a_3^2$ be zero, and find the exact value of $\alpha_2$ to satisfy the equation (\ref{2,1}). 

Now we find the value of $\alpha_2$.

From (\ref{2,1}),
\begin{align}\label{eq1}
\alpha_{2} &= -\frac{i}{\lambda}\left\lbrack \left( a_{2}^{2} + \overline{a_{2}^{1}} \right)G_{1,2}^{1}++2a_2^1G_{1,2}^2+2\overline{a_2^0}G_{1,2}^0  + G_{1,3}^{2} \right\rbrack,
\end{align}
and from (\ref{2,0}), (\ref{1,1}) and (\ref{0,2})
\begin{equation}\label{eq2}
 a_{2}^{2}=\frac{1}{\lambda( \lambda-1)}G_{1,2}^{2},\ \  a_{2}^{1}=\frac{1}{(1 - \lambda)}G_{1,2}^{1},\ \  a_{2}^{0}=\frac{1}{\lambda^{-2}-\lambda}G_{1,2}^{0}
\end{equation}
Combine (\ref{eq1}), and (\ref{eq2}), we have
\begin{equation}\label{alpha2}
\alpha_{2} = -i\left( \frac{G_{1,2}^{2}G_{1,2}^{1}}{\lambda^2( \lambda-1)} + \frac{\left|{G_{1,2}^{1}}\right|^{2}}{\lambda-1 } +\frac{2G_{1,2}^1G_{1,2}^2}{\lambda(1-\lambda)}+\frac{2\left|{G_{1,2}^{0}}\right|^{2}}{\lambda^3-1}+ \frac{G_{1,3}^{2}}{\lambda} \right).
\end{equation}

Moreover, we can see that $a_2^0,a_2^1$, $a_2^2$ are of order $\mathcal{O}(n^{-1})$ and $a_3^0,a_3^1$, $a_3^3$ are of order $\mathcal{O}(n^{-2})$. 
Therefore, we have the observation that $G_{i,j}^k$, $\phi_j$ and $\psi_j$ are of order $\mathcal{O}(n^{-j+1})$. Notice that we have checked that all polynomials of degree less than or equal to 3 coincide in equation \ref{Birk_eq}, then all remaining terms are polynomials of degree higher than 4 and the observation about $G_{i,j}^k$, $\phi_j$ and $\psi_j$ shows that they are all $\mathcal{O}_4\left(n^{-3}\right)$.
Thus, \ref{Birk_eq} is solved using this coordinate transform and the existence of conjugacy is proved.

Finally, we can see that $H_{2,2}$ is identity up to $\mathcal{O}(n^{-1})$ accuracy, thus $H_{2,2}(\widetilde{U})\sim \widetilde{U}$, up to $\mathcal{O}(n^{-1})$ accuracy. Since 
$$U=0.99P\left(	\frac{3}{2048},  -\frac{3}{2048}+\frac{3\sqrt{2}i}{1024}, -\frac{3}{1024},  \frac{3}{2048}-\frac{3\sqrt{2}i}{1042}\right),$$

$$V=0.98PP\left(	\frac{3}{2048},  -\frac{3}{2048}+\frac{3\sqrt{2}i}{1024}, -\frac{3}{1024},  \frac{3}{2048}-\frac{3\sqrt{2}i}{1042}\right),$$
we can conclude that if for $n$ large enough, $\widetilde{V}\subset H_{2,2}(\widetilde{U})$.

\end{proof}
\begin{remark}
	Combine $(\ref{Gij})$ and $(\ref{alpha2})$, we have
	$$\alpha_{2} = -\frac{1}{n^2}\frac{4\sqrt{2}}{9} +\mathcal{O}\left(\frac{1}{n^3}\right),$$ 
     Actually, we can even conclude that $\alpha_2$ is a real number because of the self-intersection property, as we shall see in the next paragraph.
\end{remark}

The last step to prove Proposition \ref{FindBirk} is to pull everything back to real Cartesian coordinates.

\begin{lemma}
	Let ${\rm Re}\alpha_2=\alpha_2^r$ and ${\rm Im}\alpha_2=\alpha_2^i$,
	The conjugacy $T_3$ is a linear conjugacy to pull everything back to $\mathbb{R}^2$, namely
\[H_{2,3}:\ \ \xi_1 = \xi + i\eta,\ \eta_1 = \xi - i\eta.\]
Then we have
\begin{align*}
	H_{2,3}\circ &T_2 \circ H_{2.3}^{-1}:\left( \xi,\eta \right) \mapsto \left( \xi',\eta' \right)\\
	\xi' &= \xi\nu\cos w - \eta\nu\sin w + \frac{1}{n^{3}}\mathcal{O}_{4},\\
	\eta' &= \xi\nu\sin w + \eta\nu\cos w + \frac{1}{n^{3}}\mathcal{O}_{4},\\
	w &= \alpha(n) + \alpha_{2}^{r}\left( \xi^{2} + \eta^{2} \right).\\
	\nu &=\exp\left(-\alpha_{2}^{i}\left( \xi^{2} + \eta^{2} \right)\right).
\end{align*}
\end{lemma}
\begin{proof}
	Straightforward calculation.
\end{proof}
\begin{proof}[Proof of proposition \ref{FindBirk}]
	
	Let us first show that $\alpha_2^i=0$, thus $\nu=1$.
	We prove this by contradiction. We can assume that $\alpha_2^i>0$, and the other case $\alpha_2^i<0$ can be treated similarly. Then we have 
	$${\xi'}^2+{\eta'}^2=\nu^2\left({\xi}^2+{\eta}^2\right)+\mathcal{O}_{6}=\left({\xi}^2+{\eta}^2\right)-\alpha_2^i{\left({\xi}^2+{\eta}^2\right)}^2+\mathcal{O}_{6},$$
	the second equality is due to Taylor expansion. Then for $(\xi,\eta)$ very close to the origin, ${\xi'}^2+{\eta'}^2< {\xi}^2+{\eta}^2$, the image of a very small circle centered at the origin lies strictly inside it and has no intersection with itself, which is contradictory. The case of $\alpha_2^i<0$ is similar, so $\alpha_2^i=0$.
\end{proof}	
	Let $H_2=H_{2,3} \circ H_{2,2} \circ H_{2,1}$. This is the required conjugacy.
Finally, observe that
$$\xi=Re(\xi_1),\ \eta= Im(\eta_1),$$
then the domain of $(\xi,\eta)$ is 
$$0.98P\left(	\left(\frac{3}{2048},0\right),\  \left(-\frac{3}{2048},\frac{3\sqrt{2}}{1024}\right),\ \left(-\frac{3}{2048},0\right),\ \left(\frac{3}{2048},\frac{3\sqrt{2}}{1024}\right)\right),$$
then it contains
$$ {\mathcal{E}}=P\left(\left(\frac{1}{1024},0\right),\  \left(-\frac{1}{1024},\frac{\sqrt{2}}{512}\right),\ \left(-\frac{1}{1024},0\right)\  \left(\frac{1}{1024},-\frac{\sqrt{2}}{512}\right)\right).$$
Moreover, straightforward calculation shows that 
$$\left(\xi,\eta\right)=\left(\frac{y_1}{2}+\mathcal{O}(n^{-1}),\frac{\sqrt{2}x_1}{4}-\frac{\sqrt{2}y_1}{2}+\mathcal{O}(n^{-1})\right)$$
\begin{proof}[Proof of {\color{red}Proposition \ref{normal form near fixed point}}]
Composing the two conjugacies $H_1$ in Proposition \ref{twisted form near fixed points} and $H_2$ in Proposition \ref{FindBirk}, we will obtain the required conjugacy. In addition, we have the following formula.
	\begin{align*}
		\xi&=\frac{y}{4}-\frac{1}{4}+\mathcal{O}(n^{-1}),\\
		\eta&=\frac{\sqrt{2}x}{4}-\frac{3\sqrt{2}n}{4}+\frac{3\sqrt{2}}{16}+\mathcal{O}(n^{-1}).
    \end{align*}
	This gives us the exact form in equation \ref{elllptic formula} and finishes the proof.
\end{proof}

\section{Proof of theorem \ref{main no escape}}\label{Proof no escape}
The proof in this section is divided into three parts. In the first subsection, we apply a change of coordinates to simplify the dynamics of the first return map to a fundamental domain. In the second subsection, we establish the proposition that for a family of maps defined on the cylinder, which includes the first return map we find in the first section, there is zero measure of escaping orbit. In the last section, we prove the proposition.

\begin{figure}[htbp]
	\centering
	\begin{tikzpicture}[scale=0.565]
		\draw[very thick] (cos{-67.5},sin{-67.5})--(cos{247.5},sin{247.5});
		\draw[very thick] (cos{-67.5},sin{-67.5}) arc (-67.5:247.5:1);
		\draw [dashed] (0,0)--(cos{-67.5},sin{-67.5});
		\draw [dashed] (0,0)--(cos{247.5},sin{247.5});
		\draw[domain=0:10] plot(\x, {sin{-67.5}}) node[below right] {\large $\ell_1$};
		\draw[domain=cos{-67.5}:9.5] plot (\x, {tan{22.5}*(\x-cos{-67.5})+sin{-67.5}}) node[below right] {\large  $\ell_2$};
		\draw[domain=cos{247.5}:9.5] plot (\x, {tan(-22.5) *(\x-cos(247.5))+sin{247.5}}) node[below right] {\large  $\ell_3$};
		\draw plot file {Fline1.txt} node[below left]{\large  $\ell'_1$};	
		\draw plot file {Fline2.txt} node[below left]{\large  $\ell'_2$};	
		\draw plot file {Fline3.txt} node[below left]{\large  $\ell'_3$};	
		\node  at (0,-0.7) {\tiny $2\beta$};
		\node at (8,0.8) {\large  $\Rone$};
		\node at (4,4.5) {\large  $\uppercase\expandafter{\romannumeral5}$};
		\node at (-7,4.8) {\large  $\uppercase\expandafter{\romannumeral2}$};
		\node at (-8,0.8) {\large  $\uppercase\expandafter{\romannumeral3}$};
		\node at (-4,-4.5) {\large  $\uppercase\expandafter{\romannumeral6}$};
		\node at (8,-2.8) {\large  $\uppercase\expandafter{\romannumeral4}$};	
	\end{tikzpicture}
	\caption{Continuous region  for $F_\beta^2$ of a circular sector outer billiard}
	\label{domain}
\end{figure}
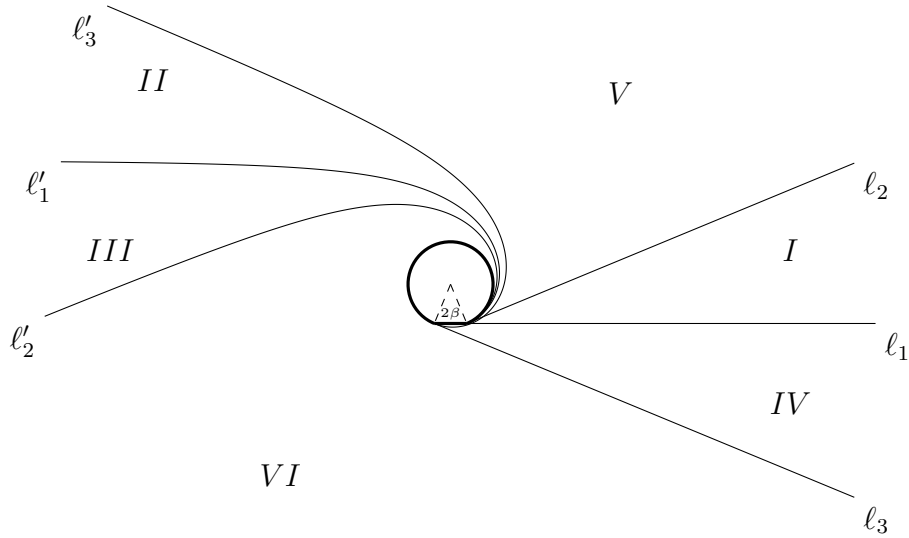
Let \(F_\beta\)
denote the outer billiard map of the circular sector with central angle $2\pi-2\beta$. 
As in figure \ref*{domain}, the plane is divided into six
regions, and in each region $F_\beta^2$ is smooth. We consider the infinite region
\(\mathcal{D}_\beta\) bounded by \(\ell_{1}\), \(F_\beta^{2}\ell_{1}\)
and \(\{ x = x_{0}\}\), for a large constant \(x_{0}\). Let \(\mathcal{F}_\beta\) denotes the first return map to the region
\(\mathcal{D}_\beta\). 

\subsection{The first return map to \texorpdfstring{$\mathcal{D}_\beta$}{}}

To prove Theorem \ref{main no escape}, we will use the following lemma. 
\begin{lemma}\label{keylemma}
	Define $A_\beta$, $B_\beta$ and $C_\beta$ to be
	\begin{align}
		&A_\beta=\frac{1}{6}\tan^3\left(\frac{\beta}{2}\right)+\frac{1}{2}\tan\left(\frac{\beta}{2}\right),\\
		&B_\beta=\left(\frac{(\pi-2\beta)}{4}+\frac{1}{3}\tan^3\left(\frac{\beta}{2}\right)+\tan\left(\frac{\beta}{2}\right)\right),\\
		&C_\beta=(2\sin\beta+\sin2\beta)B_\beta.
	\end{align}
	For $\beta\in(0,\pi/2]$, there exists a smooth change of coordinates $G:(x,y) \in \mathcal{D}_\beta \rightarrow \left(R,\phi\right) \in \textbf{C} = \left[R_0, \infty\right) \times \mathbb{T}$ of the form
	\begin{equation}
		\begin{split}
			&R = B_\beta \left[\frac{\cos\beta+1}{2}x+\frac{\sin\beta}{2}y\right]+\frac{1}{2}-B_\beta\frac{\sin\beta}{2}+\mathcal{O}\left(\frac{1}{x} \right),\\
			& \phi=\frac{y+\cos\beta}{2(\cos\beta+1)}+\mathcal{O}\left( \frac{1}{x}\right),
		\end{split}
	\end{equation}
	and a piecewise linear map on the cylinder $\textbf{C} = \left[R_0, \infty\right) \times \mathbb{T}$
	\begin{equation}
		\mathcal{L}_\beta\left(R,\phi\right)=\left(R+C_\beta\left(\left\{\phi-R\right\}-\frac{1}{2}\right), \left\{\phi-R\right\}\right).
	\end{equation}
	Moreover, the following properties hold.
	\item[$\left( 1\right)$]
	The singularities of $\mathcal{F}_\beta$ in the coordinate $(R,\phi)$ are piecewise curves in $\textbf{C}$ given by 
	\begin{equation}
		\begin{split}
			&\frac{1}{2}C_1=\left\{C_1R-\phi\right\}+\mathcal{O}\left(\frac{1}{R}\right),\\
			&1-\frac{1}{2}C_2=\left\{C_2R-\phi\right\}+\mathcal{O}\left(\frac{1}{R}\right),\\
			&0=\left\{R-\phi\right\}+\mathcal{O}\left(\frac{1}{R}\right)\\
		\end{split}
	\end{equation}
	where
	\begin{align*}
		C_1=\frac{A_\beta}{2A_\beta+\frac{\pi-2\beta}{4}},\ C_2=\frac{A_\beta+\frac{\pi-2\beta}{4}}{2A_\beta+\frac{\pi-2\beta}{4}}.\ 
	\end{align*}
	\item[$\left(2\right) $]
	If $\left(R,\phi\right)$ is not on the curves of singularities, we have
	$$\mathcal{F}_\beta\left(R,\phi\right)=\mathcal{L}_\beta^2\left(R,\phi\right)+\mathcal{O}(1/R).$$
\end{lemma}

\begin{remark}
	\item[$\left( 1\right)$] When $\beta=\pi/2$, $C_\beta=8/3$ and for $\beta\in(0,\pi/2)$, $C_\beta<4$. Therefore, the piecewise linear map $\mathcal{L}_\beta$ is always elliptic. Furthermore, $\lim_{\beta \to 0}{C_\beta}=0$. 
    \item[$\left( 2\right)$] The smooth change of coordinates is asymptotically affine when $x$ goes to infinity. Thus, the Lebesgue measure ${\rm Leb}$ in $\textbf{C}$ is equivalent to the Lebesgue measure $\mu$, which is inherited from the $\mathbb{R}^2$ plane, on $\mathcal{D}_\beta$. Therefore, they agree on the sets of zero measures.
    \item[$\left( 3\right)$] The curves of singularities have zero Lebesgue measure.
\end{remark}
The detailed proof of the lemma \ref{keylemma} is in Appendix \ref{proof of key lemma}.

\subsection{Zero measure of escaping orbits for a family of maps}
{\color{red}
\begin{main}[zero measure of escaping orbits for a family of maps]\label{key proposition}
Let $T_\Delta$ be a perturbation of piecewise linear transform on the cylinder  $\textbf{C} = \left[R_0, \infty\right) \times \mathbb{T}$:
\begin{align*}
    T_\Delta(R,\phi)&=\mathcal{L}_\Delta^k(R,\phi)+\mathcal{O}(R^{-1}),
\end{align*}
where $k$ is an positive integer, $\Delta\in(0,4)$ and
\begin{equation}
		\mathcal{L}_\Delta\left(R,\phi\right)=\left(R+\Delta\left(\left\{\phi-R\right\}-\frac{1}{2}\right), \left\{\phi-R\right\}\right).
\end{equation}
Let $E_\Delta=\left\{(R,\phi),\lim_{n\rightarrow\infty}T_\Delta^n(R,\phi)=\infty\right\}$ be the set of escaping orbits and ${\rm 
Leb}$ be the Lebesgue measure on the cylinder. Define the rotation number:
\begin{equation}
    \alpha_\Delta=\arcsin{\frac{\sqrt{\Delta(1-\Delta)}}{2}}.
\end{equation}
Then for $T_\Delta$ with $\alpha_\Delta=\pi/m$, where $m\geq 3$ is an integer,
${\rm Leb}(E_\Delta)=0.$
\end{main}
}
First let's show that Proposition \ref{key proposition} implies Theorem \ref{main no escape}.
\begin{proof}[Proof of Theorem \ref{main no escape}]
    we choose $\beta_m$ such that 
    $$\arcsin{\frac{\sqrt{C_{\beta_m}(1-C_{\beta_m})}}{2}}=\pi/m.$$ Let $\Delta=C_{\beta_m}$ and $k=2$. Then the first return map $\mathcal{F}_{\beta_m}$ given in Lemma \ref{keylemma} has the same form as $T_\Delta$ except for the singularity curves where $\mathcal{F}_{\beta_m}$ is not defined, which is of zero measure.
    Let $\bar{E}_{\beta_m}=\left\{(x,y)\in \mathcal{D}_\beta,\lim_{n\rightarrow\infty}\mathcal{F}_\beta^n(x,y)=\infty\right\}$ be the set of escaping orbits for the map $\mathcal{F}_\beta$.
    By Proposition \ref{key proposition}, we have $${\rm Leb}(E_\Delta)=\mu(\bar{E}_{\beta_m})=0.$$    
    Since $\mathcal{D}_\beta$ is a fundamental domain and for any $x\in E_\beta$, there exist $i_x\in\mathbb{N}$ that $F_\beta^{i_x}(x)\in \mathcal{D}_\beta$, then we have
\begin{align*}
	&E_\beta \subset \bigcup_{k=0}^{\infty}F_\beta^{-k}\left(\bar{E}_{\beta_m}\right),\\
	&\mu\left(E_\beta\right)=0.
\end{align*}
\end{proof}
 \begin{remark}
	Notice that every orbit of $F_\beta$ will go into the infinite region bounded by $\ell_1$ and $F_\beta\left(\ell_1\right)$. Since \(\mathcal{D}_\beta\) is defined to be the unbounded region bounded by \(\ell_{1}\), \(F_\beta^{2}\ell_{1}\)
and on the right-hand side of \(\{ x = x_{0}\}\), any escaping orbit of $F_\beta$ will go into \(\mathcal{D}_\beta\). Then if we prove that the escaping orbit of $\mathcal{F}_\beta$, namely the first return map to \(\mathcal{D}_\beta\), has zero measure, then Theorem \ref{main no escape} follows.
\end{remark}

\subsection{Proof of Proposition \ref{key proposition}}
The proof of Proposition \ref{key proposition} is made in three steps. 
\begin{itemize}
    \item The first step is to apply a linear conjugacy, and in the new coordinates, the linear part $\mathcal{L}_\Delta$ will be a piecewise rotation.
    \item The second step is to show that there is no escaping orbit for the piecewise-linear part $\mathcal{L}_\Delta$. 
    \item The last step is to observe that as $R$ goes to infinity, $\mathcal{L}_\Delta^k$ becomes a better and better approximation of $T_\Delta$ and applies the trick of Poincar{\' e} recurrence.  
\end{itemize}

\subsubsection{Coordinate transform}
We introduce the continuous region of $\mathcal{L}_\Delta$ in the cylinder $\textbf{C}$,  
\begin{align*}
	\textbf{C}_n=&\left\{(R,\phi)|\ n-1<R-\phi<n\right\},\\
	\mathcal{L}_\Delta|_{\textbf{C}_n}:R& \mapsto R+\Delta\left(\phi-R+n-\frac{1}{2}\right),\\
	\phi& \mapsto \phi-R+n.
\end{align*}
$\mathcal{L}_\Delta|_{\textbf{C}_n}$ can be written in matrix form:
\begin{align*}
	\begin{pmatrix}
		R\\
		\phi
	\end{pmatrix}
	\mapsto
	\begin{pmatrix}
		1-\Delta&\Delta\\
		-1&1
	\end{pmatrix}
	\begin{pmatrix}
		R\\
		\phi
	\end{pmatrix}
	+
	\begin{pmatrix}
		\Delta\left(n-\frac{1}{2}\right)\\
		n
	\end{pmatrix}.
\end{align*}

We then apply a linear transformation to make $\mathcal{L}_\Delta$ a piecewise rotation.
\begin{proposition}\label{linear}
	Apply the following linear conjugation on $\mathbb{R}^2$
	\begin{align}
		\begin{pmatrix}
			\overline{R}\\
			\overline{\phi}
		\end{pmatrix}
		=
		\begin{pmatrix}
			\frac{2}{\sqrt{\Delta}\,\sqrt{4-\Delta}} & -\frac{\sqrt{\Delta}}{\sqrt{4-\Delta}}\\
			0 & 1
		\end{pmatrix}
		\begin{pmatrix}
			R\\
			\phi
		\end{pmatrix}.
	\end{align}
	Then the region  $\mathbf{C}_n$ is mapped to the region
	\begin{align*}
		\mathbf{D}_n=\left\{\left(\frac{\sqrt{\Delta}\,\sqrt{4-\Delta}}{2}\overline{R}-(1+\frac{\Delta}{2})\overline{\phi}\right)\in(0,1)\right\},
	\end{align*}
	and $\mathcal{L}_\Delta$ is conjugated to a rotation $\mathcal{R}_\Delta$,  centered at $(\overline{R}_n,1/2)$,
	\begin{align*}
		&\mathcal{R}_\Delta|_{\mathbf{D}_n}
		\begin{pmatrix}
			\overline{R}\\
			\overline{\phi}
		\end{pmatrix}
		-
		\begin{pmatrix}
			\overline{R}_n\\
			\frac{1}{2}
		\end{pmatrix}
		=
		\begin{pmatrix}
			\cos\alpha_\Delta&\sin\alpha_\Delta\\
			-\sin\alpha_\Delta&\cos\alpha_\Delta
		\end{pmatrix}
		\left[
		\begin{pmatrix}
			\overline{R}\\
			\overline{\phi}
		\end{pmatrix}
		-
		\begin{pmatrix}
			\overline{R}_n\\
			\frac{1}{2}
		\end{pmatrix}
		\right],
	\end{align*}
	where $\alpha_\Delta=\arcsin(\sqrt{\Delta(4-\Delta)}/2)\in(0,\pi/2)$, $\overline{R}_n=\frac{4n-\Delta}{2\sqrt{\Delta}\sqrt{4-\Delta}}$. Furthermore, $T_\Delta$ is conjugated to $\overline{T}_\Delta$:
	\begin{align}\label{linear_rotate}
		\overline{T}_\Delta(\overline{R},\overline{\phi})=\mathcal{R}_\Delta^k(\overline{R},\overline{\phi})+\mathcal{O}\left(\frac{1}{\overline{R}}\right).
	\end{align}
\end{proposition}
\begin{proof}
	Straightforward calculation.
\end{proof}

Therefore $\mathcal{R}_\Delta$ on	$\textbf{D}_n$ is a clockwise rotation around the center of $\textbf{D}_n$, $\left(\overline{R}_n,\frac{1}{2}\right)$, as shown in Figure \ref * {rotation}. In the figure, $\textbf{D}_n$ is drawn with solid line segments and the image of $\textbf{D}_n$ is drawn with dotted line segments.
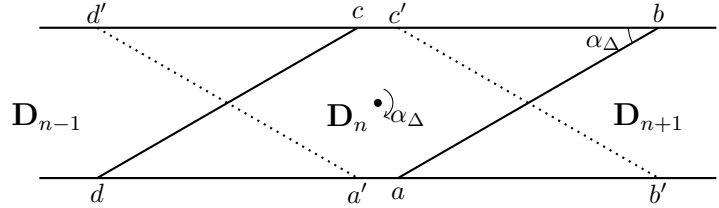
\begin{figure}[H]
	\centering
	\begin{tikzpicture}
		\draw [thick] (-4.5,-1)--(4.5,-1);
		\draw [thick] (-4.5,1)--(4.5,1);
		\draw [thick] ({2-sqrt(3)},-1)--({2+sqrt(3)},1);
		\draw [thick] ({-2-sqrt(3)},-1)--({-2+sqrt(3)},1);
		\fill (0,0) circle(1.5pt);
		\draw [dotted][thick] ({2-sqrt(3)},1)--({2+sqrt(3)},-1);
		\draw [dotted][thick] ({-2-sqrt(3)},1)--({-2+sqrt(3)},-1);
		\draw [->] (0.2*cos{67.5},0.2*sin{67.5}) arc (67.5:-67.5:0.2);
		\node at (-0.4,-0.2) {\large  $\textbf{D}_n$};
		\node at (0.4,-0.2) {\normalsize $\alpha_\Delta$};
		\node at (3.6,-0.2) {\large  $\textbf{D}_{n+1}$};
		\node at (-4.4,-0.2) {\large  $\textbf{D}_{n-1}$};
		\node at ({2-sqrt(3)},-1.2) {\normalsize $a$};
		\node at ({2+sqrt(3)},-1.2) {\normalsize $b'$};
		\node at ({2-sqrt(3)},1.2) {\normalsize $c'$};
		\node at ({2+sqrt(3)},1.2) {\normalsize $b$};
		\draw ({1.6+sqrt(3)},1) arc(180:210:0.4);
		\node at ({-2+sqrt(3)},-1.2) {\normalsize $a'$};
		\node at ({-2-sqrt(3)},-1.2) {\normalsize $d$};
		\node at ({-2+sqrt(3)},1.2) {\normalsize $c$};
		\node at ({-2-sqrt(3)},1.2) {\normalsize $d'$};
		\node at  (3,0.8) {\normalsize $\alpha_\Delta$};
	\end{tikzpicture}
	\caption{The piecewise linear map $\mathcal{R}_\Delta$}
	\label{rotation}
\end{figure}
\begin{remark}
	A straightforward calculation shows that $\textbf{D}_n$ is a rhombus with an interior angle equal to $\alpha_\beta$. In figure \ref*{rotation}, the rhombus $\diamond abcd$ is mapped to $\diamond a'b'c'd'$.
\end{remark}

\subsubsection{Invariant domains for the piecewise linear map \texorpdfstring{$\mathcal{R}_\Delta$}{TEXT}}
When $\alpha_\Delta=\pi/m$, since the edge $aa'$ is the image of itself under $2m-1$ times of rotation $\mathcal{R}_\Delta$ around $\left(\overline{R}_n,\frac{1}{2}\right)$, we can define $\textbf{O}_n$ to be the regular $2m$-polygon formed by rotation of the edge $aa'$

\begin{lemma}\label{no_excape}
	For $\left\{\Delta|\ \alpha_{\Delta}=\pi/m, m>1,m \in \mathbb{N} \right\}$,  $ \mathbf{O}_n$ are invariant under the piecewise linear map $\mathcal{R}_\Delta$.
\end{lemma}

\begin{proof}
	The proof is based on a key observation: the line
	$aa'$, after being rotated by $\mathcal{R}_\Delta$ for $m-2$ times around the center of $\textbf{D}_{n}$, coincides with the line $cc'$.
	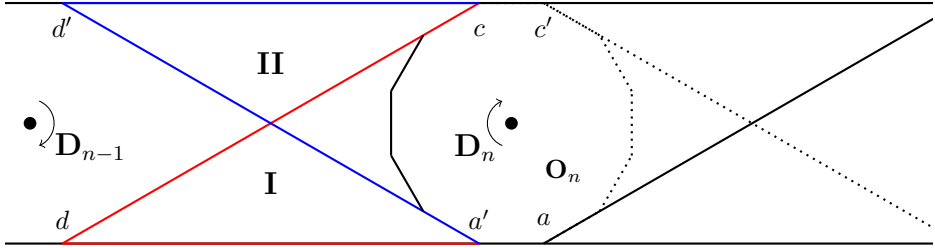
\begin{figure}[H]
		\centering
		\begin{tikzpicture}[scale=1.6]
			\clip (-6.2,-1.5) rectangle (1.5,1.5);
			\draw [thick] (-8,-1)--(3,-1);
			\draw [thick] (-8,1)--(3,1);
			\draw [thick] ({2-sqrt(3)-2},-1)--({2+sqrt(3)-2},1);
			\draw [red, thick] ({-2-sqrt(3)-2},-1)--({-2+sqrt(3)-2},-1);
			\draw [red, thick] ({-2-sqrt(3)-2},-1)--({-2+sqrt(3)-2},1);
			\fill (-2,0) circle(1.5pt);
			\fill (-6,0) circle(1.5pt);
			\draw [dotted][thick]({2-sqrt(3)-2},1)--({2+sqrt(3)-2},-1);
			\draw [blue,thick] ({-2-sqrt(3)-2},1)--({-2+sqrt(3)-2},-1);
			\draw [blue,thick] ({-2-sqrt(3)-2},1)--({-2+sqrt(3)-2},1);
			\draw [dotted][thick][rotate around={-330:(-2,0)}]({sqrt(3)-4},-1)--({-sqrt(3)},-1);
			\draw [thick][rotate around={-60:(-2,0)}]({sqrt(3)-4},-1)--({-sqrt(3)},-1);
			\draw [thick][rotate around={-90:(-2,0)}]({sqrt(3)-4},-1)--({-sqrt(3)},-1);
			\draw [thick][rotate around={-120:(-2,0)}]({sqrt(3)-4},-1)--({-sqrt(3)},-1);
			\draw [dotted][thick][rotate around={-180:(-2,0)}]({sqrt(3)-4},-1)--({-sqrt(3)},-1);
			\draw [dotted][thick][rotate around={-210:(-2,0)}]({sqrt(3)-4},-1)--({-sqrt(3)},-1);
			\draw [dotted][thick][rotate around={-240:(-2,0)}]({sqrt(3)-4},-1)--({-sqrt(3)},-1);
			\draw [dotted][thick][rotate around={-270:(-2,0)}]({sqrt(3)-4},-1)--({-sqrt(3)},-1);
			\draw [dotted][thick][rotate around={-300:(-2,0)}]({sqrt(3)-4},-1)--({-sqrt(3)},-1);
			\draw [->] ({0.2*cos{67.5}-6},0.2*sin{67.5}) arc (67.5:-67.5:0.2);
			\draw [->] ({0.2*cos{-247.5}-2},0.2*sin{-112.5}) arc (-112.5:-247.5:0.2);
			\node at (-2.3,-0.2) {\large  $\textbf{D}_n$};
			\node at (-5.5,-0.2) {\large  $\textbf{D}_{n-1}$};
			\node at ({-sqrt(3)},-0.8) {\normalsize $a$};
			\node at ({+sqrt(3)},-0.8) {\normalsize $b'$};
			\node at ({-sqrt(3)},0.8) {\normalsize $c'$};
			\node at ({+sqrt(3)},0.8) {\normalsize $b$};
			\node at ({-4+sqrt(3)},0.775) {\normalsize $c$};
			\node at ({-4+sqrt(3)},-0.8) {\normalsize $a'$};
			\node at ({sqrt(3)-3.3},-0.4) {\normalsize $\textbf{O}_{n}$} ;
			\node at (-4,-0.5) {\large $\textbf{\Rone}$};
			\node at (-4,0.5) {\large $\textbf{\Rtwo}$};
			\node at ({-4-sqrt(3)},-0.8) {\normalsize $d$};
			\node at ({-4-sqrt(3)},0.8) {\normalsize $d'$};
		\end{tikzpicture}
		\caption{Invariant polygon $\textbf{O}_n$. (The region $\Rone$ is bounded by the red edges $dc$, $da'$ and the edges of the polygon $\textbf{O}_n$; the region $\Rtwo$ is bounded by the blue edges $d'a'$, $dc$ and the edges of $\textbf{O}_n$).}
		\label{rotation1}
	\end{figure}

	As is shown in figure \ref*{rotation1},	let $\textbf{O}_{n}$ be the finite convex area in $\textbf{D}_{n}$ bounded by $\bigcup_{i=0}^{\infty}{\mathcal{R}_\Delta^i(aa')}$. Since dynamic in $\textbf{D}_{n}$ is a clockwise rotation by $\pi/m$, $a$ is a $2m$-periodic point, $\textbf{O}_{n}$ is a  regular $2m$ polygon, a $\mathcal{R}_\Delta$-invariant set. 
	\end{proof}
	
\begin{remark}
	We will see that the polygon $\textbf{O}_n$ blocks the orbits under map $\mathcal{R}_\Delta$ so that they cannot get across $\textbf{O}_n$.
	In contrast, Figure \ref{rotation2} shows the case where the rotation number of $\mathcal{R}_\Delta$ is not $\pi/m$ for some integer $m$. For example, when the rotation number is equal to $2\pi/11$, there is still a regular polygon formed rotations of $aa'$, but now this polygon is no longer invariant, since it has an intersection with $\textbf{D}_{n-1}$ and $\textbf{D}_{n+1}$. Thus, we fail to construct an invariant polygon that blocks the dynamics in this case.
	
	Thus, the theorem \ref{main no escape} and lemma \ref{no_excape} only stand in the very delicate situation when $\alpha_\beta=\pi/m$. In contrast, the result from Dolgopyat and Fayad's paper \cite{dolgopyat_unbounded_2009} shows that there are escaping balls for the semi-disc outer billiard map, and there method can be applied to the case where $\beta$ is close to $\pi/2$.
\end{remark}
	\pgfmathsetmacro{\alp}{360/11}
	\begin{figure}[H]
	\centering
	\begin{tikzpicture}[scale=1.6]	
		\clip ({-2*cos(\alp)-3.5},{-2*sin(\alp)-0.2}) rectangle ({0.5+2*cos(\alp)},{2*sin(\alp)});
        \draw [thick] ({-2*cos(\alp)-4},{2*sin(\alp)})--({2*cos(\alp)+1},{2*sin(\alp)});
        \draw [thick] ({-2*cos(\alp)-4},{-2*sin(\alp)})--({2*cos(\alp)+1},{-2*sin(\alp)});
		\draw [thick] ({-2*cos(\alp)-4},{-2*sin(\alp)})--({2*cos(\alp)-4},{2*sin(\alp)});
		\draw [thick] ({-2*cos(\alp)},{-2*sin(\alp)})--({2*cos(\alp)},{2*sin(\alp)});
		\draw [thick] ({-2*cos(\alp)-4},{2*sin(\alp)})--({2*cos(\alp)-4},{-2*sin(\alp)});
		\draw [thick] ({-2*cos(\alp)},{2*sin(\alp)})--({2*cos(\alp)},{-2*sin(\alp)});
		\fill (-2,0) circle(1.5pt);
		\fill (2,0) circle(1.5pt);
		\draw [densely dash dot][thick][rotate around={{-2*\alp}:(-2,0)}]({2*cos(\alp)-4},{-2*sin(\alp)})--({-2*cos(\alp)},{-2*sin(\alp)});
		\draw [densely dash dot][thick][rotate around={{-3*\alp}:(-2,0)}]({2*cos(\alp)-4},{-2*sin(\alp)})--({-2*cos(\alp)},{-2*sin(\alp)});
		\draw [densely dash dot][thick][rotate around={{-4*\alp}:(-2,0)}]({2*cos(\alp)-4},{-2*sin(\alp)})--({-2*cos(\alp)},{-2*sin(\alp)});
		\draw [densely dash dot][thick][rotate around={{-5*\alp}:(-2,0)}]({2*cos(\alp)-4},{-2*sin(\alp)})--({-2*cos(\alp)},{-2*sin(\alp)});
		\draw [densely dash dot][thick][rotate around={{-6*\alp}:(-2,0)}]({2*cos(\alp)-4},{-2*sin(\alp)})--({-2*cos(\alp)},{-2*sin(\alp)});
		\draw [densely dash dot][thick][rotate around={{-7*\alp}:(-2,0)}]({2*cos(\alp)-4},{-2*sin(\alp)})--({-2*cos(\alp)},{-2*sin(\alp)});
		\draw [densely dash dot][thick][rotate around={{-8*\alp}:(-2,0)}]({2*cos(\alp)-4},{-2*sin(\alp)})--({-2*cos(\alp)},{-2*sin(\alp)});
		\draw [densely dash dot][thick][rotate around={{-9*\alp}:(-2,0)}]({2*cos(\alp)-4},{-2*sin(\alp)})--({-2*cos(\alp)},{-2*sin(\alp)});
		\draw [densely dash dot][thick][rotate around={{-10*\alp}:(-2,0)}]({2*cos(\alp)-4},{-2*sin(\alp)})--({-2*cos(\alp)},{-2*sin(\alp)});
		\draw [densely dash dot][thick][rotate around={{-11*\alp}:(-2,0)}]({2*cos(\alp)-4},{-2*sin(\alp)})--({-2*cos(\alp)},{-2*sin(\alp)});
		\draw [densely dash dot][thick][rotate around={{-12*\alp}:(-2,0)}]({2*cos(\alp)-4},{-2*sin(\alp)})--({-2*cos(\alp)},{-2*sin(\alp)});
		\draw [->] ({0.2*cos{67.5}-2},0.2*sin{67.5}) arc (67.5:-67.5:0.2);
		\draw [->] ({0.2*cos{-247.5}+2},0.2*sin{-112.5}) arc (-112.5:-247.5:0.2);
		\node at (-2.3,-0.2) {\large  $\textbf{D}_n$};
		\node at (1.5,-0.2) {\large  $\textbf{D}_{n+1}$};
		\node at ({-sqrt(3)},-0.8) {\normalsize $a$};
		\node at ({+sqrt(3)},-0.8) {\normalsize $b'$};
		\node at ({-sqrt(3)},0.8) {\normalsize $c'$};
		\node at ({+sqrt(3)},0.8) {\normalsize $b$};
		\node at ({-4+sqrt(3)},0.775) {\normalsize $c$};
		\node at ({-4+sqrt(3)},-0.8) {\normalsize $a'$};
		\node at ({sqrt(3)-3.3},-0.4) {\normalsize $\textbf{O}_{n}$} ;
	\end{tikzpicture}
	\caption{Case where one fails to construct an invariant polygon that blocks the dynamics}
	\label{rotation2}
\end{figure}
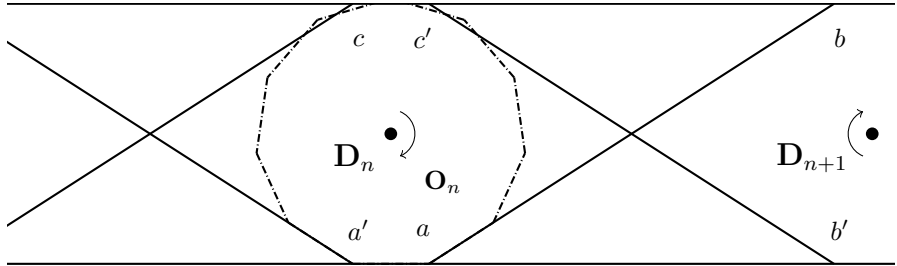

 Observe that figure \ref{rotation1}, $\textbf{C}\setminus\textbf{O}_{n}$ has two components, one on the left of $\textbf{O}_n$ and the other on the right. Denote them by $L\left(\textbf{O}_n\right)$ and $R\left(\textbf{O}_n\right)$. Then we have the following corollary:
	\begin{corollary}\label{connected argument}
	For fixed integer $n$ and any integer $n'<n$, $\mathcal{R}_\Delta\left(\textbf{D}_{n'}\right) \subset L\left(\textbf{O}_n\right)$
    \end{corollary}	
    \begin{proof}
    	Since $\mathcal{R}_\Delta$ on $\textbf{D}_{n'}$ is continous, then $\mathcal{R}_\Delta\left(\textbf{D}_{n'}\right)$ is connected. In addition, $\textbf{D}_{n'} \cap \mathcal{R}_\Delta\left(\textbf{D}_{n'}\right)\ne \varnothing$, then $L\left(\textbf{O}_n\right) \cap \mathcal{R}_\Delta\left(\textbf{D}_{n'}\right)\ne \varnothing$ and then $\mathcal{R}_\Delta\left(\textbf{D}_{n'}\right) \subset L\left(\textbf{O}_n\right)$.
   	\end{proof}
%
%
	\begin{lemma}\label{linear nonescape}
		The linearized map $\mathcal{R}_\Delta$ on $\textbf{C}= \left[R_0, \infty\right) \times \mathbb{T}$ has no escaping orbit.
	\end{lemma}	
\begin{proof}
 As shown in the figure \ref{rotation}, the region $\Rone$, which is bounded by the red edges $dc$, $da'$ and the edges of the polygon $\textbf{O}_n$, is mapped to the region $\Rtwo$, which is bounded by the blue edges $d'a'$, $dc$ and the edges of $\textbf{O}_n$.
	
	Now suppose that there is an escaping point $x\in L\left(\textbf{O}_n\right)$, then there exists a positive integer $n$ such that $\mathcal{R}_\Delta^n(x)\in R\left(\textbf{O}_n\right)$. Since we can choose $n$ to be the minimal integer that satisfies this property, we can additionally assume that $\mathcal{R}_\Delta^{n-1}(x)\notin R\left(\textbf{O}_n\right)$. Since $\textbf{O}_n$ is an invariant polygon for the map $\mathcal{R}_\Delta$, we conclude that $\mathcal{R}_\Delta^{n-1}(x)\notin \textbf{O}_n$, then $\mathcal{R}_\Delta^{n-1}(x)\in L\left(\textbf{O}_n\right)$. From the last corollary, $\mathcal{R}_\Delta^{n-1}(x)\in L\left(\textbf{O}_n\right)\cap \textbf{D}_n$. 
	
  However, we can see from figure \ref{rotation}, $\mathcal{R}_\Delta^{n-1}(x)\in L\left(\textbf{O}_n\right)\cap \textbf{D}_n$, namely the region $\Rone$, and it is mapped to the region $\Rtwo$, so we can conclude that $\mathcal{R}_\Delta^{n}(x)\in\Rtwo$, which has empty intersection with $R\left(\textbf{O}_n\right)$, contradiction.
	
\end{proof}

\subsubsection{Poincar\'{e} recurrence to a finite measure set}
In this section, we will show that the Lebesgue measure ${\rm Leb}(\overline{E}_{\Delta})=0$ for $\left\{\Delta|\ \alpha_{\Delta}=\pi/m, m>1,m \in \mathbb{N} \right\}$, where \(\overline{E}_{\Delta}\) = \{\(x\in\textbf{C}\): $\lim_{n \to \infty}\overline{T}_{\Delta}^n(x)=\infty$\}, then Proposition \ref{key proposition} follows.

When $\alpha_\beta=\pi/m$, we construct a concentric $2m$-regular polygon $\textbf{Q}_{4^n}$ in each $\left\{\textbf{O}_{4^n},n\geq N\right\}$ for some large $N$ that scales down $\textbf{O}_{4^n}$ by $1-1/2^n$. 

Then we define
\begin{align*}
	G_n&\coloneqq \textbf{O}_{4^n}\setminus\textbf{Q}_{4^n},\\
	G&\coloneqq \bigcup_{n=N}^{\infty}G_n.
\end{align*}
Obviously, $m(G)=m\left(\textbf{O}_n\right)\sum_{i=N}^{\infty}1-\left({1-\frac{1}{2^{i}}}\right)^2<\infty$. From (\ref{linear_rotate}) we know that there for large $N$, if $(\overline{R},\overline{\phi})\in\textbf{D}_{4^n},n\geq N$:
\begin{align}\label{estimate}
	\left|\overline{T}_\Delta(\overline{R},\overline{\phi})-\mathcal{R}^k_\Delta(\overline{R},\overline{\phi})\right|=\mathcal{O}\left({\frac{1}{4^n}}\right)<\frac{1}{2}\frac{1}{2^n}\frac{1-\cos\alpha_\Delta}{\sin\alpha_\Delta},
\end{align}
since $\left(1-\cos\alpha_\Delta\right)/{\sin\alpha_\Delta}$ is a positive constant not depending $n$.
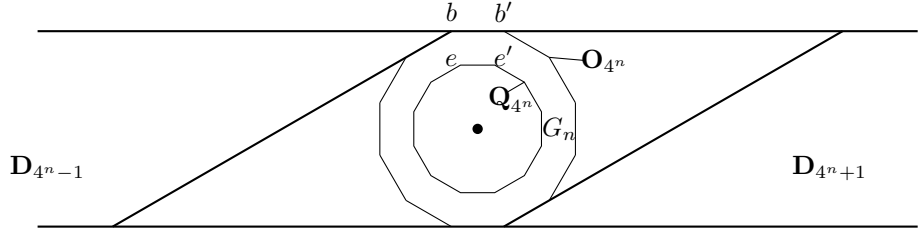
\begin{figure}[H]
	\centering
	\begin{tikzpicture}[scale=1.3]
		\draw [thick] (-4.5,-1)--(4.5,-1);
		\draw [thick] (-4.5,1)--(4.5,1);
		\draw [thick] ({2-sqrt(3)},-1)--({2+sqrt(3)},1);
		\node at ({2-sqrt(3)},1.2) {$b'$};
		\node at ({-2+sqrt(3)},1.2) {$b$};
		\node at ({2-sqrt(3)},0.75) {$e'$};
		\node at ({-2+sqrt(3)},0.7) {$e$};
		\draw [thick] ({-2-sqrt(3)},-1)--({-2+sqrt(3)},1);
		\fill (0,0) circle(1.5pt);
		\node at (3.6,-0.4) { $\textbf{D}_{4^n+1}$};
		\node at (-4.4,-0.4) { $\textbf{D}_{4^n-1}$};
		\node at (1.3,0.7) { $\textbf{O}_{4^n}$};
	    \node at (0.35,0.3) { $\textbf{Q}_{4^n}$};
	     \node at (0.84,0) { $G_{n}$};
	     \draw ({sqrt(3)-1},{sqrt(3)-1})--(1.1,0.7);
	     \draw (0.48,0.48)--(0.31,0.38);
		\draw  [rotate around={-30:(0,0)}]({-2+sqrt(3)},-1)--({2-sqrt(3)},-1) ;\draw  [rotate around={-60:(0,0)}]({-2+sqrt(3)},-1)--({2-sqrt(3)},-1) ;\draw  [rotate around={-90:(0,0)}]({-2+sqrt(3)},-1)--({2-sqrt(3)},-1) ;\draw [rotate around={-120:(0,0)}]({-2+sqrt(3)},-1)--({2-sqrt(3)},-1) ;
		\draw [rotate around={60:(0,0)}]({-2+sqrt(3)},-1)--({2-sqrt(3)},-1) ;\draw  [rotate around={90:(0,0)}]({-2+sqrt(3)},-1)--({2-sqrt(3)},-1) ;\draw  [rotate around={120:(0,0)}]({-2+sqrt(3)},-1)--({2-sqrt(3)},-1) ;\draw  [rotate around={150:(0,0)}]({-2+sqrt(3)},-1)--({2-sqrt(3)},-1) ;
		\node[draw,minimum size={50},regular polygon,regular polygon sides=12] (a)  {};
	\end{tikzpicture}
	\caption{Graph of $G_n$}
	\label{Poincaresection}
\end{figure}
\begin{lemma}\label{across}
	Every orbit going to infinity has nonempty intersection with the set $G$, that is, for every $x \in \overline{E}_\Delta$, there exist positive integers $n_x,k_x$, such that $T^{k_x}_\Delta(x)\in G_{n_x}\cap\overline{E}_\Delta$. 
\end{lemma}

\begin{proof} 
	From the fact that $\overline{E}_\Delta$ is $\overline{T}_\Delta$ invariant, we only need to show that $x$'s orbit goes into some $G_{n_x}$.
	
	For $x \in \overline{E}_\Delta$, choose a positive integer $n_x\geq N$, such that $x$ is on the left of $\textbf{O}_{2^{n_x}}$. Since the orbit of $x$ will go to infinity, it must cross the polygon $\textbf{O}_{2^{n_x}}$. However, from (\ref{estimate}), the orbit of $x$ under the map $\overline{T}_\Delta$ will jump less than $\frac{1}{2}\frac{1}{2^{n_x}}\frac{1-\cos\alpha_\Delta}{\sin\alpha_\Delta}$, which is the half of the width of the band $G_n$. Then the orbit of $X$ cannot pass across the polygon $\textbf{O}_{n_x}$ without comming inside ${G}_{n_x}$.
\end{proof}

\begin{proof}[Proof of Proposition \ref{key proposition}]\
Now consider the set 
\begin{align*}
    \overline{G}_n & \coloneqq G_n\cap \overline{E}_\Delta,    \\
    \overline{G}&\coloneqq \bigcup_{n=N}^{\infty}\overline{G}_n.
\end{align*}
If $m\left(\bar{G}\right)>0$, then Lemma \ref{across} guarantees that we can define the first return map of $\overline{T}_\Delta$ to $\overline{G}$. Now Poincar\'{e}'s recurrence theorem shows that almost every point of any $\bar{G}_n$ will go back to itself, which contradicts the definition of $\overline{G}_n$. Then ${\rm Leb}\left(\overline{G}\right)=0$.
 Lemma \ref{across} gives us
$$\overline{E}_\Delta\subset\bigcup_{k=1}^{\infty}\overline{T}_\Delta^{-k}(G),$$
 which implies 
 $${\rm Leb}\left(\overline{E}_\Delta\right)=0.$$
Since $\overline{T}_\Delta$ is obtained by simply applying a linear conjugacy to $T_\Delta$, thus
$${\rm Leb}\left({E}_\Delta\right)=0.$$
\end{proof}

\section{Application to Fermi Ulam model}

Dolgopyat and De Simoi studied the dynamics of piecewise smooth Fermi-Ulam models in their work \cite{de2012dynamics}. This model was first studied by Fermi and Ulam to explain the occurrence of cosmic rays. This is a model that describes the motion of an elastic ball that bounces between two walls, with distance $\ell(t)$ between them, where $\ell(t)$ is a 1-periodic function. Fermi and Ulam conjectured that there exist orbits whose energy grows to infinity with time. From KAM theory, it is proved that this conjecture is actually false for smooth enough $\ell(t)$.
Dolgopyat and De Simoi studied a class of piecewise smooth Fermi-Ulam models and found a normal form that describes the high-energy dynamics. They found a key variable $\Delta$ that characterizes the dynamics of the ping pong for large energies. To be more detailed, they proved that for $\Delta \in \left(0,4\right)$, the system looks regular for large energies, whereas for $\Delta \notin \left[0,4\right]$, the system is chaotic for large energies. Our result in this paper can be applied to the case where $\Delta \in \left(0,4\right)$.

By studying the orbit only at the moments of collisions with the moving wall, they simplified the model into a collision map $f$ in the collision space $\mathcal{M}$, which preserves a volume form $\omega$. Then they studied the system by finding $F$, the first return map to an infinite strip region $R$ inside $\mathcal{M}$ (here we use the same notation as in \cite{de2012dynamics}). Here we cite their results about the first return map as follows:

Let $t$ denote the collision time and since $\ell$ is periodic, $t$ can be defined on the torus $\mathbb{T}=\mathbb{R}/\mathbb{ Z}$. Let $v$ be the velocity of the ball after the collision. Let $\mathbb{A}=\mathbb{T} \times \mathbb{R}$, then the collision space is given by

$$\mathcal{M}=\left\{\left(t,v\right)\in\mathbb{A} \ {\rm s.t.} \ v>-\dot{\ell}\left(t\right)\right\}$$

Then the collision map is defined as $f:\mathcal{M} \rightarrow \mathcal{M},$
$$f\left(t_n,v_n\right)=\left(t_n+\delta t\left(t_n,v_n\right),v_n-2\dot{\ell}\left(t_n+\delta t\left(t_n,v_n\right)\right)\right)=\left(t_{n+1},v_{n+1}\right),$$
where $\delta t$ solves the functional equation
$$\delta t\left(t,v\right)=\frac{\ell(t)+\ell(t+\delta t(t,v))}{v}.$$
Then it can be observed that the volume form $\omega=\left(v+\dot{\ell}(t)\right){\rm d}t\wedge {\rm d}v$ is preserved.

\begin{proposition}\label{cite}
	There exist smooth coordinates $\left(\tau,I\right)$ on $R$ such that the first return map of $f$ on $R$ is given by
	$$F\left(\tau,I\right)=\hat{F}\left(\tau,I\right)+F_1\left(\tau,I\right)+r\left(\tau,I\right)$$
	where $\hat{F}\left(\tau,I\right)=\left(\bar{\tau},\bar{I}\right)$ with
	\begin{align*}
		\bar{\tau}=\tau-I \ \ \ {\rm mod 1,} \ \ \ \ \ \ \ \bar{I}=I+\Delta\left(\bar{\tau}-1/2\right),  
	\end{align*}
   $F_1$ is a correction of order $\mathcal{O}\left(I^{-1}\right)$ of the form 
   $$F_1\left(\tau,I\right)=I^{-1}\left(0,\Delta_1\left(\left(\bar{\tau}-1/2\right)^2-1/12\right)\right)$$
   and $r=\mathcal{O}_5\left(I^{-2}\right)$. Finally, $\omega=d\tau \wedge dI$.
\end{proposition} 

\begin{remark}
	The parameter $\Delta$ is determined by the following formula:
	$$\Delta=\ell(0)\left(\dot{\ell^+}(0)-\dot{\ell^-}(0)\right)\int_{0}^{1}\ell^{-2}(s)ds.$$ 
	For more details about Proposition \ref{cite}, see \cite{de2012dynamics}.
	It should also be pointed out that $\hat{F}$ covers a map on the torus $\tilde{F}:\mathbb{T}^2\rightarrow\mathbb{T}^2$. 
\end{remark}

Then they introduce the set of escaping orbits, namely
$$\varepsilon=\left\{\left(t_0,v_0\right):v_n\rightarrow\infty \ {\rm as}\ n\rightarrow +\infty\right\}.$$
We then collect some main theorems from Dolgopyat and De simoi's paper \cite{de_Simoi_2012}.

\begin{proposition}\label{cite1}
	If $\Delta \in \left(0,4\right)$, and $\tilde{F}$ has a stable accelerating orbit, then $m\left(\varepsilon\right)$=0.
\end{proposition}

\begin{proposition}\label{cite2}
		If $\Delta \notin \left[0,4\right]$, then $m\left(\varepsilon\right)$=0.
\end{proposition} 

Therefore, when $\Delta \in \left(0,4\right)$, the linear part of $F$ is elliptic and when $\Delta \notin \left[0,4\right]$ the linear part is hyperbolic.

Observe that the map $F$ is exactly of the form the maps we studied in Prop \ref{key proposition} with $k=1$ and
$$\alpha_\Delta=\arcsin{\frac{\sqrt{\Delta }\,\sqrt{4-\Delta }}{2}}.$$
And the invariant measure is given by $\omega=d\tau \wedge dI$, which is exactly the Lebesgue measure of the cylinder.

Therefore, we can apply Proposition \ref{key proposition} and obtain the following result:
 
 \begin{theorem}
 	For any integer $m \geq 3$, when $\alpha_\Delta=\pi/m$, namely,
 	$$\Delta=2-2\cos\left(\frac{\pi}{m}\right),$$
 	the escaping orbit for the corresponding Fermi-Ulam model has zero measure of escaping orbits.
 \end{theorem}

\section{Conclusion and Further Open Questions}\label{Open questions}

In this paper, we first showed that there are countable many elliptic islands filling up a positive proportion of the plane by finding a sequence of elliptic periodic points and verifying that at these points the twist term in the Birkhoff normal form is nonzero. Therefore, it is natural to ask if we can generalize this result to generic circular sector outer billiard maps $\mathcal{F}_\beta$.

It is easy to find a sequence of elliptic periodic points for any $\beta<\pi/2$, namely $(R,\phi)=(n,1/2)+\mathcal{O}(n^{-1})$ in Lemma \ref{keylemma}, then one may use the same method to verify if the twist term at these points are nonzero to find elliptic islands near them. There is no additional technique that needs to be introduced to address this problem.

\begin{conjecture}
	Is it true that for generic $\beta$, in the outer billiard map of the circular sector with angle $\beta$, there are elliptic islands that fill a positive proportion of the plane?
\end{conjecture}

Our second finding is that there is zero measure of escaping orbits in a class of circular sector outer billiard maps, where the corresponding $\alpha_\beta$ is $\pi/m$ for some positive integer $m$. This is a rare situation for $\beta<\pi/2$. On the other hand, in the paper \cite{dolgopyat_unbounded_2009}, Dolgopyat and Fayad show that there is an escaping ball in the semi-disc outer billiard map.
To find an escaping ball in circular sector outer billiards, one may also follow the method introduced in \cite{dolgopyat_unbounded_2009}, and find an escaping orbit that shifts by $1$ under some irritations of the linearized map. It should be noticed that it may take many more steps for an orbit to shift by 1 since the rotation number will be smaller and smaller when $\beta$ approaches $0$, which makes it much harder to find such an orbit.
However, if this is true, one can follow the method in \cite{dolgopyat_unbounded_2009} to  prove that there is an open set escaping to infinity in the corresponding outer billiard map for the circular sector. 
\begin{conjecture}
When $\alpha_\beta\neq\pi/n, \beta\leq\pi/2$, can we find an escaping ball for the outer billiard map $F_\beta$?
\end{conjecture}

\section*{Appendix}\label{appendix}

\appendix
\section{Proof of Proposition \texorpdfstring{\ref{twisted form near fixed points}}{}}\label{proof of proposition B}
\subsection{Coordinate exchange in regions $\Rone\sim\Rfour$}
We are going to find new coordinates $(\rho,\psi)$ that simplify the dynamics of $F^2$ in the regions $\Rone\sim\Rfour$. The main idea is from \cite{dolgopyat_unbounded_2009}. We first find the expression of the map $F^2$ in the polar coordinate in the following lemma:

We first consider the outer billiard map in the coordinates \((r,\theta)\), 
which are polar coordinates in the upper regions $\uppercase\expandafter{\romannumeral1}$, $\uppercase\expandafter{\romannumeral5}$, and
$\uppercase\expandafter{\romannumeral2}$, and polar coordinates shifted by \(\pi\) in the lower regions
$\uppercase\expandafter{\romannumeral3}$, $\uppercase\expandafter{\romannumeral6}$, $\uppercase\expandafter{\romannumeral4}$. Thus, both in the upper plane and in the lower plane $0\leq\theta\leq\pi$.
The following lemma gives the behavior of $F^2$ in the coordinates $(r,\theta)$.

\begin{lemma}\label{r-theta expansion}
	We consider the polar coordinates\((r,\theta)\) in the upper region and the polar coordinate
shifted by \(\pi\) in the lower region. Then the map
	\(\text{F}^{\text{2}}\) takes the following form.
	
	\begin{equation}\label{map_F^2}
		F^{2}(r,\theta) = \left( \begin{array}{r}
			r + a(\theta) + \frac{a_{1}(\theta)}{r} + \frac{a_{2}(\theta)}{r^{2}}+\frac{a_3(\theta)}{r^3}\mathcal{+ O}\left( \frac{1}{r^{4}} \right), \\
			\theta + \frac{b(\theta)}{r} + \frac{b_{1}(\theta)}{r^{2}} + \frac{b_{2}(\theta)}{r^{3}}+\frac{b_3(\theta)}{r^4}\mathcal{+ O}\left( \frac{1}{r^{5}} \right)
		\end{array} \right)
	\end{equation}
	where
	\begin{equation*}
		a(\theta) = - 2\cos\theta,\ b(\theta) = 2\left( 1 + \sin\theta \right),\ b_{1}(\theta) = 4\cos{\theta\left( 1 + \sin\theta \right)}.
	\end{equation*}
	And in regions $\uppercase\expandafter{\romannumeral1}$ and $\uppercase\expandafter{\romannumeral4}$, 
	\begin{equation*}
		a_{1}(\theta) = 2\sin^{2}\theta, a_{2}(\theta) = 4\cos\theta - 4\cos^{3}\theta, b_{2}(\theta) = 6\cos{2\theta} + \frac{8}{3}\sin3\theta + \frac{5}{3},
	\end{equation*}
	in regions $\uppercase\expandafter{\romannumeral2}$ and $\uppercase\expandafter{\romannumeral3}$,
	\begin{align*}
		a_{1}(\theta) = 2\sin^{2}\theta + 4\sin\theta,\ & a_{2}(\theta) = 8\cos\theta + 4\sin(2\theta) - 4\cos^{3}\theta,\\
		b_{2}(\theta) = 8\cos{2\theta} &+ \frac{8}{3}sin3\theta + 4sin\theta - \frac{1}{3},
	\end{align*}
	in regions $\uppercase\expandafter{\romannumeral5}$,
	\begin{align*}
		&a(\theta) = - 4\cos\theta, \ a_{1}(\theta) = 8\sin^{2}\theta,\\
		&a_{2}(\theta) = 32\cos\theta - 32\cos^{3}\theta, \ b(\theta) = 4\sin\theta,\\ 
		&b_{1}(\theta) = 8\sin{2\theta}, \ b_{2}(\theta) = \frac{64}{3}\sin3\theta.
	\end{align*}
    For $a_3(\theta)$ and $b_3(\theta)$ we do not need the explicit formula.
\end{lemma}

\begin{proof}
	This follows from directly expanding the outer billiard map into power series of \(1/r\)
to a higher order. It can be noticed that the form of \(a(\theta)\),
	\(b(\theta)\), \(a_{1}(\theta)\) and \(b_{1}(\theta)\) coincide with the
results of \cite{dolgopyat_unbounded_2009}.
\end{proof}
\begin{lemma}\label{rho_psi}
	There is a change in the coordinates from the polar coordinates \((r,\theta)\) to \((\rho,\psi)\). In regions $\uppercase\expandafter{\romannumeral1}$ and $\uppercase\expandafter{\romannumeral3}$,
	\((\rho,\psi) \in (0,\infty) \times \lbrack 0,1/3 + \mathcal{O}(1/r)\rbrack\)
	and in regions $\uppercase\expandafter{\romannumeral2}$ and $\uppercase\expandafter{\romannumeral4}$,
	\((\rho,\psi) \in (0,\infty) \times \lbrack 1/3 + \mathcal{O}(1/r),2/3\rbrack\), such
	that if \(\text{F}^{\text{2}}\) sends \((\rho,\psi)\) to \((\rho',\psi')\), then in the regions $\uppercase\expandafter{\romannumeral1} \sim \uppercase\expandafter{\romannumeral4}$
	\begin{align}\label{coord_1}
		\rho' = \rho\mathcal{+ O}\left( \frac{1}{\rho^4} \right),
		\psi' = \psi + \frac{1}{\rho}\mathcal{+ O}\left( \frac{1}{\rho^5} \right),\end{align}
	The coordinates change is of the following form:
	\begin{equation}\label{coor_pho}
		\rho = r\Phi_{1}(\theta) + \Phi_{2}(\theta) + \frac{\Phi_{3}(\theta)}{r} + \frac{\Phi_{4}(\theta)}{r^{2}},\end{equation}
	\begin{equation}\label{coor_psai}
		\psi = \Psi(\theta) + \frac{\Psi_{1}(\theta)}{r} + \frac{\Psi_{2}(\theta)}{r^{2}} + \frac{\Psi_{3}(\theta)}{r^{3}}.
	\end{equation}
	The functions \(\Phi_{i}(i \geq 2)\) and \(\Psi_{i}(i \geq 1)\) can take different forms in the regions $\uppercase\expandafter{\romannumeral1}\sim \uppercase\expandafter{\romannumeral4}$
\end{lemma}

	{\color{red}
    \begin{proof}It should be pointed out that we do not require \(\Phi_{i}\) and \(\Psi_{i}\) functions to be
continuous between different regions. We solve them separately in different reions. It turns out that all functions \(\Phi_{i}\) and \(\Psi_{i}\) can be solved explicitly from some first order ODEs. 
	
	We follow the method in \cite{dolgopyat_unbounded_2009} to calculate \(\rho' - \rho\) and expand
	it into powers of \(r\).
    If $F^2$ sends $(r,\theta)$ to $(r',\theta')$, by $(\ref{map_F^2})$, one obtains the following:
    \begin{align*}
    r'=r + a(\theta) + \frac{a_{1}(\theta)}{r} + \frac{a_{2}(\theta)}{r^{2}}+\frac{a_3(\theta)}{r^3}\mathcal{+ O}\left( \frac{1}{r^{4}}\right),\\
    \theta'=	\theta + \frac{b(\theta)}{r} + \frac{b_{1}(\theta)}{r^{2}} + \frac{b_{2}(\theta)}{r^{3}}+\frac{b_3(\theta)}{r^4}\mathcal{+ O}\left( \frac{1}{r^{5}} \right)
    \end{align*}

    Thus,
    \begin{align}\label{expanding}
        & {\rho'} - \rho \notag\\ 
        &=(r'-r)\Phi_{1}(\theta) + \Phi_{3}(\theta)\left(\frac{1}{r'}-\frac{1}{r}\right) + \Phi_{4}(\theta)\left(\frac{1}{r'^2}-\frac{1}{r^2}\right)\notag \\
        &+r'\left(\Phi_{1}(\theta')-\Phi_1(\theta)\right) + \Phi_{2}(\theta')-\Phi_2(\theta) + \frac{\Phi_{3}(\theta')-\Phi_3(\theta)}{r'} + \frac{\Phi_{4}(\theta')-\Phi_4(\theta)}{r'^{2}} 
        \end{align}
        Using Taylor expansion, we obtain
        \begin{align}
        \rho'-\rho= & \Phi_{1}^{'}b + \Phi_{1}a + \left( \Phi_{2}^{'}b+\left(\Phi_{1}^{'}b_{1} + \frac{\Phi_{1}^{''}b^{2}}{2} + a\Phi_{1}^{'}b + \Phi_{1}a_{1}  \right) \right)\frac{1}{r} \notag\\
         &+\frac{1}{r^2}\biggl(\Phi_3'b+\biggl(\Phi_1'''b^2+\Phi_1''bb_1+\Phi_1'\left(a_1b+ab_1+b_2\right)+\Phi_1a_2+\frac{\Phi_2''b^2}{2}\notag\\
        & +\Phi_2'b_1-\Phi_3a\biggr)\biggr) 
         +\mathcal{O}\left(\frac{1}{r^3}\right).
    \end{align}
Therefore we obtain the equation for $\Phi_1,\Phi_2$ and $\Phi_3$:
    \begin{align}\label{equations about phi}
     &\Phi_1'b+\Phi_1a=0,\notag \\
         &\Phi_{2}^{'}b+\left(\Phi_{1}^{'}b_{1} + \frac{\Phi_{1}^{''}b^{2}}{2} + a\Phi_{1}^{'}b + \Phi_{1}a_{1}\right) =0,\notag \\ 
         &\Phi_3'b+\biggl(\frac{\Phi_1'''b^3}{6}+\Phi_1''\left(bb_1+\frac{ab^2}{2}\right)+\Phi_1'\left(a_1b+ab_1+b_2\right)+\Phi_1a_2 \notag\\
         &+\frac{\Phi_2''b^2}{2}+\Phi_2'b_1-\Phi_3a\biggr)=0.
    \end{align}
    We use the initial condition that $\Phi_1$ takes the value 1 and \(\Phi_{i}(i \geq 2)\) take the value 0 when \(\theta = 0\) in the region $\uppercase\expandafter{\romannumeral1}$ and $\uppercase\expandafter{\romannumeral3}$, and when \(\theta = \pi\) in the region $\uppercase\expandafter{\romannumeral2}$ and $\uppercase\expandafter{\romannumeral4}$. Notice that these are all first-order odes, and $b(\theta)=2(1+\sin(\theta))$ is always positive in our coordinates. The first equation is just about $\Phi_1$ and we can explicit solve that $\Phi_1=1+\sin\theta$. Once we know $\Phi_1$, we can plug it in the second equation to obtain $\Phi_2$, and then plug $\Phi_1$ and $\Phi_2$ into the third equation to obtain $\Phi_3$. 
    It can be solved that
	\begin{align*}
		\Phi_{1} &= 1 + \sin\theta,\\
		\Phi_{2} &= 1 - |\cos\theta|, \ {\rm in\ region}\ \emph{\uppercase\expandafter{\romannumeral1}}\ {\rm and}\  \emph{\uppercase\expandafter{\romannumeral4}},\\
		\Phi_{2}& = -1 + |\cos\theta|, \ {\rm in\ region}\ \emph{\uppercase\expandafter{\romannumeral2}}\ {\rm and}\  \emph{\uppercase\expandafter{\romannumeral3}}.\\
		\Phi_{3} &= \frac{2\sin^3{\theta}-\sin\theta}{6\left(\sin{\theta}+1\right)}, \ {\rm in\ region}\ \emph{\uppercase\expandafter{\romannumeral1}}\ {\rm and}\  \emph{\uppercase\expandafter{\romannumeral4}},\\
		\Phi_{3}& = \frac{2\sin^3\theta-12\sin^2\theta-\sin\theta}{6\left(\sin\theta+1\right)}, \ {\rm in\ region}\ \emph{\uppercase\expandafter{\romannumeral2}}\ {\rm and}\  \emph{\uppercase\expandafter{\romannumeral3}}.
	\end{align*}
    If we expand $\rho'-\rho$ further so that the error term is $\mathcal{O}(r^{-4})$ and let the coefficient before $r^{-3}$ be zero, we would obtain the equation for $\Phi_4$. We don't have to find the explicit form of this equation, but we claim that it will still be a first order ODE like equations in $(\ref{equations about phi})$. More specifically, we have
    \begin{equation}
     \Phi_4'b+f(\Phi_1,\Phi_1',\Phi_1'',\Phi_1''',\Phi_1^{(4)},\Phi_2,\Phi_2',\Phi_2'',\Phi_2''',\Phi_3,\Phi_3',\Phi_3'',\Phi_4,\theta)=0.
    \end{equation}
    Here $f$ is some analytic function. Since $b(\theta)$ is always positive, the claim guarantees that we can find an analytic function $\Phi_4$ satisfying the requirement of the lemma with the initial condition that \(\Phi_{4}\) take the value 0 when \(\theta = 0\) in the region $\uppercase\expandafter{\romannumeral1}$ and $\uppercase\expandafter{\romannumeral3}$, and when \(\theta = \pi\) in the region $\uppercase\expandafter{\romannumeral2}$ and $\uppercase\expandafter{\romannumeral4}$. Then in the coordinates $(\rho,\psi)$, $\rho'=\rho+\mathcal{O}(r^{-4})=\rho+\mathcal{O}(\rho^{-4})$.
    
    The claim can be proved by observation from $(\ref{expanding})$. Notice that the only contribution of derivatives of $\Phi_4$ comes from the term 
    \begin{align*}
        \frac{\Phi_4(\theta')-\Phi_4(\theta)}{r'^2}=\frac{\Phi_4'(\theta)(\theta'-\theta)}{r'^2}+\frac{\Phi_4''(\theta)(\theta'-\theta)^2}{2r'^2}+\mathcal{O}\left(\frac{(\theta'-\theta)^3}{r'^2}\right).
    \end{align*}
    Observe that $\theta'-\theta=b/r+\mathcal{O}(r^{-2})$, therefore only  $\Phi_4'b$ will appear in the coefficient before the term $r^{-3}$. The same argument shows that only the first four derivatives of $\Phi_1$, the first three derivatives of $\Phi_2$ and the first two derivatives of $\Phi_3$ will appear in this equation. The claim is then proved.

    Similarly, we have 
    \begin{align}\label{expanding2}
        \psi'-&\psi-\frac{1}{\rho}\notag\\
        =&-\frac{1}{\rho}+\Phi_{1}(\theta)\left(\frac{1}{r'}-\frac{1}{r}\right) + \Psi_{2}(\theta)\left(\frac{1}{r'^2}-\frac{1}{r^2}\right)+\Psi_{3}(\theta)\left(\frac{1}{r'^3}-\frac{1}{r^3}\right)\notag \\
        &+ \Psi(\theta')-\Psi(\theta) + \frac{\Psi_{1}(\theta')-\Psi_1(\theta)}{r'} + \frac{\Psi_{2}(\theta')-\Psi_2(\theta)}{r'^2} + \frac{\Psi_{3}(\theta')-\Psi_3(\theta)}{r'^{3}}\notag \\
        =&\frac{1}{r}\left(b\Psi'-\frac{1}{\Phi_1}\right)+\frac{1}{r^2}\left(b\Psi_1'+\left(\frac{b^2\Psi''}{2}+b_1\Psi'+\frac{\Phi_2}{\Phi_1^2}-a\Psi_1\right)\right) \notag \\
        &+\frac{1}{r^3}\biggl(b\Psi_2'+\biggl(\frac{b^3\Phi_1^3\Psi'''}{6}+bb_1\Psi''+b_2\Psi'+\frac{b^2\Psi_1''}{2}+\Psi_1'\left(b_1-ab\right)-2a\Psi_2\notag \\
        &+\frac{\Phi_1\Phi_3-\Phi_2^2}{\Phi_1^3}\biggr)\biggr)+\mathcal{O}\left(\frac{1}{r^4}\right).
    \end{align}
    Therefore we obtain the equation for $\Psi,\Psi_1$ and $\Psi_2$:
    \begin{align}\label{equations about psi}
       & b\Psi'-\frac{1}{\Phi_1}=0, \notag \\
       &b\Psi_1'+\left(\frac{b^2\Psi''}{2}+b_1\Psi'+\frac{\Phi_2}{\Phi_1^2}-a\Psi_1\right)=0, \notag \\
       &b\Psi_2'+\biggl(\frac{b^3\Phi_1^3\Psi'''}{6}+bb_1\Psi''+b_2\Psi'+\frac{b^2\Psi_1''}{2}+\Psi_1'\left(b_1-ab\right) \notag\\
       &-2a\Psi_2+\frac{\Phi_1\Phi_3-\Phi_2^2}{\Phi_1^3}\biggr)=0. 
    \end{align}
    Still, the first equation is only about $\Psi$, and after we fix an initial value condition, we can solve $\Psi$. Then we plug it into the second equation to solve $\Psi_1$, and plug $\Psi$ and $\Psi_1$ into the third equation to solve $\Psi_2$.
    We use the following initial condition for these ODEs:
    \begin{align}\label{special initial cond}
        &\Psi(0)=0\ {\rm in\ region\ } \Rone\ {\rm and\ } \Rthree, \ \ \Psi(\pi)=\frac{2}{3}\ {\rm in\ region\ } \Rtwo\ {\rm and\ } \Rfour, 
    \end{align}
     and for $j=1$ and $2$,
    \begin{align}
        &\Psi_j(0)=0\ {\rm in\ region\ } \Rone\ {\rm and\ } \Rthree, \ \ \Psi_j(\pi)=0\ {\rm in\ region\ } \Rtwo\ {\rm and\ } \Rfour. 
    \end{align}
    Then we can solve $\Psi$, $\Psi_1$, and $\Psi_2$ one by one, and obtain the following.
    \begin{align}
        \Psi(\theta)=\frac{1}{6}-\frac{\cos^2\theta+2\cos\theta+\sin\theta\cos\theta-2\sin\theta-2}{6\left(\sin\theta+1\right)^2},
    \end{align}
    The explicit expression of $\Psi_1$ and $\Psi_2$ is too lengthy to be listed here, see the supplementary material \emph{Calculation Details} for detailed results.
    
     Similarly, if we expand $\psi'-\psi$ further so that the error term is $\mathcal{O}(r^{-5})$ and let the coefficient before $r^{-4}$ be zero, we would obtain the equation for $\Psi_3$. Still, we do not list the equation here, but this equation about $\Psi_3$ will look like the equations in $(\ref{equations about psi})$. More specifically, we will have the following.
     \begin{equation}
     \Psi_3'b+g(\Psi,\Psi',\Psi'',\Psi''',\Psi^{(4)},\Psi_1,\Psi_1',\Psi_1'',\Psi_1''',\Psi_2,\Psi_2',\Psi_2'',\Psi_3,\theta)=0,
    \end{equation}
    where $g$ is some analytic function.
    It should be pointed out that here $\Phi_1 \sim\Phi_4$ are regarded as known functions of $\theta$. This claim is proved by an argument similar to the one we have for $\Phi_4$ and we omit it here. Then $\psi'-\psi=\mathcal{O}(r^{-5})=\mathcal{O}(\rho^{-5})$.
\end{proof}
}
\begin{remark}
  The idea of finding the first return map is to simplify the dynamics in the continuous regions of $F^2$ and study the behavior of the orbit when it passes through these regions. We will see in the  proof of Proposition \ref{twisted form near fixed points} in the next section that we only need the values of $\Phi_i$ and $\Psi_i$ and their derivatives near the discontinuity lines to calculate the first return map. These values are given in the table in the Appendix \ref{important values}. Notice that we do not have the exact values of $\Phi_4(\pi/2)$ and $\Psi_3(\pi/2)$ and they are denoted by $e_1\sim e_4$ and $f_1\sim f_4$. In the next section, we will see that these values are not needed to prove Proposition \ref{twisted form near fixed points}.
\end{remark}
\subsection{{\color{red}Proof of Proposition \texorpdfstring{\ref{twisted form near fixed points}}{}}}\label{calculation}
This section is to give a proof of Proposition B.
    For a positive integer $n$, we consider the point $$(x^{(0)},y^{(0)})=\left(3n-\frac{1}{4},1\right)$$ and its rectangular neighborhood in the region $\mathcal{D}$:
        $$\mathbf{R}_{n}=\left(3n-\frac{1}{4},1\right)+\mathbf{R}$$
    Recall that we have defined $\mathbf{R}$ to be the parallelogram with vertices $\left(0,\frac{3}{512}\right)$,  $\left(\frac{3}{256},-\frac{3}{512}\right)$, $\left(0,-\frac{3}{512}\right)$ and $\left(-\frac{3}{256},\frac{3}{512}\right)$.
    We will see in the following proof that the orbit of points in $\mathbf{R}_{n}$ will go into the region $\Rfive$ but not into the region $\Rsix$. In this process, we track the trajectory of the point $(x^{(0)},y^{(0)})$ and observe that it is almost a fixed point of $\mathcal{F}$ (up to $\mathcal{O}(1/n)$ error). 

    The coordinates $(\rho,\psi)$ simplify the dynamics of $F^2$ in the continuous regions $\Rone\sim\Rfour$. We only need to study the behavior of the orbit when it passes through these regions, where we  change back to the original coordinates $(r,\theta)$. Therefore, only the values of $\Phi_i$ and $\Psi_i$ and their derivatives near the discontinuity lines are needed to calculate the first return map.

    We break down the first return map $\mathcal{F}$ into four parts. Recall that we define 
\(\mathcal{D}\) to be the infinite region bounded by \(\ell_{1}\), \(F^{2}\ell_{1}\),
and \(\{ x = x_{0}\}\), where $x_0$ is some large constant. Let $y_0$ be another large constant. We define three similar fundamental domains like $\mathcal{D}$: 
\begin{itemize}
    \item the region $\mathcal{D}_1$ is the infinite region bounded by $\ell_2$, $F^2(\ell_2)$ and $\left\{y=y_0\right\}$, 
    \item the region $\mathcal{D}_2$ is the infinite region bounded by the line $y=0,x\leq0$, its image under $F^2$ and $\left\{x=-x_0\right\}$.
    \item the region $\mathcal{D}_3$ is the infinite region bounded by $\ell_3$, $F^2(\ell_3)$ and $\left\{y=-y_0\right\}$. 
\end{itemize}
Let $\mathcal{F}_1$ be the passage from $\mathcal{D}$ to $\mathcal{D}_1$, $\mathcal{F}_2$ be the passage from $\mathcal{D}_1$ to $\mathcal{D}_2$, $\mathcal{F}_3$ be the passage from $\mathcal{D}_2$ to $\mathcal{D}_3$ and $\mathcal{F}_4$ be the passage from $\mathcal{D}_3$ to $\mathcal{D}_4$. 
    
    We introduce a new variable $\phi$ for convenience of exposition. Let $\phi=\rho\psi$ in the region $\mathcal{D}$ and $\mathcal{D}_2$, and $\phi=\rho\left(\psi-1/3\right)$ in the region $\mathcal{D}_1$ and $\mathcal{D}_3$. {\color{red} We will study how $\mathcal{F}_1\sim\mathcal{F}_4$ behave in the $(\rho,\phi)$ coordinates. We carry out the detailed computation for the map $\mathcal{F}_1$, while for the remaining ones, we omit some intermediate steps and present the final results. One can go to the supplementary material \emph{Calculation Details} for a detailed version of this proof.}

    \begin{proof}[Proof of Proposition B]  
            First notice that in the region $\mathcal{D}$, $$\theta\in\left[0,\frac{2}{r}+\mathcal{O}\left(\frac{1}{r^2}\right)\right].$$
        Then
         \begin{equation}\label{xy}
        \begin{split}
        x&=r\cos\theta=r\left(1-\frac{\theta^2}{2}+\mathcal{O}\left(\theta^4\right)\right)=r+\mathcal{O}\left(\frac{1}{r}\right),\\
        y&=r\sin\theta=r\left(\theta+\mathcal{O}\left(\theta^3\right)\right)=r\theta+\mathcal{O}\left(\frac{1}{r}\right).
        \end{split}
        \end{equation}
        And
         \begin{equation}\label{rhopsi}
        \begin{split}
        \rho&=r\Phi_1(0)+r\Phi_1'(0)\theta+\Phi_2(0)+\mathcal{O}\left(\frac{1}{r}\right)=r+r\theta+\mathcal{O}\left(\frac{1}{r}\right),\\
        \psi&=\Psi(0)+\Psi'(0)\theta+\frac{\Psi_1(0)}{r}+\mathcal{O}\left(\frac{1}{r^2}\right)=\frac{1}{2}\theta+\mathcal{O}\left(\frac{1}{r^2}\right).
        \end{split}
        \end{equation}
        We introduce $\phi=\rho\psi$ in the region $\mathcal{D}$. Combining $(\ref{xy})$ and $(\ref{rhopsi})$, we obtain the following:
        \begin{equation}\label{xytorhophi}
        \begin{split}
        \rho&=x+y+\mathcal{O}\left(\frac{1}{x}\right),\\
        \phi&=\rho\psi=\frac{1}{2}y+\mathcal{O}\left(\frac{1}{x}\right).
        \end{split}
        \end{equation}
        Recall that We fix a point $$(x^{(0)},y^{(0)})=(3n-\frac{1}{4},1)$$ in
        $$\mathbf{R}_{n}=\left(3n-\frac{3}{4},1\right)+\mathbf{R},$$
        for large positive integers $n$. Then correspondingly in the $(\rho,\phi)$ coordinates
        $$(\rho^{(0)},\phi^{(0)})=(3n+\frac{1}{4},\frac{1}{2})+\mathcal{O}\left(\frac{1}{n}\right),$$ 
	and
	\begin{equation}\label{initialrhophi}
	    \mathbf{R}_n=\left(3n+\frac{1}{4}-\frac{9}{1024},3n+\frac{1}{4}+\frac{9}{1024}\right)\times\left(\frac{1}{2}-\frac{3}{2048},\frac{1}{2}+\frac{3}{2048}\right)+\mathcal{O}\left(\frac{1}{n}\right)
	\end{equation}
    
        \emph{Step} 1. \emph{The map $\mathcal{F}_1$}. 
        
        Recall that, in the coordinates \((\rho,\psi)\), the discontinuity lines of \(\text{F}^{\text{2}}\) are
given by the following equations:
	
	\begin{align}
		\ell_{2} \subset \left\{ \psi = \frac{1}{3} - \frac{1}{3\rho} + \mathcal{O}\left( \frac{1}{\rho^{2}} \right) \right\},
		\ell_{3}^{'} \subset \left\{ \psi = \frac{1}{3} - \frac{5}{6\rho} + \mathcal{O}\left( \frac{1}{\rho^{2}} \right) \right\}.
	\end{align}
	
	A point starting in $\mathbf{R}_n$
	with coordinates \((\rho,\psi)\) takes \(m\) step to leave the region
	\emph{\uppercase\expandafter{\romannumeral1}}, where
	\begin{equation}
		m = \left\lbrack \frac{\rho}{3} - \rho\psi \right\rbrack=\left\lbrack \frac{\rho}{3} - \phi \right\rbrack,
	\end{equation}
     and we have
	\begin{equation}\label{cross1}
		\begin{split}
			\rho_{m} &= \rho + \mathcal{O}\left( \frac{1}{\rho^{3}} \right),\\
			\psi_{m} &= \frac{1}{3} - \frac{v}{\rho} + \mathcal{O}\left( \frac{1}{\rho^{4}} \right).
		\end{split}
	\end{equation}
    where 
    \begin{equation}\label{initialv}
         v=\left\{\frac{\rho}{3}-\phi\right\}\in\left(\frac{7}{12}-\frac{9}{2048},\frac{7}{12}+\frac{9}{2048}\right)+\mathcal{O}\left(\frac{1}{n}\right)\subset\left(\frac{1}{3},\frac{5}{6}\right)
    \end{equation}
    for large $n$. Then every point in $\mathbf{R}_n$ will enter the region $\Rfive$. And 
    \begin{equation}\label{initialv0}
         v^{(0)}=\left\{\frac{\rho^{(0)}}{3}-\phi^{(0)}\right\}=\frac{7}{12}+\mathcal{O}\left(\frac{1}{n}\right).
    \end{equation}
    Switching to $(r,\theta)$ coordinates and iterating once more, the orbit will visit the region $\mathcal{D}_1$. {\color{red} Notice now that $\psi_m$ is $\mathcal{O}(\rho^{-1})$ close to $1/3$. Since $\Psi(\pi/2)=1/3$, this implies that $\theta$ is $\mathcal{O}(r^{-1})$ close to $\pi/2$. Thus, Taylor expansion at $\theta=\pi/2$ will help us to return to the $(r,\theta)$ coordinates.

    To change from \((\rho_{m},\psi_{m})\) back to \((r_{m},\theta_{m})\), we assume that $(r_m,\theta_m)$ has the form 
	\begin{equation}
		\begin{split}
			r_{m} &= \rho +c_2 + \frac{c_3}{\rho} + \frac{c_4}{\rho^{2}}\mathcal{+ O}\left( \frac{1}{\rho^{3}} \right),\\
			\theta_{m} &= \frac{\pi}{2} + \frac{d_2}{\rho} + \frac{d_3}{\rho^{2}} + \frac{d_4}{\rho^{3}} + \mathcal{O}\left( \frac{1}{\rho^{4}} \right),
		\end{split}
	\end{equation}
where $c_2\sim c_4$ and $d_2\sim d_4$ are undetermined constants. Plug it in the Taylor expansion near \(\{\theta = \frac{\pi}{2}\}\):
	
	\begin{equation}\label{Taylorrho}
    \begin{split}
		\rho_{m} =& r_{m}\Phi_{1}(\frac{\pi}{2}) + r_{m}\Phi_{1}^{'}(\frac{\pi}{2})\left( \theta_{m} - \frac{\pi}{2} \right) \\
		&+ \frac{1}{2}r_{m}{\Phi}_{1}^{(2)}(\frac{\pi}{2})\left( \theta_{m} - \frac{\pi}{2} \right)^{2}+ \frac{1}{6}r_{m}{\Phi}_{1}^{(3)}(\frac{\pi}{2})\left( \theta_{m} - \frac{\pi}{2} \right)^{3} \\
		& + \Phi_{2}(\frac{\pi}{2}) + \Phi_{2}^{'}(\frac{\pi}{2})\left( \theta_{m} - \frac{\pi}{2} \right)+ \frac{1}{2}{\Phi}_{2}^{''}(\frac{\pi}{2})\left( \theta_{m} - \frac{\pi}{2} \right)^{2} \\
		&+\frac{\Phi_{3}(\frac{\pi}{2})}{r_{m}} + \frac{\Phi_{3}^{'}(\frac{\pi}{2})\left( \theta_{m} - \frac{\pi}{2} \right)}{r_{m}} + \frac{\Phi_{4}(\frac{\pi}{2})}{r_{m}^{2}}\mathcal{+ O}\left( \frac{1}{r_{m}^{3}} \right),
        \end{split}
	\end{equation}
	\begin{equation}\label{Taylorpsi}
    \begin{split}
		\psi_{m} = &\Psi(\frac{\pi}{2}) + \Psi^{'}(\frac{\pi}{2})\left( \theta_{m} - \frac{\pi}{2} \right) + \frac{\Psi^{(2)}(\frac{\pi}{2})}{2}\left( \theta_{m} - \frac{\pi}{2} \right)^{2} \\
		&+ \frac{\Psi^{(3)}(\frac{\pi}{2})}{6}\left( \theta_{m} - \frac{\pi}{2} \right)^{3} 
		+ \frac{\Psi_{1}(\frac{\pi}{2})}{r_{m}} + \frac{\Psi_{1}^{'}(\frac{\pi}{2})\left( \theta_{m} - \frac{\pi}{2} \right)}{r_{m}} \\
		& + \frac{\Psi_{1}^{''}(\frac{\pi}{2})}{2}\frac{\left( \theta_{m} - \frac{\pi}{2} \right)^{2}}{r_{m}}  
		+\frac{\Psi_{2}(\frac{\pi}{2})}{r_{m}^{2}} \\
		&+ \frac{\Psi_{2}^{'}(\frac{\pi}{2})\left( \theta_{m} - \frac{\pi}{2} \right)}{r_{m}^{2}} + \frac{\Psi_{3}(\frac{\pi}{2})}{r_{m}^{3}}\mathcal{+ O}\left( \frac{1}{r_{m}^{4}} \right).
        \end{split}
	\end{equation}
	The result should be consistent with $(\ref{cross1})$. Plugging in the values from Table \ref{values of Phi} and Table \ref{values of Psi} in Appendix \ref{important values}, we obtain the equations for $c_2 \sim c_4$ and $d_2\sim d_4$:
\begin{align*}
    &c_2+\frac{1}{2}=0,\ 3d_2+24v=2,\\
    &24c_3=3d_2^2-12d_2-2,\ 9d_3=9d_2-12c_2-10,\\
    &12c_4=3c_2d_2^2+3d_3d_2+2c_2-24e_1-6d_3,\\
    &18d_4 =1152f_1 +24c_3 +33d_2  +36c_2 d_2 +3d_2 ^3-48c_2 ^2 -80c_2 -12d_2 ^2 -18d_3.
\end{align*}
Solve these equations, then we have 

\begin{equation}\label{inverse Taylor}
		\begin{split}
			r_{m} &= \frac{1}{2}\rho - \frac{1}{2} + \frac{8v^2+\frac{8v}{3}-\frac{13}{36}}{\rho} + \frac{8v^2+\frac{32v}{9}-2e_1-\frac{23}{108}}{\rho^{2}}\mathcal{+ O}\left( \frac{1}{\rho^{3}} \right),\\
			\theta_{m} &= \frac{\pi}{2} - \frac{\frac{2}{3} - 8v}{\rho} + \frac{\frac{2}{9}-8v}{\rho^{2}} + \frac{\frac{256v^3}{3}+\frac{32v^2}{3}-\frac{92v}{9}-64f_1-\frac{94}{81}}{\rho^{3}} + \mathcal{O}\left( \frac{1}{\rho^{4}} \right),
		\end{split}
	\end{equation}
    }
    
    Now we have successfully changed back to the $(r,\theta)$ coordinates. Using $(\ref{map_F^2})$ to iterate once more, we have
	
	\begin{equation}
		\begin{split}
			r_{m + 1}& = \frac{1}{2}\rho - \frac{1}{2} + \frac{8v^2-\frac{88v}{3}+\frac{659}{36}}{\rho} + \frac{8v^2 -\frac{256v}{9}-2e_1 +\frac{1801}{108}}{\rho^{2}}\mathcal{+ O}\left( \frac{1}{\rho^{3}} \right),\\
			\theta_{m + 1} &= \frac{\pi}{2} - \frac{\frac{26}{3} - 8v}{\rho} + \frac{\frac{74}{9}-8v}{\rho^{2}} + \frac{\frac{256v^3 }{3}-\frac{1120v^2 }{3}+\frac{4516v}{9}-64f_1 -\frac{16402}{81}}{\rho^{3}} + \mathcal{O}\left( \frac{1}{\rho^{4}} \right),
		\end{split}
	\end{equation}
	Change the coordinates back to \((\rho\), \(\psi)\), using the Taylor expansion $(\ref{Taylorrho})$ and $(\ref{Taylorpsi})$ in the region $\Rtwo$, we have: \begin{equation}
		\begin{split}
			\rho_{m+1} &= \rho + \frac{1}{\rho}\left(\frac{22}{3}-16v\right) + \frac{1}{\rho^{2}}\left(4e_2-4e_1-16v+\frac{58}{9}\right)\mathcal{+ O}\left( \frac{1}{\rho^{3}} \right),\\
			\psi_{m+1}& = \frac{1}{3} + \frac{1}{\rho}\left(\frac{7}{6} - v\right) + \frac{1}{\rho^{2}}\left(\frac{\pi}{2}-\frac{7}{9}\right) \\
            &+ \frac{1}{\rho^{3}}\left(-8v^2+\frac{74v}{3}-8f_1+8f_2+\pi-\frac{601}{36}\right) + \mathcal{O}\left( \frac{1}{\rho^{4}} \right),
		\end{split}
	\end{equation}
    Since in the region $\mathcal{D}_1$, $\phi=\rho(\psi-1/3)$, we have
	\begin{equation}\label{F1}
		\begin{split}
        \mathcal{F}_1:(\rho,&\psi)\mapsto(\rho',\phi')\\
			\rho' &= \rho + \frac{1}{\rho}\left(\frac{22}{3}-16v\right) + \frac{1}{\rho^{2}}\left(4e_2-4e_1-16v+\frac{58}{9}\right)\mathcal{+ O}\left( \frac{1}{\rho^{3}} \right),\\
			\phi'& = \left(\frac{7}{6} - v\right) + \frac{1}{\rho}\left(\frac{\pi}{2}-\frac{7}{9}\right) \\
            &+ \frac{1}{\rho^{2}}\left(8v^2 -\frac{4v}{3}-8f_1 +8f_2 +\pi -\frac{293}{36}\right) + \mathcal{O}\left( \frac{1}{\rho^{3}} \right),
		\end{split}
	\end{equation}
       where $v=\left\{\rho/3-\phi\right\}$. From $(\ref{initialrhophi})$ and $(\ref{initialv})$, we know that for $(\rho',\phi')\in\mathcal{F}_1(\mathbf{R}_n)$,
       \begin{equation}\label{region1rhopsi}
       \begin{split}
           \rho'&\in\left(3n+\frac{1}{4}-\frac{9}{1024},3n+\frac{1}{4}+\frac{9}{1024}\right)+\mathcal{O}\left(\frac{1}{n}\right),\\
           \phi'&\in\left(\frac{7}{12}-\frac{9}{2048}\right),\left(\frac{7}{12}+\frac{9}{2048}\right)+\mathcal{O}\left(\frac{1}{n}\right)
        \end{split}
       \end{equation}
       Let $(\rho^{(1)},\phi^{(1)})=\mathcal{F}_1(\rho^{(0)},\phi^{(0)})$, then by $\ref{initialv0}$
       \begin{equation}\label{fixpoint 1}
           (\rho^{(1)},\phi^{(1)})=\left(3n+\frac{1}{4},\frac{7}{12}\right)+\mathcal{O}\left(\frac{1}{n}\right).
       \end{equation}

{\color{red} We now proceed to study $\mathcal{F}_2\sim\mathcal{F}_4$, and calculations similar to what we have done in the first step are omitted. See the supplementary material for details.}

\emph{Step} 2. \emph{The map $\mathcal{F}_2$}.

Then we study the map $\mathcal{F}_2$. Still, we start from a point $(\rho,\phi)$ in $\mathcal{F}_1(\mathbf{R}_n)$ in the region $\Rtwo$. Recall that in the region $\mathcal{D}_2$, $\phi=\rho(\psi-1/3)$. It takes
        $$m=\left[\frac{2\rho}{3}-\rho\psi\right]=\left[\frac{\rho}{3}-\phi\right]$$
        steps to leave the region $\Rtwo$. Similarly, letting $(\rho_m,\psi_m)=F^{2m}(\rho,\psi)$, then we have
    \begin{equation}\label{cross2}
		\begin{split}
			\rho_{m} &= \rho + \mathcal{O}\left( \frac{1}{\rho^{3}} \right),\\
			\psi_{m} &= \frac{1}{3} - \frac{v}{\rho} + \mathcal{O}\left( \frac{1}{\rho^{4}} \right),
		\end{split}
	\end{equation}
    where 
    \begin{equation}\label{region1v}
         v=\left\{\frac{\rho}{3}-\phi\right\}\in\left(\frac{1}{2}-\frac{15}{2048},\frac{1}{2}+\frac{15}{2048}\right)+\mathcal{O}\left(\frac{1}{n}\right)
    \end{equation}
    by $(\ref{region1rhopsi})$. And 
    \begin{equation}\label{region1v0}
         v^{(1)}=\left\{\frac{\rho^{(1)}}{3}-\phi^{(1)}\right\}=\frac{1}{2}+\mathcal{O}\left(\frac{1}{n}\right).
    \end{equation}
    A similar computation as in the first step yields
\begin{equation}\label{F2}
		\begin{split}
        \mathcal{F}_2:(\rho,&\phi)\mapsto(\rho',\phi')\\
			\rho' &= \rho -4v + \frac{1}{\rho^{2}}\left(\frac{v{\left(8v^2 -12v+11\right)}}{3}\right)\mathcal{+ O}\left( \frac{1}{\rho^{3}} \right),\\
			\phi'& =1- v + \frac{1}{\rho^{2}}\left(\frac{v}{2}\right)+
            \mathcal{O}\left( \frac{1}{\rho^{3}} \right),
		\end{split}
	\end{equation}
       where $v=\left\{\rho/3-\phi\right\}$. Moreover, from $(\ref{region1rhopsi})$ and $(\ref{region1v})$, we know that for $(\rho',\phi')\in\mathcal{F}_2\mathcal{F}_1(\mathbf{R}_n)$,
       \begin{equation}\label{region2rhopsi}
       \begin{split}
           \rho'&\in\left(3n-\frac{7}{4}-\frac{15}{512},3n-\frac{7}{4}+\frac{15}{512}\right)+\mathcal{O}\left(\frac{1}{n}\right),\\
           \phi'&\in\left(\frac{1}{2}-\frac{15}{2048},\frac{1}{2}+\frac{15}{2048}\right)+\mathcal{O}\left(\frac{1}{n}\right).
        \end{split}
       \end{equation}
        Let $(\rho^{(2)},\phi^{(2)})=\mathcal{F}_2(\rho^{(1)},\phi^{(1)})$, then by $(\ref{region1v0})$
       \begin{equation}\label{fixpoint 2}
           (\rho^{(2)},\phi^{(2)})=\left(3n-\frac{7}{4},\frac{1}{2}\right)+\mathcal{O}\left(\frac{1}{n}\right)
       \end{equation}
       
\emph{Step} 3. \emph{The map $\mathcal{F}_3$}.

Then we study the map $\mathcal{F}_3$. We start from $(\rho,\phi)$ in $\mathcal{F}_2\mathcal{F}_1(\mathbf{R}_n)$ in the region $\Rthree$, then it takes $m$ steps before the point leaving the region $\Rthree$, where
\begin{equation}
		m = \left\lbrack \frac{\rho}{3} - \rho\psi \right\rbrack=\left[\frac{\rho}{3}-\phi\right].
	\end{equation}
    And
\begin{equation}\label{cross3}
		\begin{split}
			\rho_{m} &= \rho + \mathcal{O}\left( \frac{1}{\rho^{3}} \right),\\
			\psi_{m} &= \frac{1}{3} - \frac{v}{\rho} + \mathcal{O}\left( \frac{1}{\rho^{4}} \right).
		\end{split}
	\end{equation}
 where 
 \begin{equation}\label{region3v}
 v=\left\{\frac{\rho}{3}-\phi\right\}\in\left(\frac{11}{12}-\frac{35}{2048},\frac{11}{12}+\frac{35}{2048}\right)+\mathcal{O}(n^{-1})
 \end{equation}
 by $(\ref{region2rhopsi})$. And 
    \begin{equation}\label{region3v0}
         v^{(2)}=\left\{\frac{\rho^{(2)}}{3}-\phi^{(2)}\right\}=\frac{11}{12}+\mathcal{O}\left(\frac{1}{n}\right).
    \end{equation}
    Recall that 
       \begin{align}
		\ell_{2}^{'} \subset \left\{ \psi = \frac{1}{3} - \frac{5}{6\rho} + \mathcal{O}\left( \frac{1}{\rho^{2}} \right) \right\},
	\end{align}
 Then $(\rho_m,\psi_m)$ is still in region $\Rthree$.
 Change back to $(r,\theta)$ coordinates, and we have
  \begin{equation}\label{inverse Taylor 3}
		\begin{split}
			r_{m} &= \frac{1}{2}\rho +\frac{1}{2}+\frac{8v^2 -\frac{8v}{3}+\frac{23}{36}}{\rho}  \\
            &+ \frac{\frac{128v}{9}-2e_3+\frac{4\pi }{3}-8\pi v-8v^2 -\frac{277}{108}}{\rho^{2}}\mathcal{+ O}\left( \frac{1}{\rho^{3}} \right),\\
			\theta_{m} &=\frac{\pi}{2} - \frac{8v+\frac{2}{3}}{\rho} + \frac{8v+4\pi -\frac{46}{9}}{\rho^{2}} \\
            &+\frac{\frac{256v^3 }{3}-\frac{32v^2 }{3}-\frac{524v}{9}-64f_3 -12\pi +\frac{1174}{81}}{\rho^{3}} + \mathcal{O}\left( \frac{1}{\rho^{4}} \right),
		\end{split}
	\end{equation}
 Iterating once more, we obtain
 \begin{equation}
		\begin{split}
			r_{m + 1}& = \frac{1}{2}\rho + \frac{1}{2} + \frac{8v^2-\frac{56v}{3}+\frac{407}{36}}{\rho}\\ 
            &+ \frac{\frac{272v}{9}-2e_3 +\frac{28\pi }{3}-8\pi v-8v^2 -\frac{2677}{108}}{\rho^{2}}\mathcal{+ O}\left( \frac{1}{\rho^{3}} \right),\\
			\theta_{m + 1} &= \frac{\pi}{2} + \frac{\frac{22}{3} - 8v}{\rho} + \frac{8v+4\pi -\frac{118}{9}}{\rho^{2}} \\
            &+ \frac{\frac{256v^3 }{3}-\frac{800v^2 }{3}+\frac{1972v}{9}-64f_3 -12\pi -\frac{1886}{81}}{\rho^{3}} + \mathcal{O}\left( \frac{1}{\rho^{4}} \right).
		\end{split}
	\end{equation}
Notice that $r_{m+1} =n/2+\mathcal{O}(1)$. From $(\ref{region3v})$, we have
\begin{equation*}
    \theta_{m+1}\geq \frac{\pi}{2}-\frac{1}{\rho}\frac{35}{768}+\mathcal{O}\left(\frac{1}{n^2}\right),
\end{equation*}
then
\begin{equation*}
    \theta_{m+1}\geq \frac{\pi}{2}-\frac{1}{r_{m+1}}\frac{35}{384}+\mathcal{O}\left(\frac{1}{r_{m+1}^2}\right).
\end{equation*}
But we know that
\begin{align}
		\ell_{3} \subset \left\{ \theta=\frac{\pi}{2}-\frac{1}{r}+\mathcal{O}\left(\frac{1}{r^2}\right)\right\}.
	\end{align}
Then for $n$ large enough, $(r_{m+1},\theta_{m+1})\in\Rfour$, and the orbit jumps from $\Rthree$ to $\Rfour$, without going into the region $\Rsix$.
Changing back to the $(\rho,\psi)$ coordinates, we have
\begin{equation}
\begin{split}
			\rho_{m+1} &= \rho + \frac{1}{\rho}\left(\frac{50}{3}-16v\right) + \frac{1}{\rho^{2}}\left(4e_4-4e_3+8\pi+16v-\frac{254}{9}\right)\mathcal{+ O}\left( \frac{1}{\rho^{3}} \right),\\
			\psi_{m+1}& = \frac{1}{3} + \frac{1}{\rho}\left(\frac{5}{6} - v\right) + \frac{1}{\rho^{2}}\left(\frac{\pi}{2}-\frac{7}{9}\right) \\
            &+ \frac{1}{\rho^{3}}\left(-8v^2+\frac{22v}{3}-8f_3+8f_4-\pi+\frac{73}{36}\right) + \mathcal{O}\left( \frac{1}{\rho^{4}} \right),
		\end{split}
	\end{equation}
    Recall that in the region $\mathcal{D}_3$, $\phi=\rho(\psi-1/3)$. then we have
       \begin{equation}\label{F3}
		\begin{split}
        \mathcal{F}_3:(\rho,&\phi)\mapsto(\rho',\phi')\\
			\rho' &= \rho + \frac{1}{\rho}\left(\frac{50}{3}-16v\right) + \frac{1}{\rho^{2}}\left(4e_4-4e_3+8\pi+16v-\frac{254}{9}\right)\mathcal{+ O}\left( \frac{1}{\rho^{3}} \right),\\
			\phi'& = \frac{5}{6} - v + \frac{1}{\rho}\left(\frac{\pi}{2}-\frac{7}{9}\right) \\
            &+ \frac{1}{\rho^{2}}\left(8v^2 -\frac{20v}{3}-8f_3 +8f_4 -\pi +\frac{31}{12}\right) + \mathcal{O}\left( \frac{1}{\rho^{3}} \right),
		\end{split}
	\end{equation}
       where $v=\left\{\rho/3-\phi\right\}$. Moreover, from $(\ref{region2rhopsi})$ and $(\ref{region3v})$, we know that for $(\rho',\phi')\in\mathcal{F}_3\mathcal{F}_2\mathcal{F}_1(\mathbf{R}_n)$,
       \begin{equation}\label{region3rhopsi}
       \begin{split}
           \rho'&\in\left(3n-\frac{7}{4}-\frac{15}{512},3n-\frac{7}{4}+\frac{15}{512}\right)+\mathcal{O}\left(\frac{1}{n}\right),\\
             \phi'&=\rho'(\psi'-1/3)\in\left(-\frac{1}{12}-\frac{35}{2048},-\frac{1}{12}+\frac{35}{2048}\right)+\mathcal{O}\left(\frac{1}{n}\right).
        \end{split}
       \end{equation}
        Let $(\rho^{(3)},\phi^{(3)})=\mathcal{F}_3(\rho^{(2)},\phi^{(2)})$, then by $(\ref{region3v0})$
       \begin{equation}\label{fixpoint 3}
           (\rho^{(3)},\phi^{(3)})=\left(3n-\frac{7}{4},-\frac{1}{12}\right)+\mathcal{O}\left(\frac{1}{n}\right)
       \end{equation}
       
       \emph{Step} 4. \emph{The map $\mathcal{F}_4$}.
       
       A point $(\rho,\phi)$ given by $(\ref{region3rhopsi})$ in the region $\Rfour$ takes 
       $$m=\left[\frac{2\rho}{3}-\rho\psi\right]=\left[\frac{\rho}{3}-\phi\right]$$
       steps before it enters the region $\Rone$. Let
       $$v=\left\{\frac{2\rho}{3}-\rho\psi\right\},$$
    and
    \begin{equation}\label{region4v0}
         v^{(3)}=\left\{\frac{\rho^{(3)}}{3}-\phi^{(3)}\right\}=\frac{1}{2}+\mathcal{O}\left(\frac{1}{n}\right).
    \end{equation}
    Then after $m$ times interation of $F^2$,
       \begin{equation}\label{cross4}
		\begin{split}
			\rho_{m} &= \rho + \mathcal{O}\left( \frac{1}{\rho^{3}} \right),\\
			\psi_{m} &= \frac{2}{3} - \frac{v}{\rho} + \mathcal{O}\left( \frac{1}{\rho^{4}} \right).
		\end{split}
	\end{equation}
       Still, similar calculation like the first step shows that:
       \begin{equation}\label{F4}
		\begin{split}
        \mathcal{F}_4:(\rho,&\phi)\mapsto(\rho',\phi')\\
			\rho' &= \rho +4-4v+ \frac{1}{\rho^{2}}\frac{8v^3-4v^2+5v-1}{3}\mathcal{+ O}\left( \frac{1}{\rho^{3}} \right),\\
			\phi'& =1 - v+ \frac{1}{\rho^{2}}\left(\frac{1}{2}-\frac{v}{2}\right) + \mathcal{O}\left( \frac{1}{\rho^{3}} \right),
		\end{split}
	\end{equation}
       where $v=\left\{\rho/3-\phi\right\}$.
     Let $(\rho^{(4)},\phi^{(4)})=\mathcal{F}_4(\rho^{(3)},\phi^{(3)})$, then by $(\ref{region4v0})$
       \begin{equation}\label{fixpoint 4}
           (\rho^{(4)},\phi^{(4)})=\left(3n+\frac{1}{4},\frac{1}{2}\right)+\mathcal{O}\left(\frac{1}{n}\right)
       \end{equation}
    Notice that
    \begin{equation}
        (\rho^{(4)},\phi^{(4)})=\mathcal{F}((\rho^{(0)},\phi^{(0)}))=(\rho^{(0)},\phi^{(0)})+\mathcal{O}\left(\frac{1}{n}\right).
    \end{equation}

\emph{Step} 5. \emph{Expand near the fixed point}.
   
   Define 
   \begin{equation}\label{rhophitoxy}
   \begin{split}
        \rho&=3n+\frac{1}{4}+\tilde{x},\\
        \phi&=1/2+\tilde{y}
        \end{split}
   \end{equation}
   $\rho=3n+1/4+\tilde{x},\phi=1/2+\tilde{y}$, then in the coordinates $(\tilde{x},\tilde{y})$, we have
   \begin{align}
   &\mathcal{F}:(\tilde{x},\tilde{y})\mapsto(\tilde{x}',\tilde{y}'),\notag\\
   &(\tilde{x}',\tilde{y}')=\mathcal{F}_4\mathcal{F}_3\mathcal{F}_2\mathcal{F}_1(3n+1/4+\tilde{x},1/2+\tilde{y})-(3n+1/4,1/2).
   \end{align}
   We calculate the composition $\mathcal{F}_4\mathcal{F}_3\mathcal{F}_2\mathcal{F}_1$ using $(\ref{F1})$, $(\ref{F2})$, $(\ref{F3})$ and $(\ref{F4})$ and expand it into powers of $1/n$, then we obtain:
   \begin{equation}\label{Ultimate1}
   \begin{split}
       \tilde{x}'&=\frac{1}{9}\tilde{x}-\frac{8}{3}\tilde{y}+\frac{1}{n}\left(\frac{2\pi }{9}-\frac{4}{81}-\frac{32\tilde{y}}{9}+\frac{128\tilde{x}}{81}\right)+\frac{1}{n^2}X+\mathcal{O}\left(\frac{1}{n^3}\right),\\
       \tilde{y}'&=\frac{4}{9}\tilde{x}-\frac{5}{3}\tilde{y}+\frac{1}{n}\left(-\frac{\pi }{9}+\frac{2}{81}+\frac{32\tilde{x}}{81}\right)+\frac{1}{n^2}Y+\mathcal{O}\left(\frac{1}{n^3}\right),
       \end{split}
   \end{equation}
   where 
    \begin{equation}
   \begin{split}
       X&=-\frac{1520\tilde{x}\tilde{y}^2 }{729}+\frac{2080\tilde{x}^2 \tilde{y}}{2187}-\frac{3392\tilde{x}^3 }{19683}+\frac{1360\tilde{y}^3 }{729}\\
       &-\frac{416\tilde{x}^2 }{729}-\frac{224\tilde{y}^2 }{81}-\frac{1922\tilde{x}}{729}+\frac{998\tilde{y}}{243}+\frac{14\pi \tilde{x}}{27}-\frac{8\pi \tilde{y}}{9}+\frac{800\tilde{x}\tilde{y}}{243}\\
&+\frac{28e_1 }{81}-\frac{28e_2 }{81}+\frac{4e_3 }{27}-\frac{4e_4 }{27}+\frac{64f_1 }{27}-\frac{64f_2 }{27}-\frac{32f_3 }{9}+\frac{32f_4 }{9}+\frac{1486}{729}-\frac{83\pi }{162},\\ 
       Y&=-\frac{32\tilde{x}\tilde{y}^2 }{81}+\frac{64\tilde{x}^2 \tilde{y}}{243}-\frac{128\tilde{x}^3 }{2187}+\frac{16\tilde{y}^3 }{81}\\
       &-\frac{80\tilde{x}^2 }{729}-\frac{176\tilde{y}^2 }{81}-\frac{1205\tilde{x}}{729}+\frac{950\tilde{y}}{243}+\frac{5\pi \tilde{x}}{27}-\frac{2\pi \tilde{y}}{9}+\frac{320\tilde{x}\tilde{y}}{243}\\
       &+\frac{4e_1 }{81}-\frac{4e_2 }{81}+\frac{4e_3 }{27}-\frac{4e_4 }{27}+\frac{40f_1 }{27}-\frac{40f_2 }{27}-\frac{8f_3 }{9}+\frac{8f_4 }{9}+\frac{817}{729}-\frac{121\pi }{324}.
       \end{split}
   \end{equation}
Now $(0,0)$ is a nondegenerate fixed point of the linear part of $(\ref{Ultimate1})$, then there exists a fixed point $(x_0,y_0)=\mathcal{O}(n^{-1})$ of $(\ref{Ultimate1})$. Then we define
\begin{equation}\label{findH1}
        (x_1,y_1)=(\tilde{x},\tilde{y})-(x_0,y_0).
\end{equation}
and the conjugacy 
\begin{align}
H_1:&(x,y)\mapsto(x_1,y_1),
\end{align}
Notice that this conjugacy is found by a sequence of coordinate changes:
$$(x,y)\rightarrow(\rho,\phi)\rightarrow(\tilde{x},\tilde{y})\rightarrow(x_1,y_1),$$
Then from $(\ref{xytorhophi})$ and $(\ref{rhophitoxy})$, we know that 
\begin{equation}
    \begin{split}
        x_1&=x+y-3n-\frac{1}{4}+\mathcal{O}\left(\frac{1}{n}\right),\\
        y_1&=\frac{1}{2}y-\frac{1}{2}+\mathcal{O}\left(\frac{1}{n}\right).
    \end{split}
\end{equation}
Plug $(\ref{findH1})$ into $(\ref{Ultimate1})$  and notice that only polynomials of degree at most 3 appear in $(\ref{Ultimate1})$. Then we have
    \begin{equation}
		\begin{split}
			H_1\circ \mathcal{F}\circ H_1^{-1}:&\left( x_{1},y_{1} \right)\mapsto ({\widehat{x}}_{1},{\widehat{y}}_{1})\\
			{\widehat{x}}_{1} = A_{1,1}x_{1} +& A_{1,2}y_{1} + \sum_{j = 2}^{3}{F_{1,j}\left( x_{1},y_{1} \right)}+\mathcal{O}_4\left( n^{- 3} \right),\\
			{\widehat{y}}_{1} = A_{2,1}x_{1} + &A_{2,2}y_{1} + \sum_{j = 2}^{3}{F_{2,j}\left( x_{1},y_{1} \right)}+  \mathcal{O}_4\left( n^{- 3} \right),\end{split}
	\end{equation}
    where $F_{i,j}$ is a homogeneous polynomial of degree $j$ for any $i,j$.

    From $(\ref{Ultimate1})$, we know that
    \[A_{1,1} = \frac{1}{9}+\mathcal{ O}\left( n^{- 1} \right),\ \ A_{1,2} = - \frac{16}{9}+\mathcal{O}\left( n^{- 1} \right),\]
\[A_{1,1} = \frac{2}{3}+\mathcal{ O}\left( n^{- 1} \right),\ \ A_{2,2} = - \frac{5}{3}+\mathcal{O}\left( n^{- 1} \right).\]
And since the outer billiard map $F$ is area-reserving, the absolute value of the Jacobian determinant of $\mathcal{F}$ is $1$ at the fixed point. This property is preserved under any conjugacy, thus ${\rm det}(A_{ij})=1$.  

Then we have the following.
\begin{align*}
      F_{i,2}&=0+\mathcal{O}\left(n^{-2}\right),\ \ i=1,2,
      \end{align*}
because all terms higher than degree 2 are of $\mathcal{O}(n^{-2})$ in $(\ref{Ultimate1})$. And 
\begin{align*}
      F_{1,3}&=-\frac{1520xy^2 }{729}+\frac{2080{x}^2 {y}}{2187}-\frac{3392{x}^3 }{19683}+\frac{1360{y}^3 }{729}+\mathcal{O}\left(n^{-3}\right),\\
      F_{2,3}&=-\frac{32{x}{y}^2 }{81}+\frac{64{x}^2 {y}}{243}-\frac{128{x}^3 }{2187}+\frac{16{y}^3 }{81}+\mathcal{O}\left(n^{-3}\right),
\end{align*}
because no polynomials of degree higher than 3 appear in $(\ref{Ultimate1})$, then the linear transform $(\tilde{x},\tilde{y})\mapsto (x,y)$ will not change the monomials of degree 3.

Then Proposition B is proved.
\end{proof}
\begin{remark}
    It can be observed that $e_1\sim e_4$ and $f_1\sim f_4$ appear only as constant terms in $(\ref{Ultimate1})$, and that's why it doesn't influence the $F_{1,3}$ and $F_{2,3}$ up to $\mathcal{O}(n^{-3})$.
\end{remark}

\section{Proof of Lemma \ref{keylemma}}\label{proof of key lemma}
The process of finding the normal form for a circular outer billiard system is similar to the process we have done to find the normal form for the semi-disc outer billiard system. There are 3 parts: find the expression of $F_\beta^2$ in the polar coordinate, find a new coordinate to simplify the dynamics in each continuous region of $F_\beta^2$ and calculate the first return map.
\subsection{Circular sector outer billiard}
In this section, we calculate the asymptotic expression of circular sector outer billiard systems to prepare ourselves for finding the normal form.

Here we use the coordinates in which the circular sector is given by
$$\left\{ x^2+y^2\leq1, y\geq -\cos\beta \right\}.$$
$F_\beta$ has discontinuities at $\ell_1=\left\{x\geq \sin\beta ,y=-\cos\beta \right\}$ and the half tangent line $\ell_2$, $\ell_3$ at the points $\left(\sin\beta, -\cos\beta\right)$ and $\left(-\sin\beta, -\cos\beta \right)$. Away from the discontinuities, the map $F_\beta$ has the following forms:

between $\ell_1$ and $\ell_2$ is the reflection about point $O_1=\left(\sin\beta, -\cos\beta\right)$, denote this reflection by $R_1$;

between $\ell_2$ and $\ell_3$ is the reflection about a tangency point to the circular part, denote this reflection by $T$;

between $\ell_3$ and $\ell_1$ is the reflection about point $O_2=\left(-\sin\beta, -\cos\beta\right)$, denote this reflection by $R_2$.

Let $\ell'_j={F_\beta}^{-1}\ell_j$. $F_\beta^2$ will have six continuity domains.

Region $\uppercase\expandafter{\romannumeral1}$: $F_\beta^2=TR_1$;

Region $\uppercase\expandafter{\romannumeral2}$: $F_\beta^2=R_2T$;

Region $\uppercase\expandafter{\romannumeral3}$: $F_\beta^2=R_1T$;

Region $\uppercase\expandafter{\romannumeral4}$: $F_\beta^2=TR_2$;

Region $\uppercase\expandafter{\romannumeral5}$ and $\Rsix$: $F_\beta^2=T^2$;

For the upper region $\uppercase\expandafter{\romannumeral1}$,  $\uppercase\expandafter{\romannumeral2}$ and $\uppercase\expandafter{\romannumeral5}$  polar coordinate $\left(r,\theta\right)$ is considered, and in the lower region $\uppercase\expandafter{\romannumeral3}$,  $\uppercase\expandafter{\romannumeral4}$ and $\uppercase\expandafter{\romannumeral6}$ polar coordinate shifted by $\pi$ is considered. Therefore, both in the upper and lower region $\theta \in \left[0,\pi\right]$

\begin{proposition}\label{prop3}
	The discontinuous lines are given by the following formula:
	\begin{align*}
		&\ell_1\subset \left\{\theta=-\frac{\cos\beta}{r}+\mathcal{O}\left(\frac{1}{r^2}\right)\right\}, \ell_1'\subset \left\{\theta=\pi-\frac{2+\cos\beta}{r}+\mathcal{O}\left(\frac{1}{r^2}\right)\right\}, \\
		&\ell_2\subset \left\{\theta=\beta-\frac{1}{r}+\mathcal{O}\left(\frac{1}{r^2}\right)\right\}, \ell_2'\subset \left\{\theta=\pi-\beta-\frac{3}{r}+\mathcal{O}\left(\frac{1}{r^2}\right)\right\}, \\
		&\ell_3\subset \left\{\theta=\pi-\beta-\frac{1}{r}+\mathcal{O}\left(\frac{1}{r^2}\right)\right\}, \ell_3'\subset \left\{\theta=\beta-\frac{3}{r}+\mathcal{O}\left(\frac{1}{r^2}\right)\right\}, \\
	\end{align*}
\end{proposition}

\begin{proof}
	Straightforward calculation. It should be noticed that lines $\ell_1$ and $\ell_1'$ may have $\theta$ shifted by $\pi$ when the equation is written in the lower region coordinate.
\end{proof}

\begin{proposition}
	In the coordinate we introduce, $R_1$, $R_2$ and $T$ have the following forms:
	\begin{align*}
		R_1\left(r, \theta \right):&\\
		&r\mapsto  r-2\sin\left(\beta-\theta\right)+\frac{1+\cos\left(2\beta-2\theta\right)}{r}+\mathcal{O}\left(\frac{1}{r^2}\right),\\
		&\theta\mapsto\theta+\frac{2\cos\left(\beta-\theta\right)}{r}+\frac{2\sin\left(2\beta-2\theta\right)}{r^2}+\mathcal{O}\left(\frac{1}{r^3}\right),\\		
		R_2\left(r, \theta \right):&\\
		&r\mapsto  r-2\sin\left(\beta+\theta\right)+\frac{1+\cos\left(2\beta+2\theta\right)}{r}+\mathcal{O}\left(\frac{1}{r^2}\right),\\
		&\theta\mapsto\theta+\frac{2\cos\left(\beta+\theta\right)}{r}+\frac{2\sin\left(2\beta+2\theta\right)}{r^2}+\mathcal{O}\left(\frac{1}{r^3}\right),		\\
		T\left(r, \theta \right):&\\
		&r\mapsto r\\
		&\theta\mapsto\theta+\frac{2}{r}+\mathcal{O}\left(\frac{1}{r^3}\right),		
	\end{align*}
\end{proposition}

\begin{proof}
	Take $R_1$ for example. $R_1\left(x,y\right)=\left(-x+2\cos\beta,-y-2\sin\beta\right)$,  mapping points in upper region to the lower region. Change to our coordinates, we will obtain the required form. The other two maps are similar.
\end{proof}

\begin{proposition}\label{mappolar}
	The map $F_\beta^2$ takes the following form
	$$F_\beta^2\left(r,\theta\right)=\left(r+a\left(\theta\right)+\frac{a_1\left(\theta\right)}{r}+\mathcal{O}\left(\frac{1}{r^2}\right),\theta+\frac{b\left(\theta\right)}{r}+\frac{b_1\left(\theta\right)}{r^2}+\mathcal{O}\left(\frac{1}{r^3}\right)\right),$$
	where in region $\Rone$
	\begin{align*}
		&a(\theta)=-2\sin(\beta-\theta),\ a_1(\theta)=1+\cos(2\beta-2\theta),\\
		&b(\theta)=2+2\cos(\beta-\theta),\ b_1(\theta)=4\sin(\beta-\theta)+2\sin(2\beta-2\theta).
	\end{align*}
	In region $\Rtwo$
	\begin{align*}
		&a(\theta)=-2\sin(\beta+\theta),\ a_1(\theta)=1+\cos(2\beta+2\theta)-4\cos(\theta+\beta),\\
		&b(\theta)=2-2\cos(\beta+\theta),\ b_1(\theta)=4\sin(\beta+\theta)-2\sin(2\beta+2\theta).
	\end{align*}
	In region $\Rthree$
	\begin{align*}
		&a(\theta)=-2\sin(\beta-\theta),\ a_1(\theta)=1+\cos(2\beta-2\theta)+4\cos(\beta-\theta),\\
		&b(\theta)=2+2\cos(\beta-\theta),\ b_1(\theta)=4\sin(\beta-\theta)+2\sin(2\beta-2\theta).
	\end{align*}
	In region $\Rfour$
	\begin{align*}
		&a(\theta)=-2\sin(\beta+\theta),\ a_1(\theta)=1+\cos(2\beta+2\theta),\\
		&b(\theta)=2-2\cos(\beta+\theta),\ b_1(\theta)=4\sin(\beta+\theta)-2\sin(2\beta+2\theta).
	\end{align*}
	In region $\Rfive$ and $\Rsix$
	\begin{align*}
		a(\theta)=0=a_1(\theta),\ b(\theta)=4,\ b_1(\theta)=0.
	\end{align*}
\end{proposition}

\begin{proof}
	Direct computation from the previous proposition. It should be noticed that in the region $\Rfive$ and the region $\Rsix$, $F_\beta^2$ is integrable and $r$ is conserved.
\end{proof}

\subsection{Coordinate to simplify the dynamics of \texorpdfstring{$F_\beta^2$}{TEXT}}
We introduce a coordinate change in region $\Rone$,$\Rtwo$,$\Rthree$ and $\Rfour$ to better describe the dynamics in these regions. 

\begin{proposition}\label{Finding rho and psi}
	There exists a smooth change of coordinates $(r,\theta)\mapsto(\rho,\psi)$ in the region $\Rone$,$\Rtwo$,$\Rthree$ and $\Rfour$, so that in these regions $F_\beta^2$ takes the following normal form in variables $\left(\rho,\psi\right)$. If $F^2_\beta$ takes $(r,\theta)$ to $(r',\theta')$ and $(\rho,\psi)$ to $(\rho',\psi')$, then
	\begin{align}
		&\rho'=\rho+\mathcal{O}\left(\frac{1}{\rho^2}\right),\\
		&\psi'=\psi+\frac{1}{\rho}+\mathcal{O}\left(\frac{1}{\rho^3}\right).
	\end{align}
	The coordinate transformation has the following form
	\begin{align*}
		\rho=r\Phi_1(\theta)+\Phi_2(\theta),\\
		\psi=\Psi(\theta)+\frac{\Psi_1(\theta)}{r},
	\end{align*}
	And in region $\Rone$
	\begin{equation}\label{explicit_form1}
		\begin{split}
			& \Phi_1(\theta)=\cos^2\left(\frac{\beta}{2}-\frac{\theta}{2}\right), \ \Phi_2(\theta)=-\frac{\sin(\beta-\theta)}{2},\\
			&\Psi(\theta)=-\frac{1}{6}\tan^3\left(\frac{\beta}{2}-\frac{\theta}{2}\right)-\frac{1}{2}\tan\left(\frac{\beta}{2}-\frac{\theta}{2}\right),\\
                &\Psi_1(\theta)=-\frac{\cos \left(\frac{\beta }{2}-\frac{\theta}{2}\right)-\cos \left(\frac{3\,\beta }{2}-\frac{3\,\theta}{2}\right)}{16\,{\cos \left(\frac{\beta }{2}-\frac{\theta}{2}\right)}^5 }.
		\end{split}
	\end{equation}
        In region $\Rtwo$
	\begin{equation}\label{explicit_form2}
		\begin{split}
			& \Phi_1(\theta)=\sin^2\left(\frac{\beta}{2}+\frac{\theta}{2}\right), \ \Phi_2(\theta)=\frac{\sin(\beta+\theta)}{2},\\
			&\Psi(\theta)=-\frac{1}{6}\cot^3\left(\frac{\beta}{2}+\frac{\theta}{2}\right)-\frac{1}{2}\cot\left(\frac{\beta}{2}+\frac{\theta}{2}\right),\\
                &\Psi_1(\theta)=\frac{\cos \left(\beta +\theta\right)+1}{\cos \left(2\,\beta +2\,\theta\right)-4\,\cos \left(\beta +\theta\right)+3}.
		\end{split}
	\end{equation}
        In region $\Rthree$
	\begin{equation}\label{explicit_form3}
		\begin{split}
			& \Phi_1(\theta)=\cos^2\left(\frac{\beta}{2}-\frac{\theta}{2}\right), \ \Phi_2(\theta)=\frac{\sin(\beta-\theta)}{2},\\
			&\Psi(\theta)=-\frac{1}{6}\tan^3\left(\frac{\beta}{2}-\frac{\theta}{2}\right)-\frac{1}{2}\tan\left(\frac{\beta}{2}-\frac{\theta}{2}\right),\\
                &\Psi_1(\theta)=\frac{\cos \left(\frac{\beta }{2}-\frac{\theta}{2}\right)-\cos \left(\frac{3\,\beta }{2}-\frac{3\,\theta}{2}\right)}{16\,{\cos \left(\frac{\beta }{2}-\frac{\theta}{2}\right)}^5 }.
		\end{split}
	\end{equation}
        In region $\Rfour$
	\begin{equation}\label{explicit_form4}
		\begin{split}
			& \Phi_1(\theta)=\sin^2\left(\frac{\beta}{2}+\frac{\theta}{2}\right), \ \Phi_2(\theta)=-\frac{\sin(\beta+\theta)}{2},\\
			&\Psi(\theta)=-\frac{1}{6}\cot^3\left(\frac{\beta}{2}+\frac{\theta}{2}\right)-\frac{1}{2}\cot\left(\frac{\beta}{2}+\frac{\theta}{2}\right),\\
                &\Psi_1(\theta)=-\frac{\cos \left(\beta +\theta\right)+1}{\cos \left(2\,\beta +2\,\theta\right)-4\,\cos \left(\beta +\theta\right)+3}.
		\end{split}
	\end{equation}
\end{proposition}

\begin{proof}
        This proposition generalizes Lemma \ref{rho_psi}. The method of proving is the same. 
	To find $\Phi_1$ we expand $\rho$ in powers of r:
	$${\rho'} - \rho = \Phi_{1}^{'}b + \Phi_{1}a + \left( \Phi_{1}^{'}b_{1} + \frac{\Phi_{1}^{''}b^{2}}{2} + a\Phi_{1}^{'}b + \Phi_{1}a_{1} + \Phi_{2}^{'}b \right)\frac{1}{r}\mathcal{+ O}\left( \frac{1}{r^{2}} \right),$$
	where $a,b,a_1,b_1$ are derived in proposition \ref*{mappolar}. If we require that
	\begin{align}
		\frac{\Phi'_1}{\Phi_1}=-\frac{a}{b},
	\end{align}
	\begin{align}\label{Phi2}
		\Phi'_2=-\frac{1}{b}\left(a_1\Phi_1+\Phi'_1(ab+b_1)+\frac{\Phi''_1b^2}{2}\right),
	\end{align}
	we will obtain that $\rho'-\rho=\mathcal{O}(1/r^2)$. We further give the initial value that $\Phi_1(\beta)=1,\ \Phi_2(\beta)=0$ at region $\Rone, \Rthree$, and $\Phi_1(\pi-\beta)=1,\ \Phi_2(\pi-\beta)=0$ at region $\Rtwo, \Rfour$. This means that when an orbit goes from the integrable region $\Rfive,\ \Rsix$ to the regions $\Rone,\ \Rtwo,\ \Rthree,\ \Rfour$, $r$ can be transformed into $\rho$ continuously. For example, in the region $\Rone$, we have $a(\theta)=-2\sin(\beta-\theta),\ b(\theta)=2+2\cos(\beta-\theta)$, so
	\begin{align*}
		&\frac{\Phi'_1}{\Phi_1}=-\frac{-2\sin(\beta-\theta)}{2+2\cos(\beta-\theta)}, \Phi_1(\beta)=1,
	\end{align*}
	It can be observed that $\Phi(\theta)=\cos^2\left(\beta/2-\theta/2\right)$ is the solution. Similarly for $\beta_2$ we can solve\begin{align*}
		\Phi_2(\theta)=-\frac{\sin(\beta-\theta)}{2}.
	\end{align*}
	
	Likewise, we expand $\psi'-\psi$ and we will get
    \begin{equation}\label{equation to find Psi}
        \psi'-\psi=\Psi'\frac{b}{r}+\Psi'\frac{b_1}{r^2}+\frac{\Psi''}{2}{\left(\frac{b}{r}\right)}^2+\frac{\Psi'_1b}{r^2}-\frac{a\Psi_1}{r^2}+\mathcal{O}\left(\frac{1}{r^3}\right).
    \end{equation}
        Therefore, 
	$$\Psi'=\frac{1}{b\Phi_1}.$$
	Similarly, we give the initial value condition that $\Psi'(\beta)=0$ at region $\Rone, \Rthree$, and $\Psi'(\pi-\beta)=0$ at region $\Rtwo, \Rfour$. Leaving the term $1/r^2$ to be zero in equation \ref{equation to find Psi}, we obtain the equation for $\Psi_1$. We use the same initial value condition for $\Psi_1$ as for $\Psi$. Then we obtain exactly the solution in (\ref{explicit_form1})$\sim$(\ref{explicit_form4}).
    \end{proof}

\begin{remark}
	The coordinate transformation is designed to simplify the map in the region $\Rone \sim \Rfour$, where $F_\beta^2$ is not integrable. This proposition shows that $\Phi_i$ and $\Psi_i$ are monotonic functions, thus when $r$ is large enough, this transformation is a well-defined analytic homeomorphism. Using this coordinate transformation, the dynamics inside each region is simple, and significant changes occur when an orbit passes from one region to another. It should be mentioned that the initial value condition we use here is just chosen for simplicity.
\end{remark}

\begin{proposition}\label{disc_lines}
	The discontinuous lines in coordinate $(\rho, \beta)$ are given by the following equations, in region $\Rone$
	\begin{equation}\label{region1}
		\begin{split}
			&\ell_1 \subset\left\{\psi=-A_\beta-\frac{1}{4\rho}+\mathcal{O}\left(\frac{1}{\rho^2}\right)\right\},\\
			&\ell_{2}\subset\left\{\psi=-\frac{1}{4\rho}+\mathcal{O}\left(\frac{1}{\rho^2}\right)\right\}.
		\end{split}
	\end{equation}
	In region $\Rtwo$
	\begin{equation}\label{region2}
		\begin{split}
			&\ell_{3}'\subset\left\{\psi=-\frac{3}{4\rho}+\mathcal{O}\left(\frac{1}{\rho^2}\right)\right\},\\
			&\ell_1' \subset\left\{\psi=A_\beta-\frac{3}{4\rho}+\mathcal{O}\left(\frac{1}{\rho^2}\right)\right\}.
		\end{split}
	\end{equation}
	In region $\Rthree$
	\begin{equation}\label{region3}
		\begin{split}
			&\ell_1' \subset\left\{\psi=-A_\beta-\frac{3}{4\rho}+\mathcal{O}\left(\frac{1}{\rho^2}\right)\right\},\\
			&\ell_2'\subset\left\{\psi=-\frac{3}{4\rho}+\mathcal{O}\left(\frac{1}{\rho^2}\right)\right\}.
		\end{split}
	\end{equation}
	In region $\Rfour$
	\begin{equation}\label{region4}
		\begin{split}
			&\ell_3\subset\left\{\psi=-\frac{1}{4\rho}+\mathcal{O}\left(\frac{1}{\rho^2}\right)\right\},\\
			&\ell_1 \subset\left\{\psi=A_\beta-\frac{1}{4\rho}+\mathcal{O}\left(\frac{1}{\rho^2}\right)\right\},
		\end{split}
	\end{equation}
	where 
	$$A_\beta=\frac{1}{6}\tan^3\left(\frac{\beta}{2}\right)+\frac{1}{2}\tan\left(\frac{\beta}{2}\right)$$
\end{proposition}

\begin{proof}
	Straightforward calculation combining Proposition \ref{prop3} and Proposition \ref{Finding rho and psi}.
\end{proof}
\subsection{Proof of Lemma \ref{keylemma}}
We consider infinite region $\widetilde{\mathcal{D}}$ bounded by $\ell_1'$, $F_\beta^2(\ell_1')$ and $\left\{x=-\widetilde{x_0}\right\}$, then represent $\mathcal{F}_\beta$ as $\mathcal{F}_1\mathcal{F}_2$ where $\mathcal{F}_1$ is the passage from $\mathcal{D}$ to $\widetilde{\mathcal{D}}$ and $\mathcal{F}_1$ is the passage from $\widetilde{\mathcal{D}}$ to $\mathcal{D}$. The region $\mathcal{D}$ we are interested in lies in region $\Rone$, between $\ell_1$ and $F_\beta^2(\ell_1)$.Now we introduce the change of coordinates that make $\mathcal{D}$ diffeomorphic to a half cylinder up to the identification of its boundary lines $\ell_1$ and $F_\beta^2(\ell_1)$. 
$$\phi=\rho(\psi+A_\beta)+\frac{1}{4}+\mathcal{O}\left(\frac{1}{\rho^2}\right)$$
Therefore from (\ref{region1}) the  $(\rho,\psi)$ coordinates of a point in such region can be written as
\begin{align*}
	\left(\rho, -A_\beta+\frac{\phi}{\rho}-\frac{1}{4\rho}+\mathcal{O}\left(\frac{1}{\rho^2}\right)\right),
\end{align*}
where $A_\beta$ is the constant depending only on $\beta$ that we derive in proposition \ref*{disc_lines}, $\phi \in (0,1)$ and $\rho$ is large enough. Similarly, for region $\widetilde{\mathcal{D}}$ we consider
$$\widetilde{\phi}=\widetilde{\rho}(\psi+A_\beta)+\frac{3}{4}+\mathcal{O}\left(\frac{1}{\rho^2}\right)$$
then  $(\widetilde{\rho},\psi')$ coordinates of a point in such region can be written as
\begin{align*}
	\left(\widetilde{\rho}, -A_\beta+\frac{\widetilde{\phi}}{\rho}-\frac{3}{4\widetilde{\rho}}+\mathcal{O}\left(\frac{1}{\widetilde{\rho}^2}\right)\right),
\end{align*}
We now study the first passage map from $\mathcal{D}$  to region $\widetilde{\mathcal{D}}$ in coordinate $(\rho, \phi)$ and $(\widetilde{\rho},\widetilde{\phi})$

Until entering region $\Rfive$ we have 
\begin{align*}
	&\rho_k=\rho+\mathcal{O}\left(\frac{1}{\rho^2}\right),\\
	&\psi_k=-A_\beta-\frac{1}{4\rho}+\frac{\phi+k}{\rho}+\mathcal{O}\left(\frac{1}{\rho^2}\right).
\end{align*}
From proposition \ref*{disc_lines} equation (\ref*{region1}), we know in $(\rho,\psi)$ coordinates in region $\Rone$, 
$$\ell_{2}\subset\left\{\psi=-\frac{1}{4\rho}+\mathcal{O}\left(\frac{1}{\rho^2}\right)\right\}.$$
So after $n_1=\left[A_\beta\rho-\phi\right]+1$ times of map $F_\beta^2$, the orbit will pass the discontinuous line $\ell_2$, and we have 
\begin{equation}
	\begin{split}
		&\rho_{n_1}=\rho+\mathcal{O}\left(\frac{1}{\rho}\right),\\
		&\psi_{n_1}=\frac{3}{4\rho}-\frac{v_1}{\rho}+\mathcal{O}\left(\frac{1}{\rho^2}\right),
	\end{split}
\end{equation}
where $v_1=\left\{A_\beta\rho-\phi\right\}$. Changing back to $(r,\theta)$ coordinates, using the value from appendix, we have:
\begin{equation}
	\begin{split}
		&r_{n_1}=\rho+\mathcal{O}\left(\frac{1}{\rho}\right),\\
		&\theta_{n_1}=\beta+\frac{3}{\rho}-\frac{4v_1}{\rho}+\mathcal{O}\left(\frac{1}{\rho^2}\right),
	\end{split}
\end{equation}

Then we repeat this process for the region $\Rfive$, in which $r$ is conserved and $\theta$ grows by almost $4/r$ each time when we apply $F_\beta^2$. From proposition \ref*{prop3}, we have
$$\ell_2\subset \left\{\theta=\beta-\frac{1}{r}+\mathcal{O}\left(\frac{1}{r^2}\right)\right\},$$
and after $n_2=\left[\frac{(\pi-2\beta)\rho}{4}-\frac{3}{2}+v_1\right]+1$  time of map $F_\beta^2$, the orbit will go into region $\Rtwo$, then
\begin{equation}
	\begin{split}
		&r_{n_1+n_2}=\rho+\mathcal{O}\left(\frac{1}{\rho}\right),\\
		&\theta_{n_1+n_2}=\pi-\beta+\frac{1}{\rho}-\frac{4v_2}{\rho}+\mathcal{O}\left(\frac{1}{\rho^2}\right),
	\end{split}
\end{equation}
where $v_2=\left\{\frac{(\pi-2\beta)\rho}{4}-\frac{3}{2}+v_1\right\}$. Similarly, we change back to coordinates $(\rho,\beta)$ in region $\Rtwo$ and obtain
\begin{equation}
	\begin{split}
		&\rho_{n_1+n_2}=\rho+\mathcal{O}\left(\frac{1}{\rho}\right),\\
		&\psi_{n_1+n_2}=\frac{1}{4\rho}-\frac{v_2}{\rho}+\mathcal{O}\left(\frac{1}{\rho^2}\right).
	\end{split}
\end{equation}
Use (\ref*{region2}),
$$\ell_1' \subset\left\{\psi=A_\beta-\frac{3}{4\rho}+\mathcal{O}\left(\frac{1}{\rho^2}\right)\right\},$$
then we get after $n_3=\left[A_\beta\rho-1+v_2\right]+1$ times of map $F_\beta^2$, the orbit will pass $\ell_1'$, and 
\begin{equation}
	\begin{split}
		\rho_{n_1+n_+n_3}&=\rho+\mathcal{O}\left(\frac{1}{\rho}\right),\\
		\psi_{n_1+n_2+n_3}&=A_\beta+\frac{1}{4\rho}-\frac{v_3}{\rho}+\mathcal{O}\left(\frac{1}{\rho^2}\right),\\
		&=A_\beta-\frac{3}{4\rho}+\frac{v_4}{\rho}+\mathcal{O}\left(\frac{1}{\rho^2}\right),
	\end{split}
\end{equation}
where $v_3=\left\{A_\beta\rho-1+v_2\right\}$ and $v_4=1-v_3$. This time, we should change to coordinates $(r,\theta)$, shift by $\pi$ to obtain polar coordinates in the lower region, and change back to coordinates $(\rho,\psi)$ in the region $\Rthree$. Firstly, 
\begin{equation}
	\begin{split}
		&r_{n_1+n_2+n_3}=\frac{2}{\cos\beta+1}\rho+2\left(v_4-\frac{1}{2}\right)\sin\beta+\mathcal{O}\left(\frac{1}{\rho}\right),\\
		&\theta_{n_1+n_2+n_3}=\pi+\frac{\cos\beta+1}{\rho}\left[\left(v_4-\frac{1}{2}\right)\cos\beta+v_4-1\right]+\mathcal{O}\left(\frac{1}{\rho^2}\right).	\end{split}
\end{equation}
Shift by $\pi$, we have
\begin{equation}
	\theta_{n_1+n_2+n_3}=\frac{\cos\beta+1}{\rho}\left[\left(v_4-\frac{1}{2}\right)\cos\beta+v_4-1\right]+\mathcal{O}\left(\frac{1}{\rho^2}\right).
\end{equation}
Finally, change back to $(\rho,\psi)$:
\begin{equation}
	\begin{split}
		&\widetilde{\rho}=\rho+\left(v_4-\frac{1}{2}\right)\left(2\sin\beta+\sin2\beta\right)+\mathcal{O}\left(\frac{1}{\rho}\right),\\
		&\psi'=-A_\beta-\frac{3}{4\rho}+\frac{v_4}{\rho}+\mathcal{O}\left(\frac{1}{\rho^2}\right).
	\end{split}
\end{equation}
Now the orbit goes into the area between $\ell_1'$ and $F_\beta^2(\ell_1')$, which is almost central-symmetric to region $\mathcal{D}$. let $\widetilde{\phi}=v_4$ and in $(\widetilde{\rho},\widetilde{\phi})$ coordinates we have
\begin{equation}
	\begin{split}
		&\widetilde{\rho}=\rho+\left(\widetilde{\phi}-\frac{1}{2}\right)\left(2\sin\beta+\sin2\beta\right)+\mathcal{O}\left(\frac{1}{\rho}\right),\\
		&\widetilde{\phi}=v_4+\mathcal{O}\left(\frac{1}{\rho}\right).
	\end{split}
\end{equation}
Observe that 
\begin{align*}
	&v_1=\left\{A_\beta\rho-\phi\right\},\\
	&v_2=\left\{\frac{(\pi-2\beta)\rho}{4}-\frac{3}{2}+v_1\right\},\\
	&v_3=\left\{A_\beta\rho-1+v_2\right\},\\
	&v_4=1-v_3,
\end{align*}
it yields 
\begin{align*}
	&v_1=\left\{A_\beta\rho-\phi\right\},\\
	&v_2=\left\{\frac{(\pi-2\beta)\rho}{4}+A_\beta\rho+\frac{1}{2}-\phi\right\},\\
	&v_3=\left\{\frac{(\pi-2\beta)\rho}{4}+2A_\beta\rho+\frac{1}{2}-\phi\right\},\\
	&v_4=\left\{\phi+\frac{1}{2}-\rho B_\beta\right\},
\end{align*}
where $B_\beta=\left((\pi-2\beta)/4+2A_\beta\right)$ is constant depending only on $\beta$.
In summary, we have $\mathcal{F}_1$:
\begin{equation}
	\begin{split}
		&\widetilde{\rho}=\rho+\left(\widetilde{\phi}-\frac{1}{2}\right)\left(2\sin\beta+\sin2\beta\right)+\mathcal{O}\left(\frac{1}{\rho}\right),\\
		&\widetilde{\phi}=\left\{\phi+\frac{1}{2}-\rho B_\beta\right\}+\mathcal{O}\left(\frac{1}{\rho}\right).
	\end{split}
\end{equation}
Likewise, repeat all calculations for point in region $\mathcal{\widetilde{D}}$ and we get $\mathcal{F}_2$:
\begin{equation}
	\begin{split}
		&\hat{\rho}=\widetilde{\rho}+\left(\phi-\frac{1}{2}\right)\left(2\sin\beta+\sin2\beta\right)+\mathcal{O}\left(\frac{1}{\widetilde{\rho}}\right),\\
		&\hat{\phi}=\left\{\widetilde{\phi}+\frac{1}{2}-\widetilde{\rho} B_\beta\right\}+\mathcal{O}\left(\frac{1}{\widetilde{\rho}}\right).
	\end{split}
\end{equation}

It should be pointed out that map $\mathcal{F}_1$ and $\mathcal{F}_2$ have the same form in coordinates $(\rho,\phi)$ and $(\widetilde{\rho},\widetilde{\phi})$ and if we let
\begin{align*}
	R=\rho\left(\frac{(\pi-2\beta)}{4}+2A_\beta\right)+\frac{1}{2}=\rho B_\beta+\frac{1}{2},\\ \widetilde{R}=\widetilde{\rho}\left(\frac{(\pi-2\beta)}{4}+2A_\beta\right)+\frac{1}{2}=\widetilde{\rho} B_\beta+\frac{1}{2},
\end{align*}

\begin{align*}
	C_\beta&=(2\sin\beta+\sin2\beta)\left(\frac{(\pi-2\beta)}{4}+2A_\beta\right),\\
	&=(2\sin\beta+\sin2\beta)\left(\frac{(\pi-2\beta)}{4}+\frac{1}{3}\tan^3\left(\frac{\beta}{2}\right)+\tan\left(\frac{\beta}{2}\right)\right)
\end{align*}

we have $\mathcal{F}_\beta=\mathcal{L}^2+\mathcal{O}\left(\frac{1}{R^2}\right)$ and $\mathcal{L}$: 
\begin{align*}
	&R\mapsto R+C_\beta\left(\left\{\phi-R\right\}-\frac{1}{2}\right),\\
	&\phi\mapsto \left\{\phi-R\right\}.
\end{align*}

The singularities are determined by preimage of discontinuous lines of $\ell_2$, $\ell_3'$ and $\ell_1'$ which are corresponding to ${v_1=0}$, ${v_2=0}$ and ${v_3=0}$, that is:
\begin{align*}
	&0=v_1=\left\{A_\beta\rho-\phi\right\}+\mathcal{O}\left(\frac{1}{R}\right),\\
	&0=v_2=\left\{\frac{(\pi-2\beta)\rho}{4}+A_\beta\rho+\frac{1}{2}-\phi\right\}+\mathcal{O}\left(\frac{1}{R}\right),\\
	&0=v_3=\left\{\frac{(\pi-2\beta)\rho}{4}+2A_\beta\rho+\frac{1}{2}-\phi\right\}+\mathcal{O}\left(\frac{1}{R}\right).
\end{align*}
And if we let
\begin{align*}
	C_1=\frac{A_\beta}{2A_\beta+\frac{\pi-2\beta}{4}},\ C_2=\frac{A_\beta+\frac{\pi-2\beta}{4}}{2A_\beta+\frac{\pi-2\beta}{4}},\ 
\end{align*}
the singularities can be written as
\begin{equation}
	\begin{split}
		&\frac{1}{2}C_1=\left\{C_1R-\phi\right\}+\mathcal{O}\left(\frac{1}{R}\right),\\
		&1-\frac{1}{2}C_2=\left\{C_2R-\phi\right\}+\mathcal{O}\left(\frac{1}{R}\right),\\
		&0=\left\{R-\phi\right\}+\mathcal{O}\left(\frac{1}{R}\right)\\
	\end{split}
\end{equation}
\begin{remark}
	It should be pointed out that when $\beta=\pi/2$, we have $C_\beta=8/3$ and $C_1=C_2=1/2$, then our piecewise linear map is the same as in \cite{dolgopyat_unbounded_2009}.
\end{remark}

\section{Calculation results}\label{important values}
We give the values of $\Phi_i$ and $\Psi_i$ and their derivatives near lines of discontinuity here.
   \begin{table}[H]
\setlength{\arraycolsep}{4pt}
   \[\begin{array}{ccccccccccc}\toprule
	&	\Phi_1 & \Phi_1' &\Phi_1'' &\Phi_1''' &\Phi_2&\Phi_2'&\Phi_2''&\Phi_3&\Phi_3'&\Phi_4\\ \midrule
	{\rm region}\ \Rone\ {\rm at}\ 0&1 & 1 & 0 & -1 & 0&0&1&0&-\frac{1}{6}&0\\ 
	{\rm region}\ \Rone\ {\rm at}\ \pi/2&2 & 0 & -1 & 0 & 1&1&0&\frac{1}{12}&0&e_1\\
	{\rm region}\ \Rtwo\ {\rm at}\ \pi&1 & -1 & 0 & 1 & 0&0&1&0&\frac{1}{6}&0\\
	{\rm region}\ \Rtwo\ {\rm at}\ \pi/2&2 & 0 & -1 & 0& 1&-1&0&-\frac{11}{12}&0&e_2\\
	{\rm region}\ \Rthree\ {\rm at}\ 0&1 & 1 & 0 & -1 & 0&0&-1&0&-\frac{1}{6}&0\\ 
	{\rm region}\ \Rthree\ {\rm at}\ \pi/2&2 & 0 & -1 & 0 & -1&-1&0&-\frac{11}{12}&0&e_3\\ 
	{\rm region}\ \Rfour\ {\rm at}\ \pi& 1 & -1 & 0 & 1 & 0&0&-1&0&\frac{1}{6}&0\\ 
	{\rm region}\ \Rfour\ {\rm at}\ \pi/2& 2 & 0 & -1 & 0 & -1 &1&0&\frac{1}{12}&0&e_4\\ \bottomrule
\end{array}\]
\caption{values of $\Phi_i$ and their derivatives near discontinuity lines} 
\label{values of Phi}
\end{table}
\begin{table}[H]
\setlength{\arraycolsep}{4pt}
    \[\begin{array}{ccccccccccc}\toprule
	&	\Psi & \Psi' &\Psi'' &\Psi''' &\Psi_1&\Psi_1'&\Psi_1''&\Psi_2&\Psi_2'&\Psi_3\\ \midrule
	{\rm region}\ \Rone\ {\rm at}\ 0&0 & \frac{1}{2} & -1 & 3 & 0&0&0&0&\frac{1}{12}&0\\ 
	{\rm region}\ \Rone\ {\rm at}\ \pi/2&\frac{1}{3} & \frac{1}{8}& 0 & \frac{1}{8} & -\frac{1}{24}&-\frac{1}{16}&-\frac{1}{12}&\frac{5}{144}&\frac{11}{192}&f_1\\
	{\rm region}\ \Rtwo\ {\rm at}\ \pi&\frac{2}{3} & \frac{1}{2} & 1 & 3 & 0&0&0&0&\frac{1}{12}&0\\
	{\rm region}\ \Rtwo\ {\rm at}\ \pi/2&\frac{1}{3} & \frac{1}{8} & 0 & \frac{1}{8}& \frac{1}{24}&-\frac{1}{16}&\frac{1}{12}&\frac{\pi }{8}-\frac{29}{144}&-\frac{25}{192}&f_2\\
	{\rm region}\ \Rthree\ {\rm at}\ 0&0 & \frac{1}{2} & -1 & 3 & 0&0&0&0&\frac{1}{12}&0\\ 
	{\rm region}\ \Rthree\ {\rm at}\ \pi/2&\frac{1}{3} & \frac{1}{8} & 0 & \frac{1}{8} & \frac{1}{24}&\frac{1}{16}&\frac{1}{12}&-\frac{\pi }{8}+\frac{29}{144}&-\frac{25}{192}&f_3\\ 
	{\rm region}\ \Rfour\ {\rm at}\ \pi& \frac{2}{3} & \frac{1}{2} & 1 & 3 & 0&0&0&0&\frac{1}{12}&0\\ 
	{\rm region}\ \Rfour\ {\rm at}\ \pi/2& \frac{1}{3} & \frac{1}{8} & 0 & \frac{1}{8} & -\frac{1}{24} &\frac{1}{16}&-\frac{1}{12}&-\frac{5}{144}&\frac{11}{192}&f_4\\ \bottomrule
\end{array}\]
\caption{values of $\Psi_i$ and their derivatives near discontinuity lines}
\label{values of Psi}
\end{table}

\bibliographystyle{plain}
\bibliography{reffordisk.bib}

\begin{thebibliography}{10}

\bibitem{BIRD1988164}
Nigel Bird and Franco Vivaldi.
\newblock Periodic orbits of the sawtooth maps.
\newblock {\em Physica D: Nonlinear Phenomena}, 30(1):164--176, 1988.

\bibitem{boyland1994dualbilliardstwistmaps}
Philip Boyland.
\newblock Dual billiards, twist maps, and impact oscillators, 1994.

\bibitem{bullett_invariant_1986}
Shaun Bullett.
\newblock Invariant circles for the piecewise linear standard map.
\newblock {\em Communications in Mathematical Physics}, 107(2):241--262, 1986.

\bibitem{PhysRevA.24.2664}
John~R. Cary and James~D. Meiss.
\newblock Rigorously diffusive deterministic map.
\newblock {\em Phys. Rev. A}, 24:2664--2668, 1981.

\bibitem{CHEN1990217}
Q.~Chen, I.~Dana, J.D. Meiss, N.W. Murray, and I.C. Percival.
\newblock Resonances and transport in the sawtooth map.
\newblock {\em Physica D: Nonlinear Phenomena}, 46(2):217--240, 1990.

\bibitem{QChen_1989}
Q~Chen and J~D Meiss.
\newblock Flux, resonances and the devil's staircase for the sawtooth map.
\newblock {\em Nonlinearity}, 2(2):347, 1989.

\bibitem{de2012dynamics}
Jacopo De~Simoi and Dmitry Dolgopyat.
\newblock Dynamics of some piecewise smooth fermi-ulam models.
\newblock {\em Chaos: An Interdisciplinary Journal of Nonlinear Science}, 22(2), 2012.

\bibitem{de_Simoi_2012}
Jacopo de~Simoi and Dmitry Dolgopyat.
\newblock Dynamics of some piecewise smooth fermi-ulam models.
\newblock {\em Chaos: An Interdisciplinary Journal of Nonlinear Science}, 22(2), June 2012.

\bibitem{dolgopyat_unbounded_2009}
Dmitry Dolgopyat and Bassam Fayad.
\newblock Unbounded orbits for semicircular outer billiard.
\newblock {\em Ann. Henri Poincar\'{e}}, 10(2):357--375, 2009.

\bibitem{Douady}
Raphael Douady.
\newblock Th\'{e}se de 3-\'{e}me cycle.
\newblock {\em Universit\'{e} de Paris 7}, 1982.

\bibitem{gutkin_dual_1992}
Eugene Gutkin and Nandor Simanyi.
\newblock Dual polygonal billiards and necklace dynamics.
\newblock {\em Communications in Mathematical Physics}, 143(3):431--449, 1992.

\bibitem{mathematikExteriorBilliards}
Oliver Knill.
\newblock {E}xterior {B}illiards --- mathematik.com.
\newblock \url{https://www.mathematik.com/Exterior/index.html}.

\bibitem{MR0442980}
J\"{u}rgen Moser.
\newblock {\em Stable and random motions in dynamical systems}, volume No. 77 of {\em Annals of Mathematics Studies}.
\newblock Princeton University Press, Princeton, NJ; University of Tokyo Press, Tokyo, 1973.

\bibitem{Moser+2001}
Jurgen Moser.
\newblock {\em Stable and Random Motions in Dynamical Systems}.
\newblock Princeton University Press, Princeton, 2001.

\bibitem{moser_invariant_1962}
Jürgen Moser.
\newblock On invariant curves of area-preserving mappings of an annulus.
\newblock {\em Nachrichten der Akademie der Wissenschaften in Göttingen. {II}. Mathematisch-Physikalische Klasse}, 1962:1--20, 1962.

\bibitem{PERCIVAL1987373}
Ian Percival and Franco Vivaldi.
\newblock A linear code for the sawtooth and cat maps.
\newblock {\em Physica D: Nonlinear Phenomena}, 27(3):373--386, 1987.

\bibitem{schwartz2007unboundedorbitsouterbilliards}
Richard~Evan Schwartz.
\newblock Unbounded orbits for outer billiards, 2007.

\bibitem{schwartz2008outerbilliardskites}
Richard~Evan Schwartz.
\newblock Outer billiards on kites, 2008.

\bibitem{vivaldi_global_1987}
Franco Vivaldi and Anna~V. Shaidenko.
\newblock Global stability of a class of discontinuous dual billiards.
\newblock {\em Communications in Mathematical Physics}, 110(4):625--640, 1987.

\bibitem{Wojtkowski_1982}
M.~Wojtkowski.
\newblock On the ergodic properties of piecewise linear perturbations of the twist map.
\newblock {\em Ergodic Theory and Dynamical Systems}, 2(3–4):525–542, 1982.

\end{thebibliography}

\end{document}